\documentclass[10pt, twoside, a4paper]{article}

\usepackage[T1]{fontenc}

\usepackage{latexsym}
\usepackage{amsmath}
\usepackage{amssymb}
\usepackage{amsfonts}
\usepackage{amsthm}
\usepackage{amscd}
\usepackage{epsfig}
\usepackage{picins}
\usepackage{verbatim}
\usepackage{fancybox}
\usepackage{moreverb}
\usepackage{psfrag}

\newcommand{\R}{{\mathbb R}}
\newcommand{\N}{{\mathbb N}}

\newcommand{\Q}{{\mathbb Q}}

\newcommand{\D}{{\partial}}

\newtheorem{theorem}{Theorem}

\begin{document}

\begin{center}
\Large{\textbf{Two-scale semi-lagrangian simulation of a charged particle beam in a periodic focusing channel}} \\
\end{center}
\begin{center}
\textit{Alexandre Mouton} \\
IRMA, Universit\'e Louis Pasteur, F-67084 Strasbourg \\ \textit{} \\ \textit{} \\ \textit{} \\
\end{center}

\begin{abstract}
This paper is devoted to numerical simulation of a charged particle beam submitted to a strong oscillating electric field. For that, we consider a two-scale numerical approach as follows: we first recall the two-scale model which is obtained by using two-scale convergence techniques; then, we numerically solve this limit model by using a backward semi-lagrangian method and we propose a new mesh of the phase space which allows us to simplify the solution of the Poisson's equation. Finally, we present some numerical results which have been obtained by the new method, and we validate its efficiency through long time simulations.
\end{abstract}

\paragraph*{AMS subjects:} 35B27, 76X05, 65M25 (65Y20).


\section{Introduction}
\setcounter{equation}{0}

Recent papers have proved that the two-scale convergence theory developed by Allaire in \cite{Allaire} and Nguetseng in \cite{NGuetseng} can be used successfully in order to develop numerical methods for solving ODEs or PDEs with oscillatory singular perturbations: for example, Fr\'enod, Salvarani and Sonnendr\"ucker have developed a two-scale PIC method in \cite{PIC-two-scale} for simulations of charged particle beams in a periodic focusing channel, and Fr\'enod, Mouton and Sonnendr\"ucker have developed a two-scale finite volume method in \cite{EVA01} for solving the weakly compressible 1D Euler equations. We can also cite the work of Ailliot, Fr\'enod and Monbet in \cite{Ailliot-Frenod-Monbet} about the simulation of tide oscillation for long term drift forecast of objects in the ocean. \\
\indent On the other hand, many papers devoted to the numerical simulation of Vlasov-type problems involve Particle-In-Cell methods (see Birdsall and Langdon \cite{Birdsall}) or Eulerian methods like semi-lagrangian schemes (see Sonnendr\"ucker \textit{et al.}\cite{Bertrand-Ghizzo-Roche-Sonnen}, Grandgirard \textit{et al.}\cite{GYSELA4D}, Filbet and Sonnendr\"ucker \cite{Comparison,Paraxial}, Cheng and Knorr \cite{Cheng-Knorr}). Since papers like \cite{EVA01} or \cite{PIC-two-scale} can be viewed as parts of a work programme which goal is the development of two-scale numerical methods for simulations of magnetic confinement fusion, it can be interesting to couple two-scale convergence results on Vlasov models such as two-scale models obtained in \cite{Larmor-radius} with a semi-lagrangian scheme. \\
\indent However, it is preferable in a first time to develop such a method on a more simple problem in order to study its behavior especially in the context of non-smooth solutions. The context of charged particle beams in a periodic focusing external field described in \cite{PIC-two-scale} offers a relatively simple framework for answering these questions. \\
\indent We recall that such a phenomenon can be successfully represented by the 3D Vlasov-Maxwell system. In the same spirit of \cite{PIC-two-scale}, we will consider non-relativistic long and thin beams, so we can replace the full three-dimensional Vlasov-Maxwell system by its paraxial approximation. To obtain this approximation, we do the following assumptions:
\begin{itemize}
\item the beam has already reached its stationary state,
\item the beam is long and thin, 
\item the beam is propagating at constant velocity $v_{z}$ along the longitudinal axis $z$,
\item the beam is sufficiently long so we can neglect the longitudinal self-consistent forces,
\item the external electric field is supposed to be $l$-periodic in $z$ and independent of the time,
\item the beam is axisymmetric,
\item the initial distribution $f_{0}$ is concentrated in angular momentum.
\end{itemize}
Under all these assumptions, the 3D Vlasov-Maxwell system reduces itself to a 2D Vlasov-Poisson system where the longitudinal position coordinate $z$ plays the role of time, and even to a 1D axisymmetric Vlasov-Poisson system of the form
\begin{equation}\label{NH-polar}
\begin{split}
& \D_{t}f^{\epsilon} + \frac{v_{r}}{\epsilon} \, \D_{r}f^{\epsilon} + \big(E_{r}^{\epsilon} + \Xi_{r}^{\epsilon}\big) \, \D_{v_{r}}f^{\epsilon} = 0 \, , \\
& f^{\epsilon}(r,v_{r},t=0) = f_{0}(r,v_{r}) \, , \\
& \frac{1}{r} \, \D_{r}(r\,E_{r}^{\epsilon}) = \int_{\R} f^{\epsilon} \, dv_{r} \, , \\
& \Xi_{r}^{\epsilon}(r,t) = -\frac{1}{\epsilon} \, H_{0} \, r + H_{1}\Big(\omega_{1} \, \frac{t}{\epsilon} \Big) \, r \, ,
\end{split}
\end{equation}
where $r \geq 0$ is the radial component of the position vector in the transverse plane to the propagation direction, $v_{r} \in \R$ is the projection of the transverse velocity in the transverse plan to the propagation direction, $\epsilon$ is the ratio between the characteristic transverse radius of the beam and the characteristic longitudinal length of the beam, $f^{\epsilon} = f^{\epsilon}(r,v_{r},t)$ is the distribution function of the particles, $E_{r}^{\epsilon} = E_{r}^{\epsilon}(r,t)$ is the radial part of the transverse self-consistent electric field, $\Xi_{r}^{\epsilon} = \Xi_{r}^{\epsilon}(r,t)$ is the radial part of the transverse external electric field defined with the dimensionless real constant $\omega_{1}$ and the tension functions $H_{0}$ and $H_{1}$ which are respectively constant and periodic. All these quantities and variables are dimensionless. This system is naturally defined for $r \geq 0$ but we can extend it to $r \in \R$ by using the conventions $f^{\epsilon}(-r,-v_{r},t) = f^{\epsilon}(r,v_{r},t)$, $E_{r}^{\epsilon}(-r,t) = -E_{r}^{\epsilon}(r,t)$, and $\Xi_{r}^{\epsilon}(-r,t) = -\Xi_{r}^{\epsilon}(r,t)$. \\

\indent The aim the paper is to simulate the system (\ref{NH-polar}) with a two-scale semi-lagrangian scheme when $\epsilon \to 0$. Inspired by Fr\'enod, Salvarani and Sonnendr\"ucker \cite{PIC-two-scale}, and Fr\'enod, Mouton and Sonnendr\"ucker \cite{EVA01}, we derive the model (\ref{NH-polar}) by using the two-scale convergence theory developed in Allaire \cite{Allaire} and Nguetseng \cite{NGuetseng} and we obtain a new model which is independent of $\epsilon$. Then, instead of discretizing the model (\ref{NH-polar}) with a classical semi-lagrangian scheme, we discretize the new model with a semi-lagrangian method in order to obtain an approximation of a function $F = F(r,v_{r},\tau,t)$, where $\tau$ is the second time scale, which verifies $f^{\epsilon}(r,v_{r},t) \sim 2 \pi F(r,v_{r},\frac{t}{\epsilon},t)$. Proceeding in such a way presents the advantage that there is no longer $\frac{1}{\epsilon}$-oscillations in the limit model, so a very small time step is no longer required in order to simulate these oscillations. However, since semi-lagrangian schemes are based on interpolation on a phase space mesh, we have to pay attention to the contributions of the second time scale $\tau$ in the limit model. For this reason, in the same spirit of the work of Lang \textit{et al.}\cite{Lang}, we introduce a $\tau$-dependent moving mesh adapted for this two-scale numerical method the goal of which is to reduce the number of interpolations during the simulation. \\

\indent The paper is organized as follows: in section 2, we recall the procedure leading from the paraxial approximation of the complete 3D Vlasov-Maxwell system to the model (\ref{NH-polar}), and we will recall some two-scale convergence results about this model. In section 3, we build the semi-lagrangian method on the limit model and we see how to simplify it by considering a particular mesh. Section 4 is devoted to numerical results obtained with the two-scale semi-lagrangian method: on one hand, we compare them to some results obtained from a classical semi-lagrangian method on the system (\ref{NH-polar}) in terms of quality of results and CPU time, and on the other hand, we see in this section the consequences of the use of the new mesh in the same terms.

\section{The two-scale model}
\setcounter{equation}{0}

Firstly, we recall in this section the way to obtain the system (\ref{NH-polar}) from the paraxial approximation of the 3D Vlasov-Maxwell equations. Then, we present a theorem about the two-scale convergence of the solutions of (\ref{NH-polar}) which has been proved by Fr\'enod, Salvarani and Sonnendr\"ucker in \cite{PIC-two-scale}. \\

\subsection{Scaling of the paraxial model}

\indent By applying the hypothesis about the considered beam which are mentioned in the introduction, the full three-dimensional Vlasov-Maxwell is reduced to
\begin{equation} \label{VP-dimension}
\begin{split}
& v_{z} \, \D_{z}f + \mathbf{v} \cdot \nabla_{\mathbf{x}} f + \frac{q}{m} \, (\mathbf{E} + \mathbf{\Xi}) \cdot \nabla_{\mathbf{v}} f = 0 \, , \\
& f(\mathbf{x},z = 0, \mathbf{v}) = f_{0}(\mathbf{x},\mathbf{v}) \, , \\
& - \nabla_{\mathbf{x}} \phi = \mathbf{E} \, , \qquad - \Delta_{\mathbf{x}} \phi = \frac{\rho}{\varepsilon_{0}} = \frac{q}{\varepsilon_{0}} \int_{\R^{2}} f \, d\mathbf{v} \, , 
\end{split}
\end{equation}
where $z$ is the longitudinal position coordinate, $\mathbf{x} = (x,y)$ is the transverse position vector, $v_{z}$ is the constant longitudinal speed, $\mathbf{v} = (v_{x},v_{y})$ is the transverse speed vector, $f = f(\mathbf{x},z,\mathbf{v})$ is the distribution function of particles whose charge is $q$ and mass is $m$, $\phi = \phi(\mathbf{x},z)$ is the potential linked with the transverse self-consistent electric field $\mathbf{E} = \mathbf{E}(\mathbf{x},z)$, and $\mathbf{\Xi} = \mathbf{\Xi}(\mathbf{x},z)$ is the transverse external electric field. More details about this derivation can be found in \cite{Degond-Raviart} and \cite{Paraxial}. In this model, we assume that the external electric field is given by
\begin{equation} \label{Xi-dimension}
\mathbf{\Xi}(\mathbf{x},z) = - H_{0} \, \mathbf{x} + H_{1} \Big( \omega_{1} \, \frac{z}{l} \Big) \, \mathbf{x} \, ,
\end{equation}
where $H_{0}$ is a positive constant tension, $H_{1}$ is a $l$-periodic tension and $\omega_{1}$ is a dimensionless real constant. \\

\indent In the same spirit of Fr\'enod, Salvarani and Sonnendr\"ucker \cite{PIC-two-scale}, we build a dimensionless version of the system (\ref{VP-dimension})-(\ref{Xi-dimension}) by introducing the dimensionless variables $\mathbf{x}',z',\mathbf{v}'$ defined by
\begin{equation}
\mathbf{x} = \lambda \, \mathbf{x}' \, , \qquad z = L \, z' \, , \qquad \mathbf{v} = v_{z} \, \mathbf{v}' \, ,
\end{equation}
where $\lambda$ is the characteristic transverse radius of the beam and $L$ is the characteristic length of the beam. Moreover, we define the dimensionless quantities $f',f_{0}',\mathbf{E}',\phi', H_{0}',H_{1}'$ by
\begin{equation}
\begin{array}{rclrcl}
f(\lambda\,\mathbf{x}',L\,z',v_{z}\,\mathbf{v}') &=& \displaystyle \frac{m\,\varepsilon_{0}}{q^{2} \,\lambda \, L} \, f'(\mathbf{x}',z',\mathbf{v}') \, , & \qquad \mathbf{E}(\lambda \, \mathbf{x}', L \, z') &=& \displaystyle \frac{m \, v_{z}^{2}}{q \, L} \, \mathbf{E}'(\mathbf{x}',z') \, , \\ \\
f_{0}(\lambda\,\mathbf{x}',v_{z}\,\mathbf{v}') &=& \displaystyle \frac{m\,\varepsilon_{0}}{q^{2} \,\lambda \, L} \, f_{0}'(\mathbf{x}',\mathbf{v}') \, , & H_{0} &=& \displaystyle \overline{H_{0}}  \, H_{0}' \, , \\ \\
\phi(\lambda \, \mathbf{x}', L \, z') &=& \displaystyle \frac{m \, \lambda \, v_{z}^{2}}{q \, L} \, \phi'(\mathbf{x}',z') \, , & H_{1}(\tau) &=& \displaystyle \overline{H_{1}} \, H_{1}'(\tau) \, .
\end{array}
\end{equation}
With these new variables, the system (\ref{VP-dimension})-(\ref{Xi-dimension}) can be rewritten under the form
\begin{equation}\label{VP-no_dimension}
\begin{split}
&\D_{z'}f' + \frac{L}{\lambda} \, \mathbf{v}'\cdot \nabla_{\mathbf{x}'}f' + \Bigg[ \mathbf{E}' - \frac{\overline{H_{0}} \, \lambda \, q \, L}{v_{z}^{2} \, m} \, H_{0}'\,\mathbf{x}' + \frac{\overline{H_{1}} \, \lambda \, q \, L}{v_{z}^{2} \, m} \, H_{1}'\Big(\frac{\omega_{1} \,L\,z'}{l}\Big) \, \mathbf{x}' \Bigg] \cdot \nabla_{\mathbf{v}'}f' = 0 \, , \\
&f'(\mathbf{x}',z' = 0, \mathbf{v}') = f_{0}'(\mathbf{x}',\mathbf{v}') \, , \\
&-\nabla_{\mathbf{x}'}\phi' = \mathbf{E}' \, , \qquad -\Delta_{\mathbf{x}'}\phi' = \int_{\R^{2}} f' \, d\mathbf{v}' \, .
\end{split}
\end{equation}
Since the beam is supposed to be long and thin, it is natural to take the ratio
\begin{equation}
\frac{\lambda}{L} = \epsilon \, . 
\end{equation}
Furthermore, as we want to simulate the beam over a large number of periods of the external electric field, we also consider the ratio
\begin{equation}
\frac{l}{L} = \epsilon \, .
\end{equation}
Finally, we suppose that the external electric field is much stronger than the self-consistent electric field and that its oscillations in $z$ direction are of the same order as $\mathbf{E}$, so we consider that
\begin{equation}
\frac{\overline{H_{0}} \, \lambda \, q \, L}{v_{z}^{2} \, m} = \frac{1}{\epsilon} \, , \qquad \frac{\overline{H_{1}} \, \lambda \, q \, L}{v_{z}^{2} \, m} = 1 \, .
\end{equation}
Under all these hypothesis the system (\ref{VP-no_dimension}) reduces to
\begin{equation} \label{NH-cartesian}
\begin{split}
&\D_{t}f^{\epsilon} + \frac{1}{\epsilon} \, \mathbf{v}\cdot \nabla_{\mathbf{x}}f^{\epsilon} + (\mathbf{E}^{\epsilon} + \mathbf{\Xi}^{\epsilon}) \cdot \nabla_{\mathbf{v}}f^{\epsilon} = 0 \, , \\
&f^{\epsilon}(\mathbf{x}, \mathbf{v}, t = 0) = f_{0}(\mathbf{x},\mathbf{v}) \, , \\
&-\nabla_{\mathbf{x}}\phi^{\epsilon} = \mathbf{E}^{\epsilon} \, , \qquad -\Delta_{\mathbf{x}}\phi^{\epsilon} = \int_{\R^{2}} f^{\epsilon} \, d\mathbf{v} \, , \\
&\mathbf{\Xi}^{\epsilon}(\mathbf{x},t) = -\frac{1}{\epsilon} \, H_{0} \, \mathbf{x} + H_{1} \Big( \frac{\omega_{1} \, t}{\epsilon} \Big) \, \mathbf{x} \, ,
\end{split}
\end{equation}
where the primed notations for the variables and the initial distribution have been eliminated, and where $z'$ has been replaced by $t$, $f'$ by $f^{\epsilon}$, $\mathbf{E}'$ by $\mathbf{E}^{\epsilon}$, $\phi'$ by $\phi^{\epsilon}$, and $\mathbf{\Xi}'$ by $\mathbf{\Xi}^{\epsilon}$. \\
\indent We introduce the polar coordinates $(r,\theta,v_{r},v_{\theta})$ linked with $(\mathbf{x},\mathbf{v})$ by the relations
\begin{equation}
\begin{array}{rclrcl}
x &=& r \, \cos\theta \, , & \qquad v_{r} &=& v_{x} \, \cos\theta + v_{y} \, \sin\theta \, ,\\
y &=& r \, \sin\theta \, , & \qquad v_{\theta} &=& v_{y} \, \cos\theta - v_{x} \, \sin\theta \, . 
\end{array}
\end{equation}
Since the beam is supposed to be axisymmetric, the system does not depend on $\theta$. Furthermore, we assume that the initial distribution is concentrated in angular momentum. Then the system (\ref{NH-cartesian}) reduces to (\ref{NH-polar}).

\subsection{Two-scale convergence results}

\indent Since the aim of the paper is to develop a two-scale numerical method in order to simulate the model (\ref{NH-polar}), we have to establish that the solution $f^{\epsilon}$ of this model two-scale converges to a function $F = F(r,v_{r},\tau,t)$ in a certain Banach space, and we have to find a system of equations verified by $F$. These results have been proved in \cite{PIC-two-scale} and are recalled in the theorem below.

\begin{theorem}
We assume that $H_{0} = 1$ and that the initial distribution $f_{0}$ of (\ref{NH-polar}) satisfies the following properties:
\begin{itemize}
\item[(i)] $f_{0} \in L^{1}\big(\R^{2};\,|r|drdv_{r}\big) \cap L^{p}\big(\R^{2};\,|r|drdv_{r}\big)$ with $p \geq 2$,
\item[(ii)] $f_{0}(r,v_{r}) \geq 0$ for all $(r,v_{r}) \in \R^{2}$,
\item[(iii)] $\displaystyle \int_{\R^{2}} \big( r^{2} + v_{r}^{2} \big) \, f_{0}(r,v_{r}) \, |r|\, dr \, dv_{r} < + \infty$ .
\end{itemize}
Then, considering a sequence of solutions $(f^{\epsilon},\mathbf{E}_{r}^{\epsilon})$ of (\ref{NH-polar}) and extracting a subsequence from it, we can say that $f^{\epsilon}$ two-scale converges to $F \in L^{\infty}\big([0,T] \times [0,2\pi]; L^{2}(\R^{2};\,|r|drdv_{r})\big)$ and $\mathbf{E}_{r}^{\epsilon}$ two-scale converges to $\mathcal{E}_{r} \in L^{\infty}\big([0,T] \times [0,2\pi]; W^{1,\,3/2}(\R;\,|r|dr)\big)$. Furthermore, there exists a function $G \in L^{\infty}\big([0,T]; L^{2}(\R^{2};\,|q|dqdu_{r})\big)$ such that
\begin{equation}\label{FG-polar}
F(r,v_{r},\tau,t) = G\big(\cos(\tau)\,r - \sin(\tau)\,v_{r}, \sin(\tau)\,r + \cos(\tau)\,v_{r}, t\big) \, ,
\end{equation}
and $(G,\mathcal{E}_{r})$ is solution of
\begin{equation}\label{H-polar}
\begin{split}
&\D_{t}G - \Bigg[ \int_{0}^{2\pi} \sin(\sigma) \mathcal{E}_{r}\big(\cos(\sigma) \, q + \sin(\sigma) \, u_{r}, \sigma, t) \, d\sigma \Bigg] \, \D_{q}G \\
&\qquad \qquad - \Bigg[ \int_{0}^{2\pi} \sin(\sigma) \, \frac{I_{\Q}(\omega_{1})}{2\pi} \, H_{1}(\omega_{1} \, \sigma) \, \big(\cos(\sigma) \, q + \sin(\sigma) \, u_{r}\big) \, d\sigma \Bigg] \, \D_{q}G \\
&\qquad \qquad + \Bigg[ \int_{0}^{2\pi} \cos(\sigma) \mathcal{E}_{r}\big(\cos(\sigma) \, q + \sin(\sigma) \, u_{r}, \sigma, t) \, d\sigma \Bigg] \, \D_{u_{r}}G \\
&\qquad \qquad + \Bigg[ \int_{0}^{2\pi} \cos(\sigma) \, \frac{I_{\Q}(\omega_{1})}{2\pi} \, H_{1}(\omega_{1} \, \sigma) \, \big(\cos(\sigma) \, q + \sin(\sigma) \, u_{r}\big) \, d\sigma \Bigg] \, \D_{u_{r}}G = 0 \, , \\
& G(q,u_{r},t = 0) = \frac{1}{2\pi} \, f_{0}(q,u_{r}) \, , \\
& \frac{1}{r} \, \D_{r} \big( r \, \mathcal{E}_{r}(r,\tau,t) \big) = \int_{\R} G\big(\cos(\tau)\,r - \sin(\tau)\,v_{r}, \sin(\tau)\,r + \cos(\tau)\,v_{r}, t\big) \, dv_{r} \, ,
\end{split}
\end{equation}
where $I_{\Q}(\omega_{1})$ is equal to 1 if $\omega_{1} \in \Q$, and 0 otherwise.
\end{theorem}

Of course, such a result exists for the solution of the system (\ref{NH-cartesian}) and can be found in \cite{PIC-two-scale}.

\section{The two-scale semi-lagrangian method}
\setcounter{equation}{0}

\indent In this section, we develop a two-scale semi-lagrangian method in order to approach the solution $f^{\epsilon}$ of (\ref{NH-polar}), in the case where $H_{0} = 1$. As it has been suggested in the introduction, the strategy is to discretize the model (\ref{FG-polar})-(\ref{H-polar}) in order to obtain a good approximation of $F$ which can be used for approaching $f^{\epsilon}(r,v_{r},t) \sim 2\pi F(r,v_{r},\frac{t}{\epsilon},t)$. As in the two-scale PIC-method developed by Fr\'enod, Salvarani and Sonnendr\"ucker in \cite{PIC-two-scale}, there is an advantage by proceeding in such a way: since there is no longer $\frac{1}{\epsilon}$-frequency oscillations in the system (\ref{H-polar}), we do not need a very small time step for good simulation. In a first time, we recall the basis of a semi-lagrangian method. Then we present a motivation for the development of a two-scale semi-lagrangian method through the description of a classical backward semi-lagrangian method on the model (\ref{NH-polar}). Finally, we describe the two-scale numerical method itself and we suggest a new mesh in order to simplify it.

\subsection{The semi-lagrangian method}

\indent In this paragraph, we recall the way to discretize the abstract model
\begin{equation} \label{abstract-model}
\D_{t}f(\mathbf{x},t) + \mathbf{U}(\mathbf{x},t) \cdot \nabla_{\mathbf{x}}f(\mathbf{x},t) = 0
\end{equation}
with a semi-lagrangian method. For this, we have to consider the characteristics of (\ref{abstract-model}), which are the solutions of
\begin{equation} \label{abstract-charac}
\D_{t} \mathbf{X}(t) = \mathbf{U}\big(\mathbf{X}(t),t\big)  \, .
\end{equation}
It is an easy game to remark that $f$ is constant along the characteristics, i.e.
\begin{equation}
\D_{t} \big( f\big(\mathbf{X}(t),t\big)\big) = 0 \, ,
\end{equation}
so we can write
\begin{equation} \label{conservation-f-abstract}
f\big(\mathbf{X}(t;\mathbf{x},s),t\big) = f(\mathbf{x},s) \, ,
\end{equation}
where $\mathbf{X}(t; \mathbf{x},s)$ is the solution (\ref{abstract-charac}) with the condition $\mathbf{X}(s) = \mathbf{x}$. \\
\indent This property of $f$ is used in the semi-lagrangian method as follows: assuming that we know the value of $f$ at time $t_{n}-\Delta t$ on the mesh points $(\mathbf{x}_{i})_{i\,=\,0,\,\dots,\,N}$, we use the property (\ref{conservation-f-abstract}) to say that
\begin{equation}
f(\mathbf{x}_{i},t_{n}+\Delta t) = f\big(\mathbf{X}(t_{n}-\Delta t; \mathbf{x}_{i}, t_{n}+\Delta t), t_{n}-\Delta t\big) \, ,
\end{equation}
so we have to compute the point $\mathbf{X}(t_{n}-\Delta t; \mathbf{x}_{i}, t_{n}+\Delta t)$ first, then compute $f\big(\mathbf{X}(t_{n}-\Delta t; \mathbf{x}_{i}, t_{n}+\Delta t), t_{n}-\Delta t\big)$ by interpolating $f(\cdot,t_{n}-\Delta t)$ on the points $(\mathbf{x}_{i})_{i\,=\,0,\,\dots,\,N}$ in order to obtain an approximation of $f(\mathbf{x}_{i},t_{n}+\Delta t)$. Sonnendr\"ucker \textit{et al.} have suggested in \cite{Bertrand-Ghizzo-Roche-Sonnen} a way to compute a second order approximation of $\mathbf{X}(t_{n}-\Delta t; \mathbf{x}_{i}, t_{n}+\Delta t)$: they discretize the equation (\ref{abstract-charac}) with a finite difference method in order to obtain the following approximation:
\begin{equation} \label{abstract-fix-pt}
\mathbf{X}(t_{n}-\Delta t; \mathbf{x}_{i}, t_{n}+\Delta t) = \mathbf{x}_{i} - 2 \,\mathbf{d}_{i} \, ,
\end{equation}
where $\mathbf{d}_{i}$ is solution of
\begin{equation} \label{abstract-fixed-point}
\mathbf{d}_{i} = \Delta t \, \mathbf{U}(\mathbf{x}_{i}-\mathbf{d}_{i},t_{n}) \, .
\end{equation}
In many cases, $\mathbf{U}(\cdot,t_{n})$ is only known at points $\mathbf{x}_{i}$. Then we have to replace (\ref{abstract-fixed-point}) by
\begin{equation} \label{abstract-fix-pt-interp}
\mathbf{d}_{i} = \Delta t \, \Pi \mathbf{U}(\mathbf{x}_{i}-\mathbf{d}_{i},t_{n}) \, ,
\end{equation}
where $\Pi$ is an interpolation operator on points $(\mathbf{x}_{i})_{i}$. Assuming that the polynomial function $\mathbf{x} \mapsto \Pi \mathbf{U}(\mathbf{x},t_{n})$ is regular enough, we write the following expansion of (\ref{abstract-fix-pt-interp}):
\begin{equation} \label{abstract-expansion}
\mathbf{d}_{i} = \Delta t \, \mathbf{U}(\mathbf{x}_{i},t_{n}) - \Delta t \, \nabla_{\mathbf{x}} (\Pi \, \mathbf{U}) (\mathbf{x}_{i},t_{n}) \, \mathbf{d}_{i} + \mathcal{O}(\Delta t^{2}) \, ,
\end{equation}
because, in the expansion of $\Pi \mathbf{U}$ in $\mathbf{x}$, we get a $\mathcal{O}\big(|\mathbf{d}_{i}|^{2}\big)$ term, which is a $\mathcal{O}(\Delta t^{2})$. Then, we obtain a second order accurate approximation of $\mathbf{d}_{i}$ given by
\begin{equation} \label{abstract-fix-pt-solution}
\mathbf{d}_{i} = \Delta t \, \big( \mathbf{Id} + \Delta t \, \nabla_{\mathbf{x}} (\Pi \, \mathbf{U}) (\mathbf{x}_{i},t_{n}) \big)^{-1} \times \mathbf{U}(\mathbf{x}_{i},t_{n}) \, .
\end{equation}

\indent Considering this approach, the semi-lagrangian method writes:
\begin{enumerate}
\item knowing $f$ at time $t_{n}-\Delta t$ and $\mathbf{U}$ at time $t_{n}$, we compute $(\mathbf{d}_{i})_{i\,=\,0,\,\dots,\,N}$ by using the relation (\ref{abstract-fix-pt-solution}) for each $i$,
\item we compute $f$ at time $t_{n}+\Delta t$ as follows:
\begin{equation}
f(\mathbf{x}_{i},t_{n}+\Delta t) = \Pi f (\mathbf{x}_{i}-2\,\mathbf{d}_{i},t_{n}-\Delta t) \, .
\end{equation}
\end{enumerate}

\subsection{Implementation of the non-homogenized model}

In this paragraph, we describe a classical semi-lagrangian method on the system (\ref{NH-polar}). Since the electric fields $E_{r}^{\epsilon}$ and $\Xi_{r}^{\epsilon}$ do not depend on $v_{r}$, we can do a time-splitting on the first equation of the model, i.e. solving separately at each time step the equations
\begin{equation} \label{time-split_1}
\D_{t} f^{\epsilon} + \frac{v_{r}}{\epsilon} \, \D_{r}f^{\epsilon} = 0 \, ,
\end{equation}
and
\begin{equation} \label{time-split_2}
\D_{t} f^{\epsilon} + \big(E_{r}^{\epsilon} + \Xi_{r}^{\epsilon}) \, \D_{v_{r}}f^{\epsilon} = 0 \, ,
\end{equation}
with a second order in time numerical scheme instead of solving the complete equation with the same scheme. As a consequence, we only do 1D interpolations instead of 2D interpolations. Then, denoting $f_{n}^{\epsilon}$ with the aproximation of $f^{\epsilon}(\cdot,\cdot,t_{n})$ and $E_{n}^{\epsilon}$ with the approximation of $E_{r}^{n}(\cdot,t_{n})$ on the uniform mesh $(r_{i},{v_{r}}_{j})_{i,j}$ where $\Delta r$ is the size of one cell in $r$ direction and $\Delta v_{r}$ is the size of one cell in $v_{r}$ direction, an iteration of the method is organized as follows:
\begin{enumerate}
\item Knowing $f_{n}^{\epsilon}$ and $E_{n}^{\epsilon}$, we do a backward advection of $\frac{\Delta t}{2}$ in $v_{r}$ direction and we define $f_{*}^{\epsilon}$ by
\begin{equation}
f_{*}^{\epsilon}(r_{i},{v_{r}}_{j}) = \Pi_{v_{r}} f_{n}^{\epsilon}\Big(r_{i},{v_{r}}_{j}-\frac{\Delta t}{2} \big(E_{n}^{\epsilon}(r_{i})+\Xi_{r}^{\epsilon}(r_{i},t_{n})\big) \Big) \, ,
\end{equation}
where $\Pi_{v_{r}}$ is a 1D cubic spline interpolation operator on the points $({v_{r}}_{j})_{j}$.
\item We do an advection of $\Delta t$ in $r$ direction and we define $f_{**}^{\epsilon}$ by
\begin{equation}
f_{**}^{\epsilon}(r_{i},{v_{r}}_{j}) = \Pi_{r} f_{*}^{\epsilon}\Big(r_{i}-\frac{\Delta t}{\epsilon} \, {v_{r}}_{j}, {v_{r}}_{j} \Big) \, ,
\end{equation}
where $\Pi_{r}$ is a 1D cubic spline interpolation operator on the points $(r_{i})_{i}$.
\item We compute $E_{n+1}^{\epsilon}$ by discretizing the formula
\begin{equation}
E_{n+1}^{\epsilon}(r_{i}) = \left\{
\begin{array}{ll}
\displaystyle \frac{1}{r_{i}} \int_{0}^{r_{i}} \int_{\R} s \, f_{**}^{\epsilon}(s,v_{r}) \, dv_{r} \, ds & \textit{if $r_{i} \neq 0$,} \\ \\
0 & \textit{otherwise,}
\end{array}
\right.
\end{equation}
with the trapezoidal rule on the points $(r_{i},{v_{r}}_{j})_{i,j}$.
\item We compute $f_{n+1}^{\epsilon}$ by doing a last advection of $\frac{\Delta t}{2}$ in $v_{r}$ direction:
\begin{equation}
f_{n+1}^{\epsilon}(r_{i},{v_{r}}_{j}) = \Pi_{v_{r}} f_{**}^{\epsilon}\Big(r_{i},{v_{r}}_{j}-\frac{\Delta t}{2} \big(E_{n+1}^{\epsilon}(r_{i})+\Xi_{r}^{\epsilon}(r_{i},t_{n+1})\big) \Big) \, .
\end{equation}
\end{enumerate}

If we use such a method, we have to guarantee the accuracy of the scheme, especially we consider the $\frac{1}{\epsilon}$-frequency oscillations of the external electric field. A solution is to assume that the time step $\Delta t$ satisfies
\begin{equation} \label{CFL-NH}
\left\{
\begin{array}{rcccl}
r_{i} - \Delta r & \leq & r_{i} - \cfrac{\Delta t}{\epsilon} \, \, {v_{r}}_{j} & \leq & r_{i} + \Delta r \, , \\
{v_{r}}_{j} - \Delta v_{r} & \leq & {v_{r}}_{j} - \cfrac{\Delta t}{2} \, \big(E_{n}^{\epsilon}(r_{i}) + \Xi^{\epsilon}(r_{i},t_{n})\big) & \leq & {v_{r}}_{j} + \Delta v_{r} \, ,
\end{array}
\right.
\end{equation}
for all $i,j,n$. Then, in order to obtain good results with this method, $\Delta t$ has to be of the same order of $\epsilon$ which penalizes the method in terms of CPU time cost when we consider a very small $\epsilon$.

\subsection{Implementation of the two-scale model}

In this paragraph, we adapt the semi-lagrangian method we have described in the paragraph 3.1 to the model (\ref{H-polar}). In this case, the characteristics of the system are the solutions of
\begin{equation} \label{charac-ts}
\left\{
\begin{array}{rcl}
\D_{t}Q(t) &=& \langle \mathcal{E}_{1} \rangle \big(Q(t),U_{r}(t),t\big) \, ,\\
\D_{t}U_{r}(t) &=& \langle \mathcal{E}_{2} \rangle \big(Q(t),U_{r}(t),t\big) \, ,
\end{array}
\right.
\end{equation}
where $\langle \mathcal{E}_{1} \rangle$ and $\langle \mathcal{E}_{2} \rangle$ are defined by
\begin{equation} \label{def-E1-E2}
\left\{
\begin{array}{rcl}
\langle \mathcal{E}_{1} \rangle (q,u_{r},t) &=& \displaystyle - \int_{0}^{2\pi} \sin(\sigma) \Big[ \mathcal{E}_{r}\big(\cos(\sigma) \, q + \sin(\sigma)\,u_{r}, \sigma, t\big) \\
& & \qquad \qquad \qquad \qquad \displaystyle + \frac{I_{\Q}(\omega_{1})}{2\pi} \, H_{1}(\omega_{1} \, \sigma) \, \big(\cos(\sigma) \, q + \sin(\sigma)\,u_{r}\big) \Big] \, d\sigma \, ,\\
\langle \mathcal{E}_{2} \rangle (q,u_{r},t) &=& \displaystyle \int_{0}^{2\pi} \cos(\sigma) \Big[ \mathcal{E}_{r}\big(\cos(\sigma) \, q + \sin(\sigma)\,u_{r}, \sigma, t\big) \\
& & \qquad \qquad \qquad \qquad \displaystyle + \frac{I_{\Q}(\omega_{1})}{2\pi} \, H_{1}(\omega_{1} \, \sigma) \, \big(\cos(\sigma) \, q + \sin(\sigma)\,u_{r}\big) \Big] \, d\sigma \, .
\end{array}
\right.
\end{equation}
As in the paragraph 3.1, we remark that the solution $G$ of (\ref{H-polar}) is constant along the characteristics, so we can write:
\begin{equation} \label{conservation-G}
G(q,u_{r},t_{n+1}) = G\big(Q(t_{n-1}; q, t_{n+1}), U_{r}(t_{n-1}; u_{r},t_{n+1}), t_{n-1}\big) \, ,
\end{equation}
where $t_{n} = n \, \Delta t$, and where $\big(Q(t_{n-1}; q, t_{n+1}), U_{r}(t_{n-1}; u_{r},t_{n+1})\big)$ is the solution of (\ref{charac-ts}) with the condition $\big(Q(t_{n+1}),U_{r}(t_{n+1})\big) = (q,u_{r})$. \\

\indent Firstly, we define the mesh in $q$ and $u_{r}$ directions by considering the points $q_{i} = i \, \Delta q$ and ${u_{r}}_{j} = j \, \Delta u_{r}$ ($i = -P_{q},\dots,P_{q}$, $j = -P_{u_{r}},\dots,P_{u_{r}}$). We also consider the uniform mesh $\tau_{m} = m \, \Delta \tau$ on $[0,2\pi]$ ($m = 0,\dots,P_{\tau}$). Finally, we fix a time step $\Delta t$ for the entire simulation. Then, denoting $G^{n}$ with the approximation of $G$ at time $t_{n}$, an iteration of the semi-lagrangian method is organized as follows:
\begin{enumerate}
\item Assuming that we know the value of $G^{n}$ and $G^{n-1}$ on the mesh $(q_{i},{u_{r}}_{j})_{i,j}$, we compute $\mathcal{E}_{r}^{n}$ with the formula
\begin{equation} \label{Poisson-tsu}
\mathcal{E}_{r}(q_{i},\tau_{m}) = \left\{
\begin{array}{ll}
\displaystyle \frac{1}{q_{i}} \int_{0}^{q_{i}}  \int_{\R} s \, G^{n}\big(\cos(\tau_{m}) s - \sin(\tau_{m}) v_{r}, & \\
\qquad \qquad \qquad \qquad \qquad \sin(\tau_{m})s + \cos(\tau_{m})v_{r} \big) \,ds\,dv_{r} & \textit{if $i \neq 0$, } \\ \\
0 & \textit{otherwise.}
\end{array}
\right.
\end{equation}
Since $G^{n}$ is only known at points $(q_{i},{u_{r}}_{j})$, we have to interpolate $G^{n}$. Assuming that the support of $G^{n}$ is included in $[-(P_{q}+1)\Delta q, (P_{q}+1)\Delta q] \times [-(P_{u_{r}}+1)\Delta u_{r}, (P_{u_{r}}+1)\Delta u_{r}]$ with $P_{q},P_{u_{r}} \in \N$ large enough, we use the trapezoidal rule in order to approach the integral above. We obtain
\begin{equation} \label{Poisson-tsu-discrete}
\begin{split}
\mathcal{E}_{r}^{n}(q_{i},\tau_{m}) &\approx \frac{\Delta q \, \Delta u_{r}}{2} \hspace{-0.2cm} \sum_{j \, = \, -P_{u_{r}}}^{P_{u_{r}}} \hspace{-0.2cm} \Bigg( \Pi_{2}G^{n}(\cos(\tau_{m}) q_{i} - \sin(\tau_{m}) {u_{r}}_{j}, \sin(\tau_{m})q_{i} + \cos(\tau_{m}){u_{r}}_{j} \big) \\
&\qquad + \frac{2}{i} \, \sum_{k \,=\,1}^{i-1} k \, \Pi_{2}G^{n}(\cos(\tau_{m}) q_{k} - \sin(\tau_{m}) {u_{r}}_{j}, \sin(\tau_{m})q_{k} + \cos(\tau_{m}){u_{r}}_{j} \big) \Bigg) \, ,
\end{split}
\end{equation}
where $\Pi_{2}$ is a cubic spline interpolation operator on the points $(q_{i},{u_{r}}_{j})$.

\item We compute $\langle \mathcal{E}_{1}^{n} \rangle$ and $\langle \mathcal{E}_{2}^{n} \rangle$ at points $(q_{i},{u_{r}}_{j})$: since $\mathcal{E}_{r}^{n}$ is only known at points $(q_{i},\tau_{m})$, we have to interpolate $\mathcal{E}_{r}^{n}$. By using the trapezoidal rule for approximating the integrals in (\ref{def-E1-E2}), we obtain
\begin{equation} \label{E1-E2-tsu}
\left\{
\begin{array}{rcl}
\langle \mathcal{E}_{1}^{n} \rangle (q_{i},{u_{r}}_{j}) & \approx & \displaystyle - \Delta\tau \sum_{m\,=\,0}^{P_{\tau}} \sin(\tau_{m}) \Big[ \Pi_{1} \mathcal{E}_{r}^{n}\big(\cos(\tau_{m}) \, q_{i} + \sin(\tau_{m})\,{u_{r}}_{j}, \tau_{m}\big) \\
& & \qquad \qquad \quad \displaystyle + \frac{I_{\Q}(\omega_{1})}{2\pi} \, H_{1}(\omega_{1} \, \tau_{m}) \, \big(\cos(\tau_{m}) \, q_{i} + \sin(\tau_{m})\,{u_{r}}_{j}\big) \Big] \, , \\
\langle \mathcal{E}_{2}^{n} \rangle (q_{i},{u_{r}}_{j}) & \approx & \displaystyle \Delta\tau \sum_{m\,=\,0}^{P_{\tau}} \cos(\tau_{m}) \Big[ \Pi_{1} \mathcal{E}_{r}^{n}\big(\cos(\tau_{m}) \, q_{i} + \sin(\tau_{m})\,{u_{r}}_{j}, \tau_{m} \big) \\
& & \qquad \qquad \quad \displaystyle + \frac{I_{\Q}(\omega_{1})}{2\pi} \, H_{1}(\omega_{1} \, \tau_{m}) \, \big(\cos(\tau_{m}) \, q_{i} + \sin(\tau_{m})\,{u_{r}}_{j}\big) \Big] \, ,
\end{array}
\right.
\end{equation}
where $\Pi_{1}$ is a cubic spline interpolation operator on the points $(q_{i})_{i}$.

\item We compute the the shifts $(d_{i,j}^{1},d_{i,j}^{2})$ by using the following formula:
\begin{equation} \label{tsu-fix-point}
\left(
\begin{array}{c}
d_{i,j}^{1} \\ d_{i,j}^{2}
\end{array}
\right) = \Delta t \, \mathbf{A}_{i,j}^{-1} \, \left(
\begin{array}{c}
\langle \mathcal{E}_{n}^{1} \rangle (q_{i},{u_{r}}_{j}) \\
\langle \mathcal{E}_{n}^{2} \rangle (q_{i},{u_{r}}_{j})
\end{array}
\right) \, ,
\end{equation}
where the matrix $\mathbf{A}_{i,j}$ is defined by
\begin{equation}
\mathbf{A}_{i,j} = \mathbf{Id} + \Delta t \, \left(
\begin{array}{cc}
\D_{q}\big(\Pi_{2} \langle \mathcal{E}_{n}^{1} \rangle \big) (q_{i},{u_{r}}_{j}) & \D_{u_{r}} \big(\Pi_{2} \langle \mathcal{E}_{n}^{1} \rangle \big) (q_{i},{u_{r}}_{j}) \\ \\
\D_{q} \big(\Pi_{2} \langle \mathcal{E}_{n}^{2} \rangle\big) (q_{i},{u_{r}}_{j}) & \D_{u_{r}}\big(\Pi_{2} \langle \mathcal{E}_{n}^{2} \rangle\big) (q_{i},{u_{r}}_{j})
\end{array}
\right) \, .
\end{equation}
\item We compute $G^{n+1}$ by interpolating $G^{n-1}$ at the points $(q_{i} - 2 \, d_{i,j}^{1}, {u_{r}}_{j} - 2 \, d_{i,j}^{2})$:
\begin{equation}
G^{n+1}(q_{i},{u_{r}}_{j}) = \Pi_{2} G^{n-1}(q_{i}-2\,d_{i,j}^{1}, {u_{r}}_{j}-2\,d_{i,j}^{2}) \, .
\end{equation}
\item We save the approximation of $f^{\epsilon}$ at time $t_{n+1}$ given by
\begin{equation} \label{approx_tsu}
\textstyle f^{\epsilon}(r,v_{r},t_{n+1}) \sim 2\,\pi \, \Pi_{2}G^{n+1}\Big( \cos\big(\frac{t_{n+1}}{\epsilon}\big) \, r - \sin\big(\frac{t_{n+1}}{\epsilon}\big) \, v_{r} , \sin\big(\frac{t_{n+1}}{\epsilon}\big) \, r + \cos\big(\frac{t_{n+1}}{\epsilon}\big) \, v_{r} \Big) \, .
\end{equation}
\end{enumerate}

In order to initialize this two time step advance, we have to compute $G^{1}$ from $G^{0}$ which is given as an initial data: for this purpose, we perform a complete iteration such as decribed above where $\Delta t$ is replaced by $\frac{\Delta t}{2}$ and where we assume that $G^{1/2} = G^{0}$.

\subsection{The two-scale mesh}

\indent In most of test cases, we assume that the support of the initial distribution $f_{0}$ is compact and is included in some $\Omega = [-R,R] \times [-v_{R},v_{R}] \subset \R^{2}$ for $R > 0$, $v_{R} > 0$ large enough. Then, if we follow the algorithm presented in the previous paragraph, the first thing we have to compute is $\mathcal{E}_{r}^{0}$ by approximating the integral in (\ref{Poisson-tsu}): for numerical reasons, we have to reduce the integral on $\R$ to an integral on a compact interval. Furthermore, since we can only say that the support of $(r,v_{r}) \mapsto f_{0}\big(\cos(\tau)\,r-\sin(\tau)\,v_{r},\sin(\tau)\,r+\cos(\tau)\,v_{r}\big)$ is included in $\Omega' = [-R-v_{R},R+v_{R}] \times [-R-v_{R},R+v_{R}]$ as it is illustrated in Figure 1, the integral on $\R$ in (\ref{Poisson-tsu}) is reduced to an integral on $[-R-v_{R},R+v_{R}]$. As a consequence, we have to do all the simulation on $\Omega'$ instead of $\Omega$ in order to avoid losing some data and we have to increase the number of mesh points in order to keep good interpolation results, even if the distribution function will be equal to 0 at many new points. \\

\begin{center}
\psfrag{r}{$R$}
\psfrag{vr}{$v_{R}$}
\psfrag{mr}{$-R$}
\psfrag{mvr}{$-v_{R}$}
\psfrag{mrmvr}{$-R-v_{R}$}
\psfrag{ravr}{$R+v_{R}$}
\psfrag{vrar}{$v_{R}+R$}
\psfrag{mvrmr}{$-v_{R}-R$}
\psfrag{rr}{$r$}
\psfrag{vrvr}{$v_{r}$}
\includegraphics{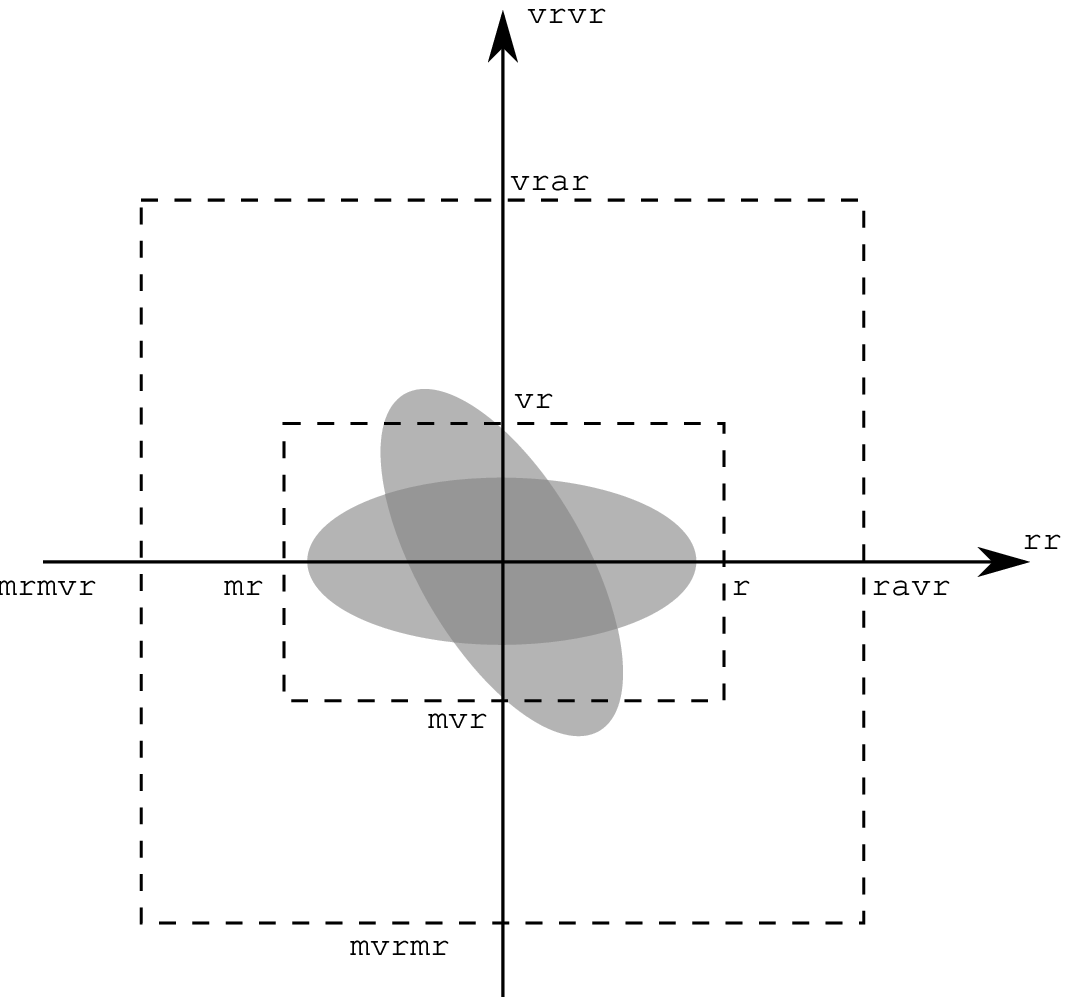}
\caption{Support of $(r,v_{r}) \mapsto f_{0}(r,v_{r})$ and $(r,v_{r}) \mapsto f_{0}\big(\cos(\tau)\,r-\sin(\tau)\,v_{r},\sin(\tau)\,r+\cos(\tau)\,v_{r}\big)$ with $\tau = \frac{2\pi}{3}$.}
\end{center}

\indent In this paragraph, we present a different approach in order to avoid this extension of the simulation domain and its mesh. Before explaining the main idea of this new method, we consider the following meshes on $\Omega$ and $[0,2\pi]$:
\begin{equation}
\begin{array}{c}
M(\Omega) = \Big\{ (r_{i},{v_{r}}_{j}) = (i\,\Delta r, j \, \Delta v_{r}) \, : \, i = -P_{r},\dots,P_{r} \, , \, j = -P_{v_{r}},\dots,P_{v_{r}} \Big\} \, , \\ \\
M\big([0,2\pi]\big) = \Big\{ \tau_{m} = m \, \Delta \tau \, : \, m = 0,\dots,P_{\tau} \Big\} \, , \\ \\
M\big([-R,R]\big) = \Big\{ r_{i} = i \, \Delta r \, : \, i = -P_{r},\dots,P_{r} \Big\} \, ,
\end{array}
\end{equation}
where $\Delta r = \frac{R}{P_{r}+1}$, $\Delta v_{r} = \frac{v_{R}}{P_{v_{r}}+1}$ and $\Delta \tau = \frac{2\pi}{P_{\tau}+1}$. Considering the function $\mathbf{\gamma}$ defined by
\begin{equation}
\begin{array}{rccc}
\mathbf{\gamma} \, : \, & \R^{2} \times [0,2\pi] & \longrightarrow & \R^{2} \\
               & (r,v_{r},\tau)         & \longmapsto     & \big( \cos(\tau) \, r - \sin(\tau) \, v_{r} , \sin(\tau) \, r + \cos(\tau) \, v_{r} \big)
\end{array}
\, ,
\end{equation}
we define $\Omega(\tau)$ and $M\big(\Omega(\tau)\big)$ by
\begin{equation}
\begin{array}{c}
\Omega(\tau) = \mathbf{\gamma}\big(\Omega \times \{\tau\}\big) \subset \R^{2} \, , \\ \\
M\big(\Omega(\tau)\big) = \Big\{ \mathbf{\gamma}(r_{i},{v_{r}}_{j},\tau) \, : \, i = -P_{r},\dots,P_{r} \, , \, j = -P_{v_{r}},\dots,P_{v_{r}} \Big\} \, .
\end{array}
\end{equation}

\begin{center}
\psfrag{r}{$r$}
\psfrag{vr}{$v_{r}$}
\includegraphics[scale=0.8]{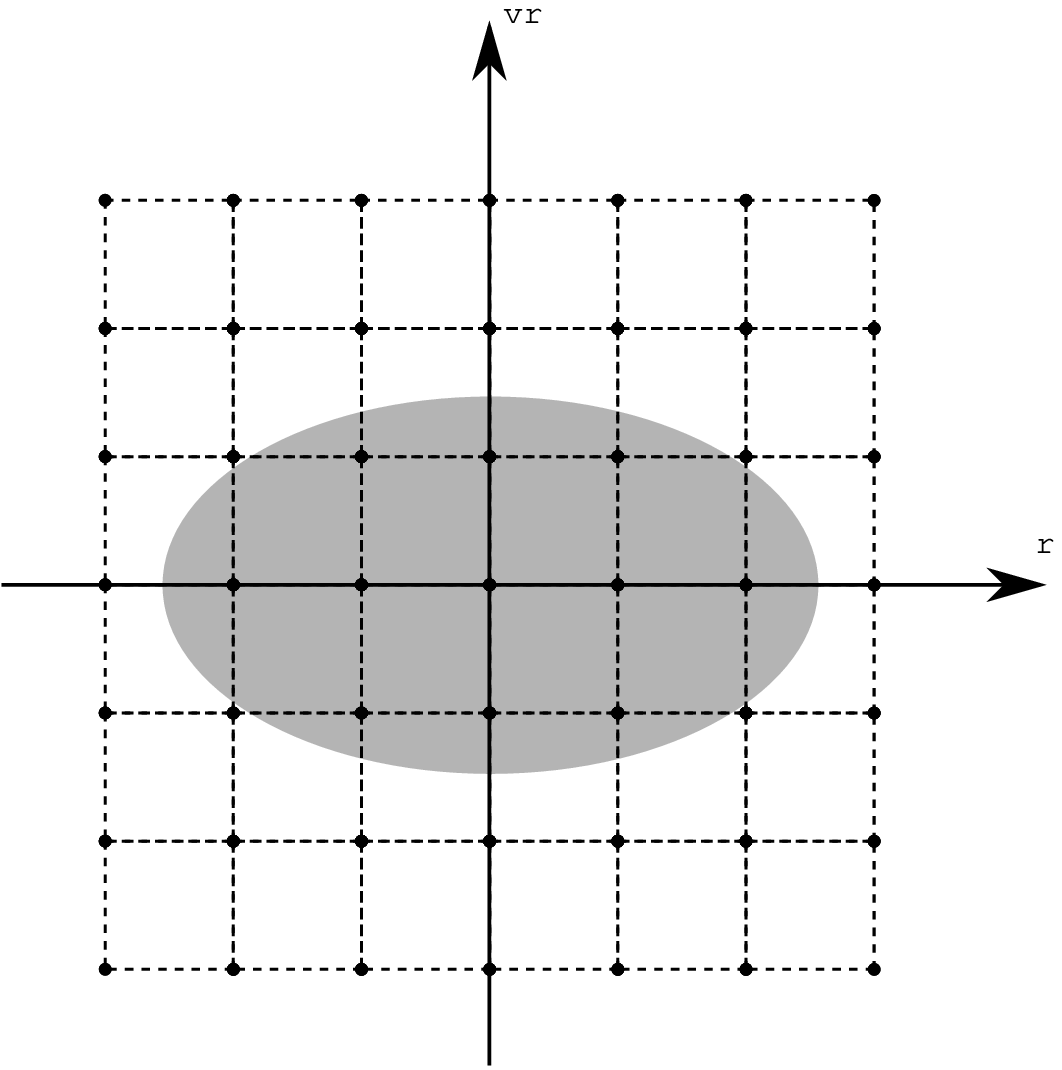}
\caption{Mesh $M(\Omega) = M\big(\Omega(0)\big)$ and support of $(r,v_{r}) \mapsto f_{0}(r,v_{r})$.}
\end{center}
\begin{center}
\psfrag{r}{$r$}
\psfrag{vr}{$v_{r}$}
\includegraphics[scale=0.8]{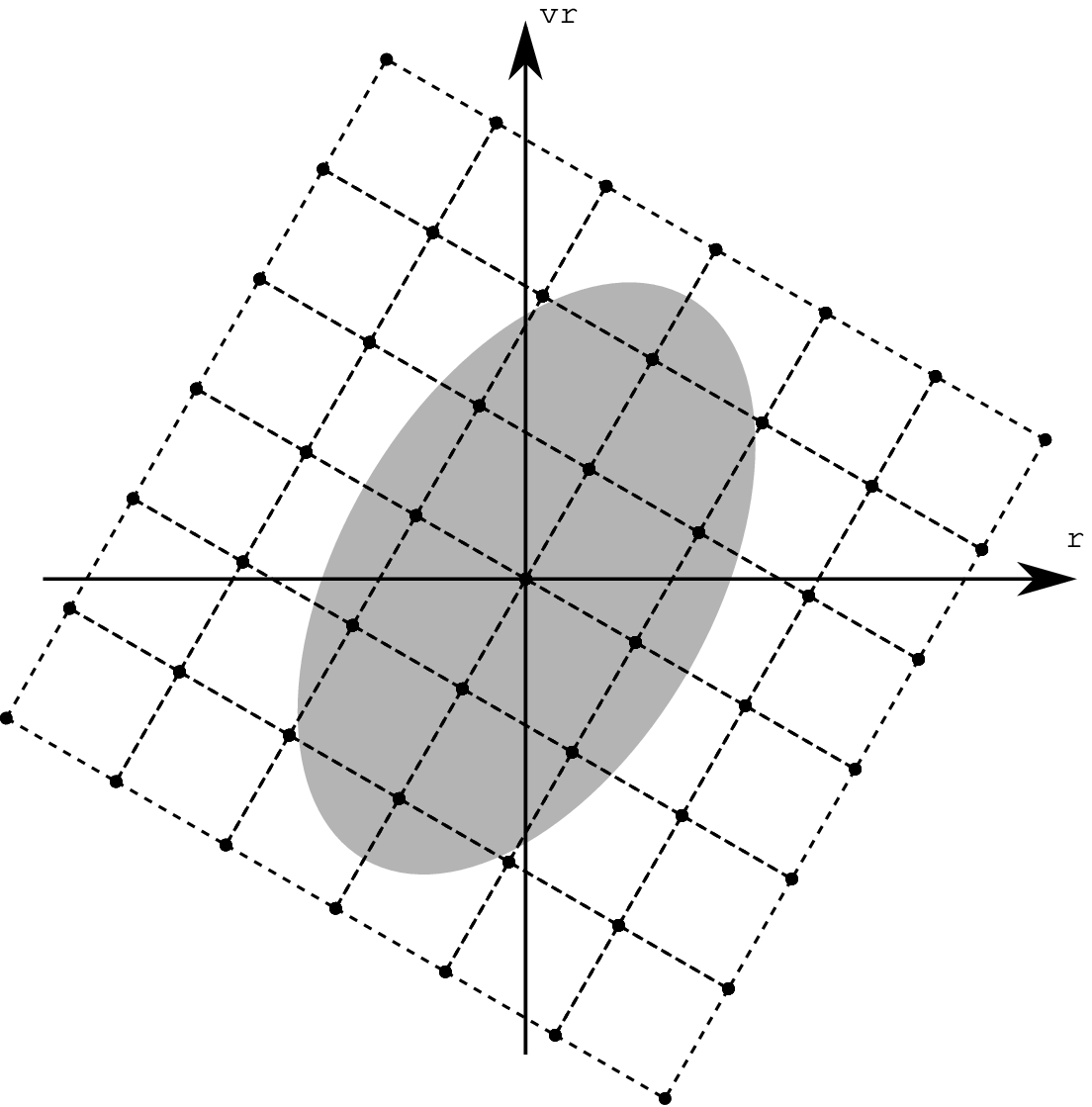}
\caption{Mesh $M\big(\Omega(\frac{\pi}{3})\big)$ and support of $(r,v_{r}) \mapsto f_{0}\circ \mathbf{\gamma}(r,v_{r},\frac{\pi}{3})$.}
\end{center}
\begin{center}
\psfrag{r}{$r$}
\psfrag{vr}{$v_{r}$}
\includegraphics[scale=0.8]{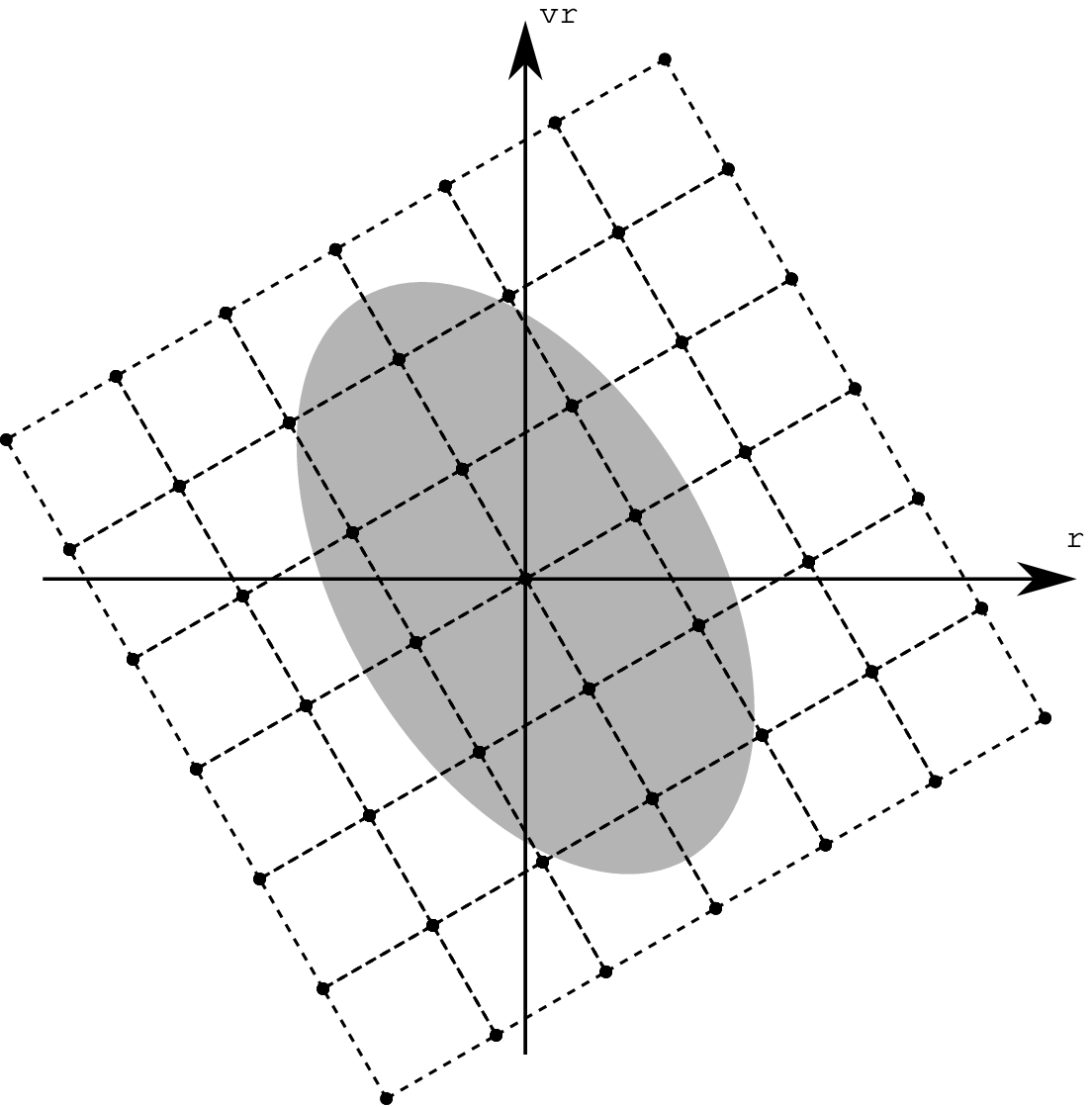}
\caption{Mesh $M\big(\Omega(\frac{2\pi}{3})\big)$ and support of $(r,v_{r}) \mapsto f_{0}\circ \mathbf{\gamma}(r,v_{r},\frac{2\pi}{3})$.}
\end{center}

The main idea of our new method is to compute $(r,v_{r}) \mapsto F(r,v_{r},\tau_{m},t_{n})$ at points $(r_{i},{v_{r}}_{j}) \in M(\Omega)$ and to compute $r \mapsto \mathcal{E}_{r}(r,\tau_{m},t_{n})$ at points $r_{i} \in M\big([-R,R]\big)$, whereas the function $(q,u_{r}) \mapsto G(q,u_{r},t_{n})$ is computed at points $\mathbf{\gamma}(r_{i},{v_{r}}_{j},\tau_{m}) \in M\big(\Omega(\tau_{m})\big)$ for every $\tau_{m} \in M\big([0,2\pi]\big)$. This approach is similar as the time-dependent moving grid described by Lang \textit{et al.} in \cite{Lang}, even if, in our case, the mesh $M\big(\Omega(\tau)\big)$ only depends on $M(\Omega)$ and $\tau$, and is completely defined before the beginning of the simulation. \\
\indent As a first consequence, considering that the support of $(r,v_{r}) \mapsto f_{0}(r,v_{r})$ is included in $\Omega$ is equivalent to considering that the support of $(r,v_{r}) \mapsto f_{0}\big(\mathbf{\gamma}(r,v_{r},\tau)\big)$ is included in $\Omega(\tau)$ for any $\tau \in [0,2\pi]$ such as illustrated in Figures 2,3, and 4, where the support of a Kapchinsky-Vladimirsky distribution is represented (see \cite{KV}, \cite{Reiser} or \cite{Paraxial} for more details about this distribution). Then we do not have to extend $\Omega$ in order to avoid losing some data. Furthermore, the equation (\ref{FG-polar}) reads
\begin{equation} \label{FG-tsts}
F(r,v_{r},\tau,t) = G \big(\mathbf{\gamma}(r,v_{r},\tau), t\big) \, .
\end{equation}

\indent Considering this approach, we fix a time step $\Delta t$ for all the simulation. Then, an iteration of this new semi-lagrangian method is organized as follows:
\begin{enumerate}
\item Assuming that $G^{n-1}$ and $G^{n}$ are known on the mesh $M\big(\Omega(\tau_{m})\big)$ for all $\tau_{m} \in M\big([0,2\pi]\big)$, we compute $\mathcal{E}_{r}^{n}$ at points $(r_{i},\tau_{m}) \in M\big([-R,R]\big) \times M\big([0,2\pi]\big)$. With the new notations, the equation (\ref{Poisson-tsu}) simplifies itself to
\begin{equation} \label{Poisson-tsts}
\mathcal{E}_{r}^{n}(r_{i},\tau_{m}) = \left\{
\begin{array}{ll}
\displaystyle \frac{1}{r_{i}} \int_{0}^{r_{i}} \int_{\R} s \, G^{n} \big( \mathbf{\gamma}(s,v_{r},\tau_{m}) \big) \, dv_{r} \, ds & \textit{if $i \neq 0$,} \\ \\
0 & \textit{otherwise.}
\end{array}
\right.
\end{equation}
Assuming that the support of $G^{n}\big(\mathbf{\gamma}(\cdot,\cdot,\tau_{m})\big)$ is included in $\Omega(\tau_{m})$ for each $\tau_{m}$ (which is equivalent to assume that the support of $G^{n}\big(\mathbf{\gamma}(\cdot,\cdot,0)\big)$ is included in $\Omega$), we use the trapezoidal rule to approximate the integral above:
\begin{equation}
\begin{split}
\mathcal{E}_{r}^{n}(r_{i},\tau_{m}) &\approx \frac{\Delta r \, \Delta v_{r}}{2} \, \sum_{j \, = \, -P_{v_{r}}}^{P_{v_{r}}} \Bigg( G^{n} \big( \mathbf{\gamma} (r_{i},{v_{r}}_{j},\tau_{m}) \big) + \frac{2}{i} \, \sum_{k \,=\,1}^{i-1} k \, G^{n} \big(\mathbf{\gamma} (r_{k},{v_{r}}_{j},\tau_{m}) \big) \Bigg) \, .
\end{split}
\end{equation}
We remark here that, contrary to the computation done in (\ref{Poisson-tsu-discrete}), we do not have to interpolate $G^{n}$.

\item We compute $\langle \mathcal{E}_{1}^{n} \rangle$ and $\langle \mathcal{E}_{2}^{n} \rangle$ at points $\mathbf{\gamma}(r_{i},{v_{r}}_{j},\tau_{m}) \in M\big(\Omega(\tau_{m})\big)$:
\begin{equation} \label{E1-E2-tsts}
\hspace{-0.4cm} \left\{
\begin{array}{rcl}
\hspace{-0.2cm} \langle \mathcal{E}_{1}^{n} \rangle \big( \mathbf{\gamma} (r_{i},{v_{r}}_{j},\tau_{m}) \big) \hspace{-0.1cm} &=& \hspace{-0.2cm} \displaystyle - \int_{0}^{2\pi} \sin(\sigma) \Big[ \mathcal{E}_{r}^{n}\big(\cos(\sigma-\tau_{m}) \, r_{i} + \sin(\sigma-\tau_{m})\,{v_{r}}_{j}, \sigma\big) \\
& & \displaystyle + \frac{I_{\Q}(\omega_{1})}{2\pi} \, H_{1}(\omega_{1} \, \sigma) \, \big(\cos(\sigma-\tau_{m}) \, r_{i} + \sin(\sigma-\tau_{m})\,{v_{r}}_{j}\big) \Big] \, d\sigma \, ,\\
\hspace{-0.1cm} \langle \mathcal{E}_{2}^{n} \rangle \big( \mathbf{\gamma} (r_{i},{v_{r}}_{j},\tau_{m}) \big) \hspace{-0.2cm} &=& \hspace{-0.2cm} \displaystyle \int_{0}^{2\pi} \cos(\sigma) \Big[ \mathcal{E}_{r}^{n}\big(\cos(\sigma-\tau_{m}) \, r_{i} + \sin(\sigma-\tau_{m})\,{v_{r}}_{j}, \sigma\big) \\
& & \displaystyle + \frac{I_{\Q}(\omega_{1})}{2\pi} \, H_{1}(\omega_{1} \, \sigma) \, \big(\cos(\sigma-\tau_{m}) \, r_{i} + \sin(\sigma-\tau_{m})\,{v_{r}}_{j}\big) \Big] \, d\sigma \, .
\end{array}
\right.
\end{equation}
We approximate the integrals with the trapezoidal rule. However, since $\mathcal{E}_{r}^{n}$ is only known at points of $M\big([-R,R]\big) \times M\big([0,2\pi]\big)$, we have to interpolate it: for that we choose a cubic spline interpolation operator on the mesh $M\big([-R,R]\big)$ and we denote it with $\Pi_{1}$. Then we have the following approximations:
\begin{equation} \label{E1-E2-tsts-discrete}
\hspace{-0.4cm} \left\{
\begin{array}{rcl}
\hspace{-0.2cm} \langle \mathcal{E}_{1}^{n} \rangle \big( \mathbf{\gamma} (r_{i},{v_{r}}_{j},\tau_{m}) \big) \hspace{-0.3cm} &\approx & \hspace{-0.3cm} \displaystyle - \Delta \tau \hspace{-0.1cm} \sum_{k\,=\,0}^{P_{\tau}} \sin(\tau_{k}) \Big[ \Pi_{1} \mathcal{E}_{r}^{n}\big(\cos(\tau_{k}-\tau_{m}) \, r_{i} + \sin(\tau_{k}-\tau_{m})\,{v_{r}}_{j}, \tau_{k}\big) \\
& & \displaystyle \quad + \frac{I_{\Q}(\omega_{1})}{2\pi} \, H_{1}(\omega_{1} \, \tau_{k}) \, \big(\cos(\tau_{k}-\tau_{m}) \, r_{i} + \sin(\tau_{k}-\tau_{m})\,{v_{r}}_{j}\big) \Big] \, ,\\
\hspace{-0.2cm} \langle \mathcal{E}_{2}^{n} \rangle \big( \mathbf{\gamma} (r_{i},{v_{r}}_{j},\tau_{m}) \big) \hspace{-0.3cm} & \approx & \hspace{-0.3cm} \displaystyle \Delta\tau \hspace{-0.1cm} \sum_{k\,=\,0}^{P_{\tau}} \cos(\tau_{k}) \Big[ \Pi_{1} \mathcal{E}_{r}^{n}\big(\cos(\tau_{k}-\tau_{m}) \, r_{i} + \sin(\tau_{k}-\tau_{m})\,{v_{r}}_{j}, \tau_{k}\big) \\
& & \displaystyle \quad + \frac{I_{\Q}(\omega_{1})}{2\pi} \, H_{1}(\omega_{1} \, \tau_{k}) \, \big(\cos(\tau_{k}-\tau_{m}) \, r_{i} + \sin(\tau_{k}-\tau_{m})\,{v_{r}}_{j}\big) \Big] \, .
\end{array}
\right.
\end{equation}

\item We compute the shifts $\mathbf{d}(r_{i},{v_{r}}_{j},\tau_{m})$ verifying
\begin{equation}
\mathbf{d}(r_{i},{v_{r}}_{j},\tau_{m}) = \Delta t \, \left(
\begin{array}{c}
\langle \mathcal{E}_{1}^{n} \rangle \big( \mathbf{\gamma} (r_{i},{v_{r}}_{j},\tau_{m}) - \mathbf{d}(r_{i},{v_{r}}_{j},\tau_{m}) \big) \\
\langle \mathcal{E}_{2}^{n} \rangle \big( \mathbf{\gamma} (r_{i},{v_{r}}_{j},\tau_{m}) - \mathbf{d}(r_{i},{v_{r}}_{j},\tau_{m}) \big)
\end{array}
\right) \, .
\end{equation}
For that, we consider the cubic spline interpolation operator $\Pi_{2}^{m}$ on the mesh $M\big(\Omega(\tau_{m})\big)$ and, inspired by (\ref{abstract-fix-pt-solution}), we consider the following approximation of $\mathbf{d}(r_{i},{v_{r}}_{j},\tau_{m})$:
\begin{equation}
\mathbf{d}(r_{i},{v_{r}}_{j},\tau_{m}) = \Delta t \, \mathbf{A}_{i,j,m}^{-1} \, \left(
\begin{array}
{c}
\langle \mathcal{E}_{1}^{n} \rangle \big( \mathbf{\gamma} (r_{i},{v_{r}}_{j},\tau_{m}) \big) \\
\langle \mathcal{E}_{2}^{n} \rangle \big( \mathbf{\gamma} (r_{i},{v_{r}}_{j},\tau_{m}) \big)
\end{array}
\right) \, ,
\end{equation}
where the matrix $\mathbf{A}_{i,j,m}$ is defined by
\begin{equation}
\mathbf{A}_{i,j,m} = \mathbf{Id} + \Delta t \, \left(
\begin{array}{cc}
\D_{q} \big(\Pi_{2}^{m} \langle \mathcal{E}_{1}^{n} \rangle \big) \big(\mathbf{\gamma} (r_{i},{v_{r}}_{j},\tau_{m}) \big) & \D_{u_{r}} \big(\Pi_{2}^{m} \langle \mathcal{E}_{1}^{n} \rangle \big)\big(\mathbf{\gamma} (r_{i},{v_{r}}_{j},\tau_{m}) \big) \\ \\
\D_{q} \big(\Pi_{2}^{m} \langle \mathcal{E}_{2}^{n} \rangle \big) \big(\mathbf{\gamma} (r_{i},{v_{r}}_{j},\tau_{m}) \big) & \D_{u_{r}} \big(\Pi_{2}^{m} \langle \mathcal{E}_{2}^{n} \rangle\big) \big(\mathbf{\gamma} (r_{i},{v_{r}}_{j},\tau_{m}) \big)
\end{array}
\right) \, .
\end{equation}

\item We compute $G^{n+1}$ on the meshes $M\big(\Omega(\tau_{m})\big)$ for all $\tau_{m} \in M\big([0,2\pi]\big)$:
\begin{equation}
G^{n+1}\big(\mathbf{\gamma}(r_{i},{v_{r}}_{j},\tau_{m})\big) = \Pi_{2}^{m} G^{n-1}\big(\mathbf{\gamma}(r_{i},{v_{r}}_{j},\tau_{m}) - 2 \, \mathbf{d}(r_{i},{v_{r}}_{j}, \tau_{m})\big)
\end{equation}

\item Assuming that there exists a fixed integer $K \in \N$ such that
\begin{equation} \label{def_Delta_t}
\Delta t = \epsilon \, \Delta \tau \, K \, ,
\end{equation}
we save the approximation of $f^{\epsilon}$ on $M(\Omega)$ given by
\begin{equation} \label{approx_feps_tsts}
f^{\epsilon}(r_{i},{v_{r}}_{j},t_{n+1}) \sim 2\pi \, G^{n+1}\big(\gamma(r_{i},{v_{r}}_{j},\tau_{(n+1)\,K})\big) \, .
\end{equation}
Contrary to the scheme described in the previous paragraph, we do not have to interpolate $G^{n+1}$ to obtain an approximation of $f^{\epsilon}$ at time $t_{n+1}$. \\
\end{enumerate}
\indent Concerning the initialization of this two time step advance, we compute $G^{0} = \frac{1}{2\pi} \, f_{0}$ on the meshes $M\big(\Omega(\tau_{m})\big)$ for all $\tau_{m} \in M\big([0,2\pi]\big)$ and, assuming that $G^{1/2} = G^{0}$, we compute $G^{1}$ by performing a complete iteration such as described above with $\Delta t$ replaced by $\frac{\Delta t}{2}$. \\

\indent We have built another two-scale semi-lagrangian method based on a mesh depending on the variable $\tau$. Compared to the semi-lagrangian method described in the paragraph 3.3, this technique allows us to avoid some interpolations not only when we compute the two-scale limit of the electric field $\mathcal{E}$ but also when we build the approximation of $f^{\epsilon}$ defined by (\ref{approx_feps_tsts}). As a consequence, we reduces the global numerical diffusion within the simulation. However, we have to compute $G$ and $\langle \mathcal{E} \rangle$ on each mesh $M\big(\Omega(\tau_{m})\big)$, which can be expensive in CPU time.

\section{Numerical results}
\setcounter{equation}{0}

In this section, we present some numerical results obtained with the two-scale semi-lagrangian methods described in the paragraphs 3.3 and 3.4. Following the approach of \cite{PIC-two-scale} for validating our new methods, we study in a first time some linear cases where we can find an analytic expression of the solution $G$ of the system (\ref{H-polar}). Then we test the methods on non-linear cases.

\subsection{Linear cases}

\indent In order to validate the two-scale semi-lagrangian methods described in the paragraphs 3.3 and 3.4, we simulate some linear cases, i.e. with a self-consistent electric field set to $0$. This assumption allows us to compute analytically the solution $G$ of the system (\ref{H-polar}) under an adequate choice of $\omega_{1}$ and $H_{1}$. \\

\begin{center}
\begin{tabular}{cc}
\includegraphics[scale=0.15]{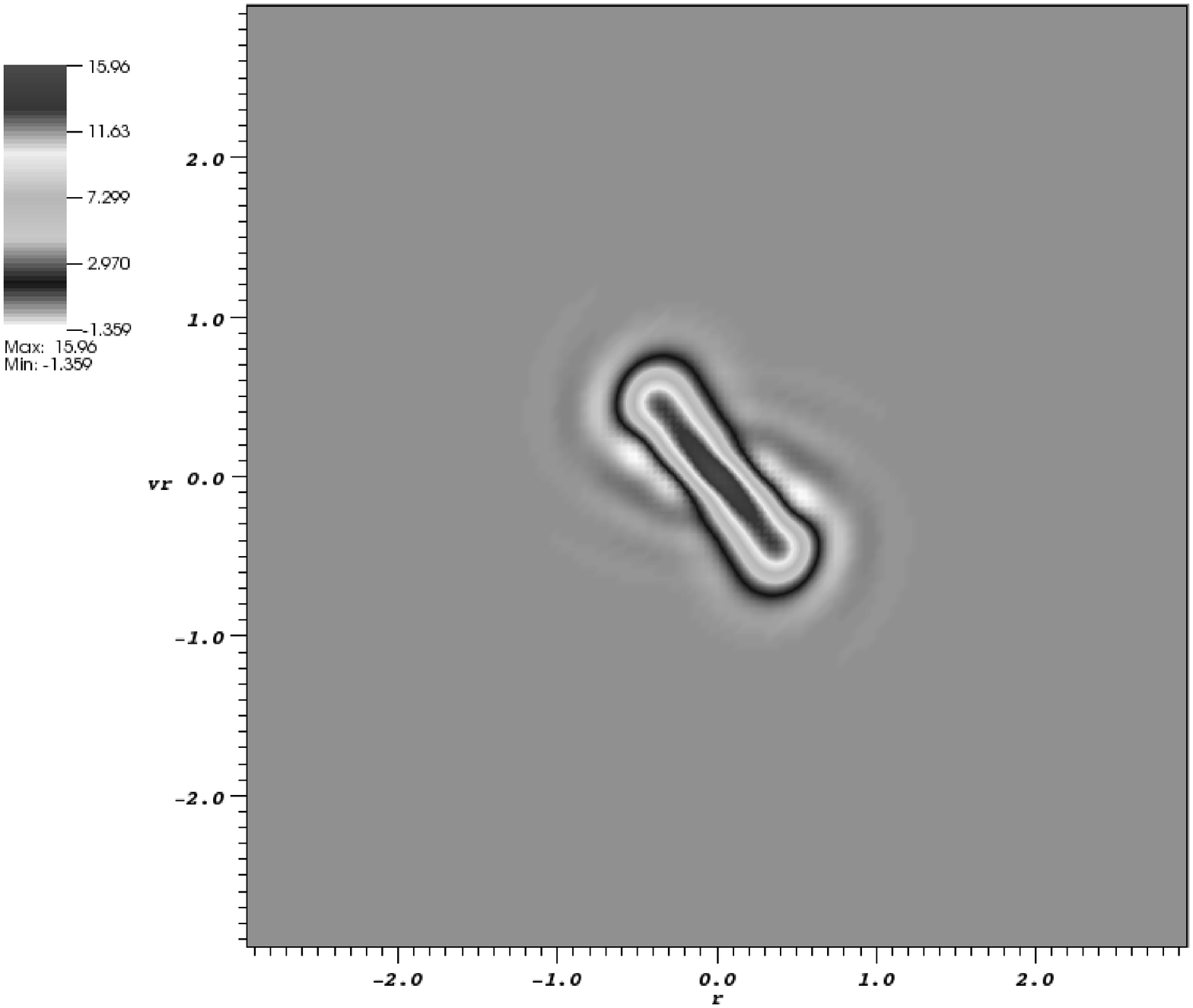} & \includegraphics[scale=0.15]{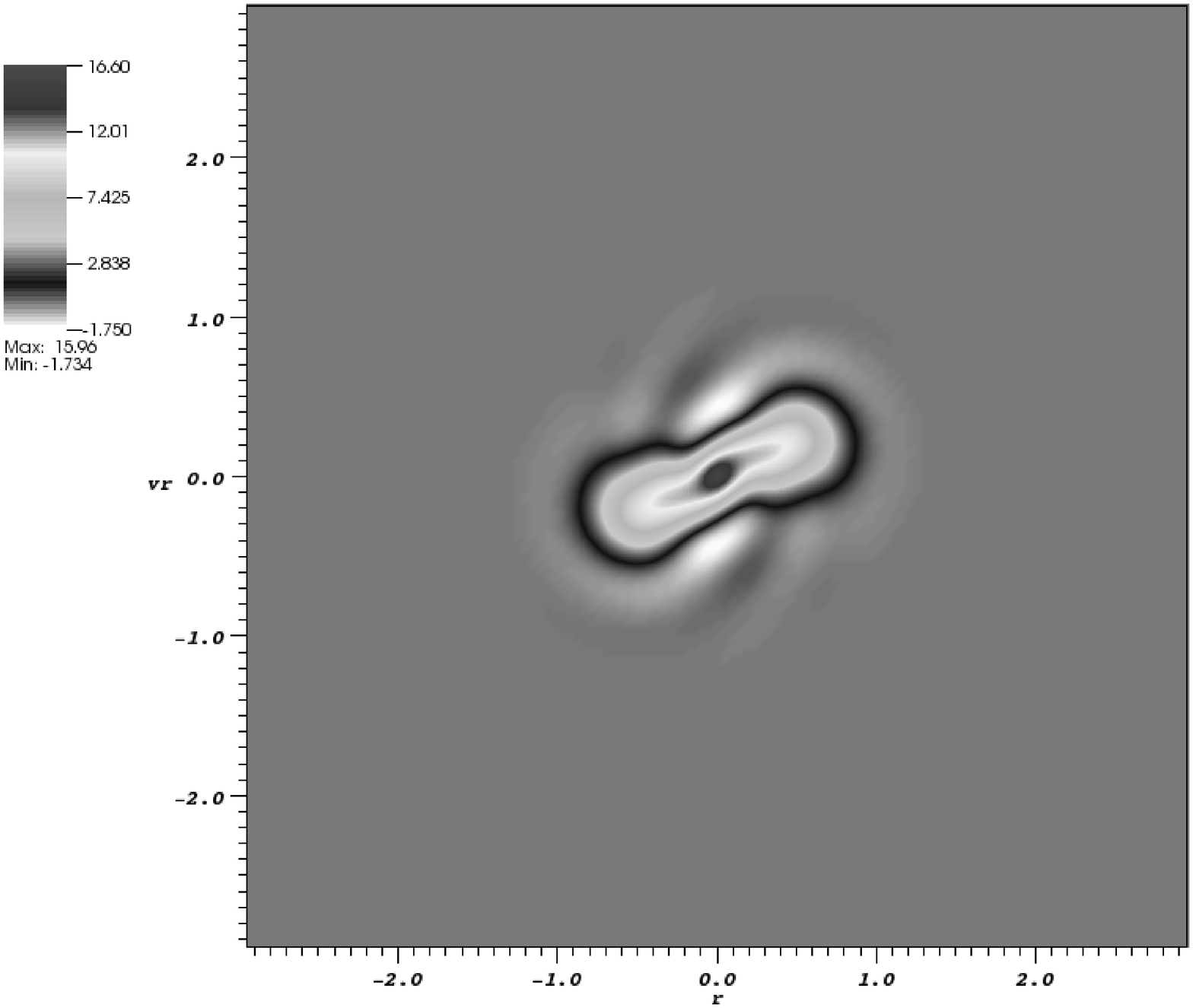} \\
\includegraphics[scale=0.15]{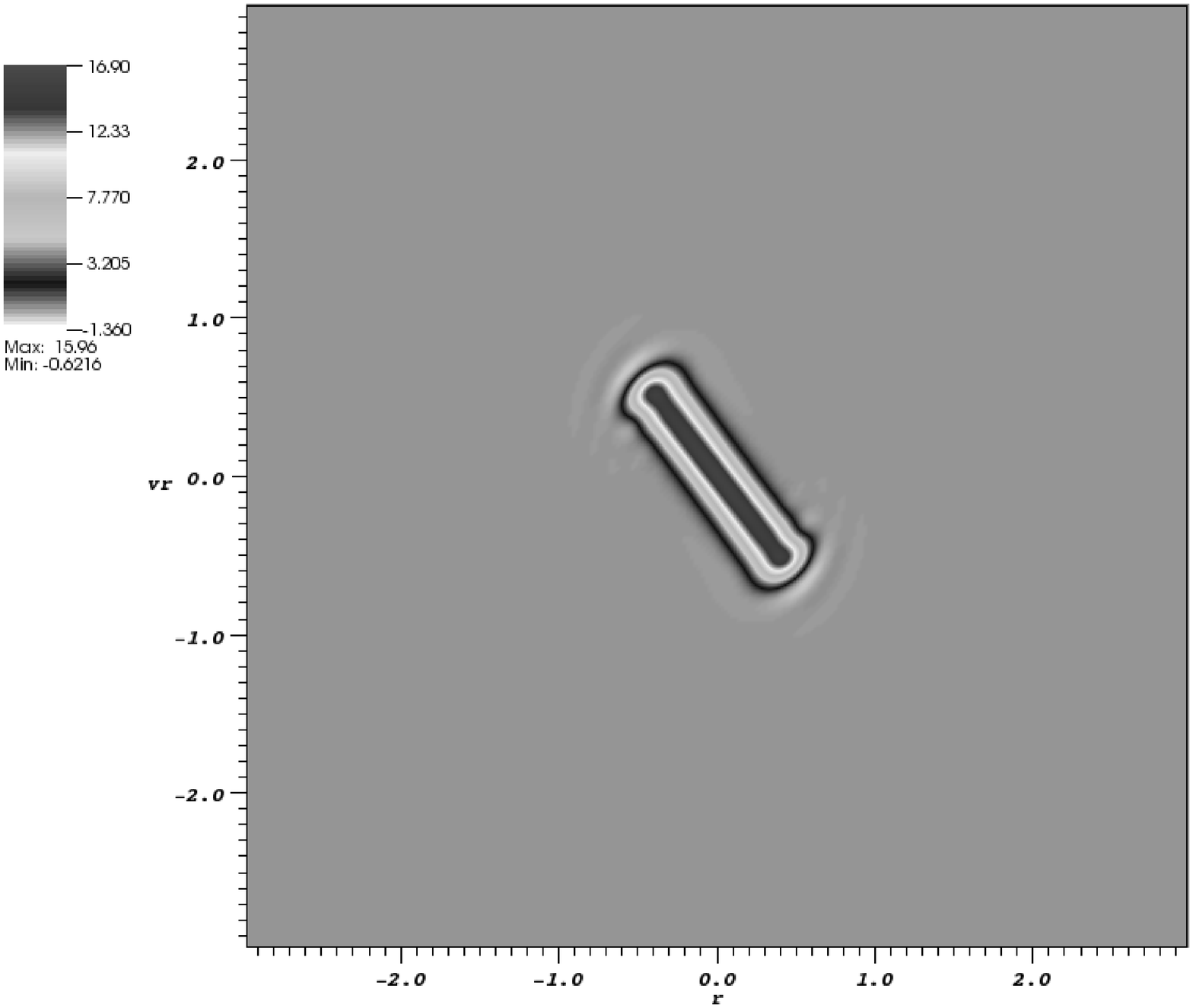} & \includegraphics[scale=0.15]{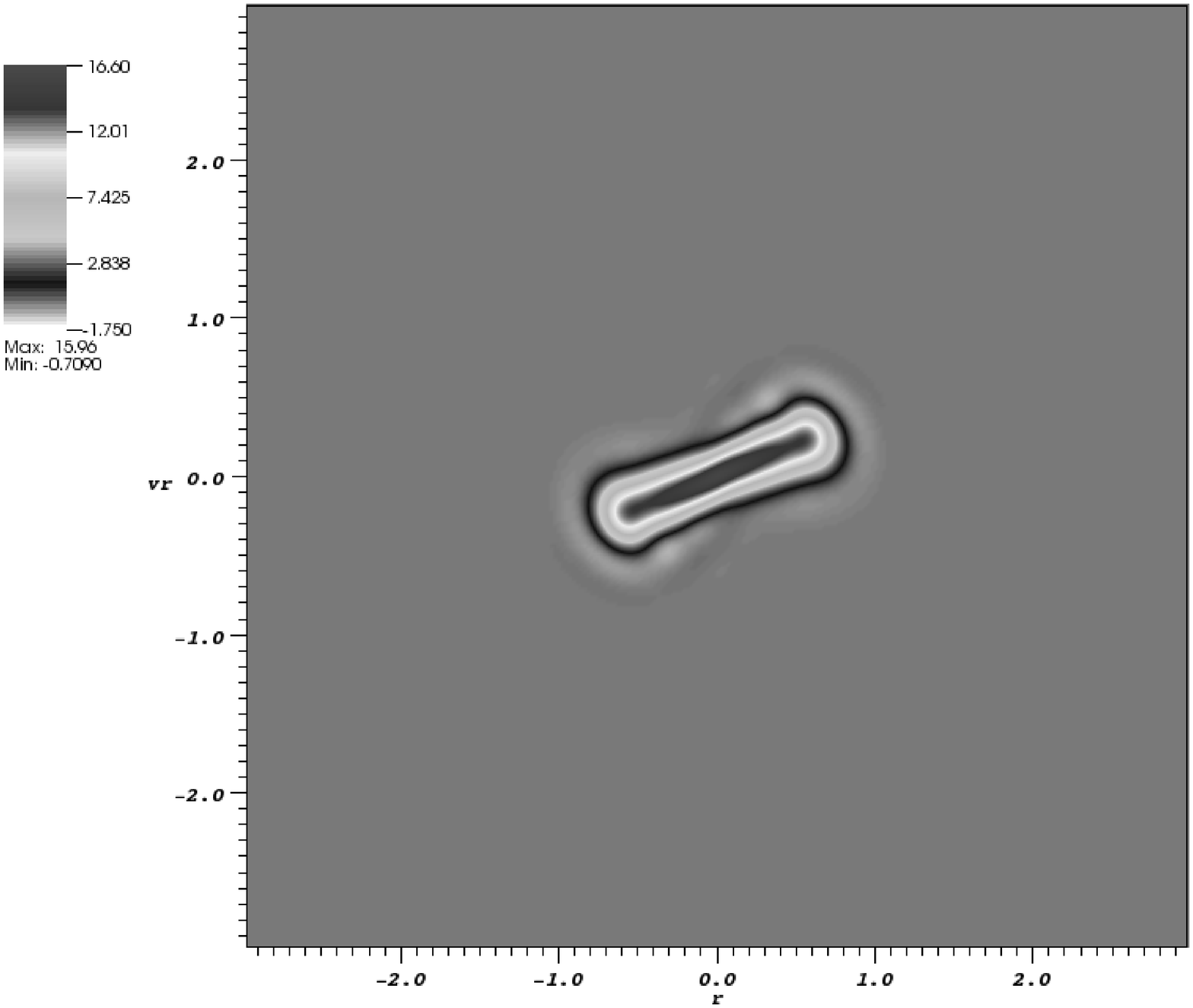} \\
\includegraphics[scale=0.15]{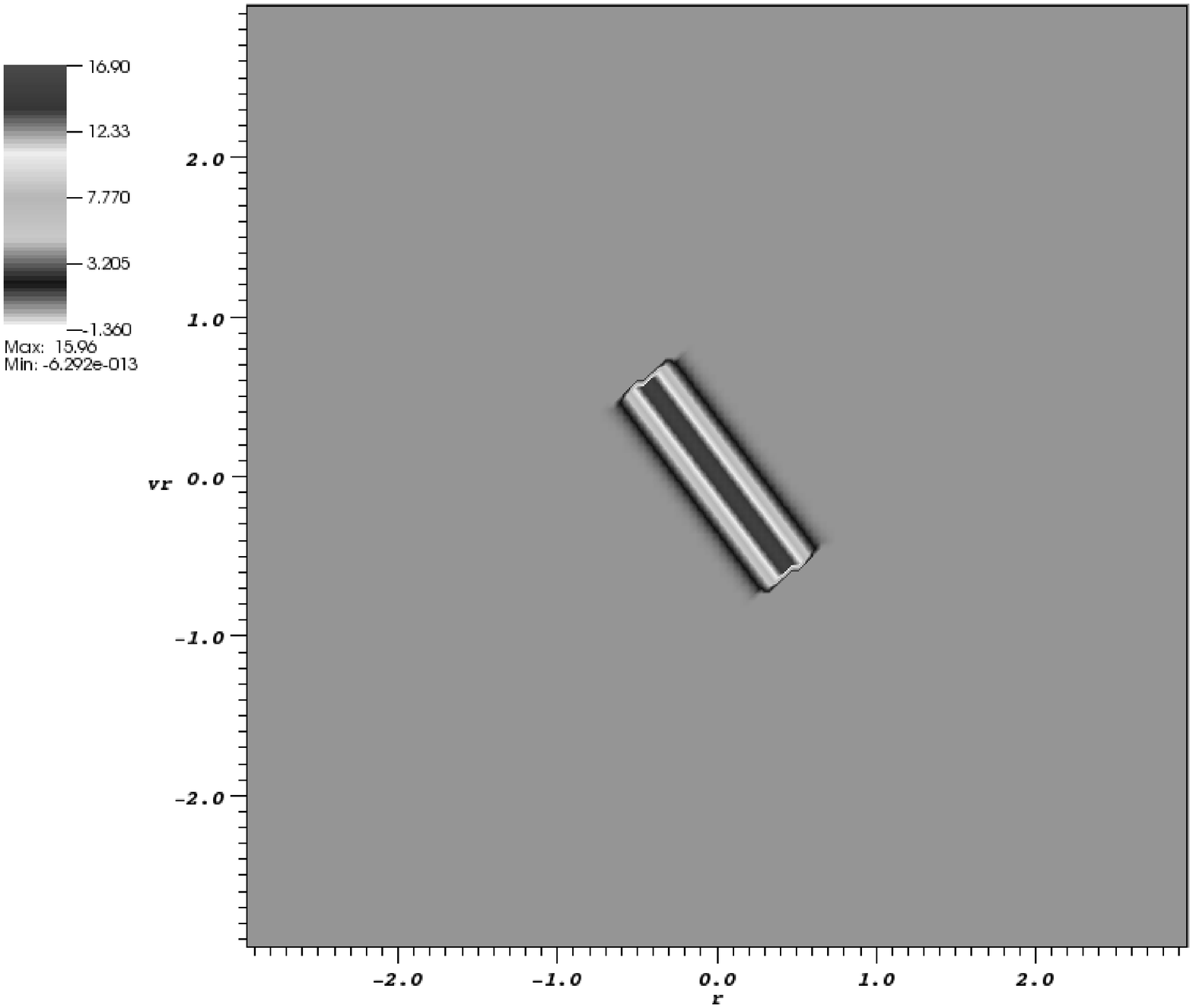} & \includegraphics[scale=0.15]{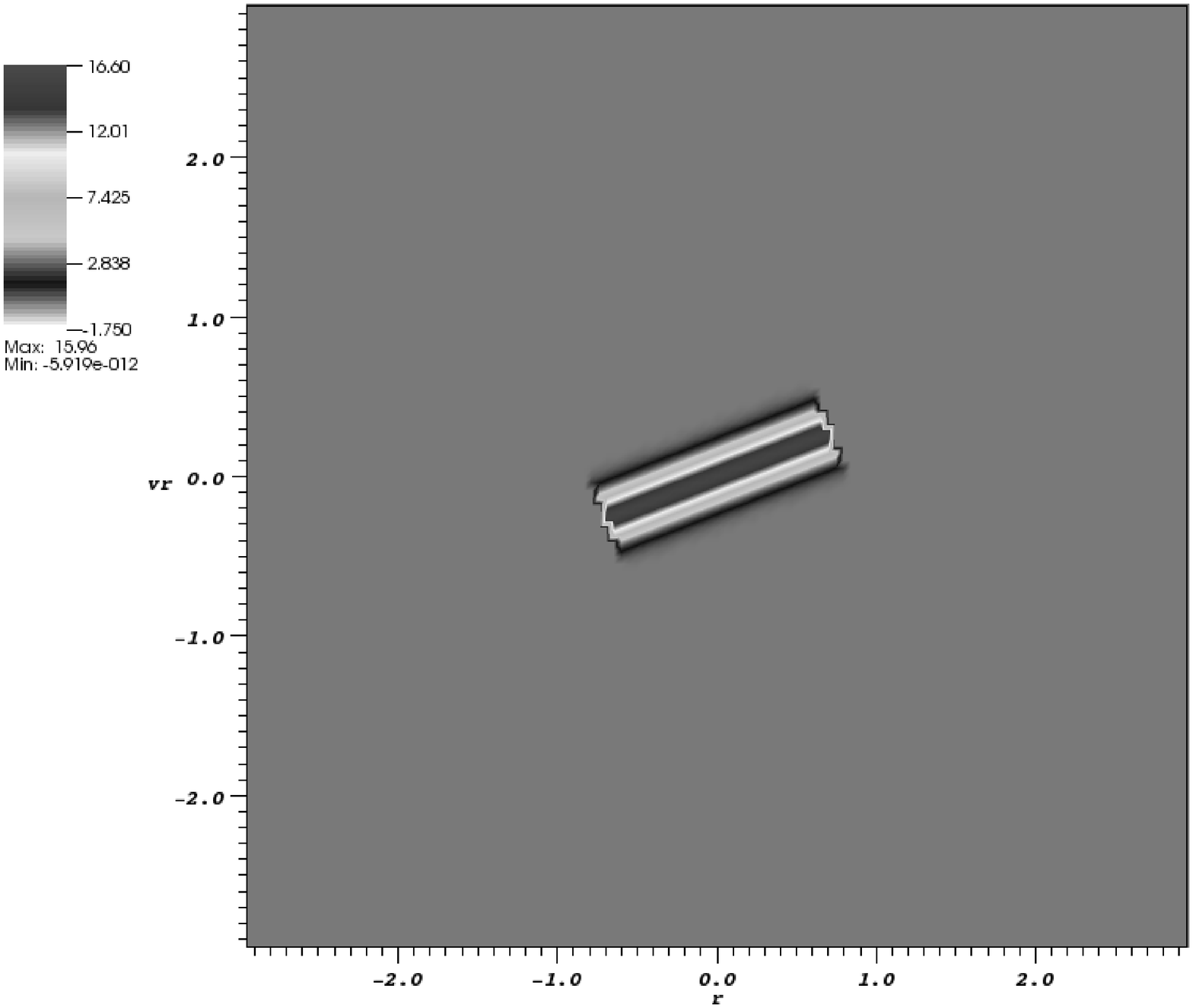} \\
\includegraphics[scale=0.15]{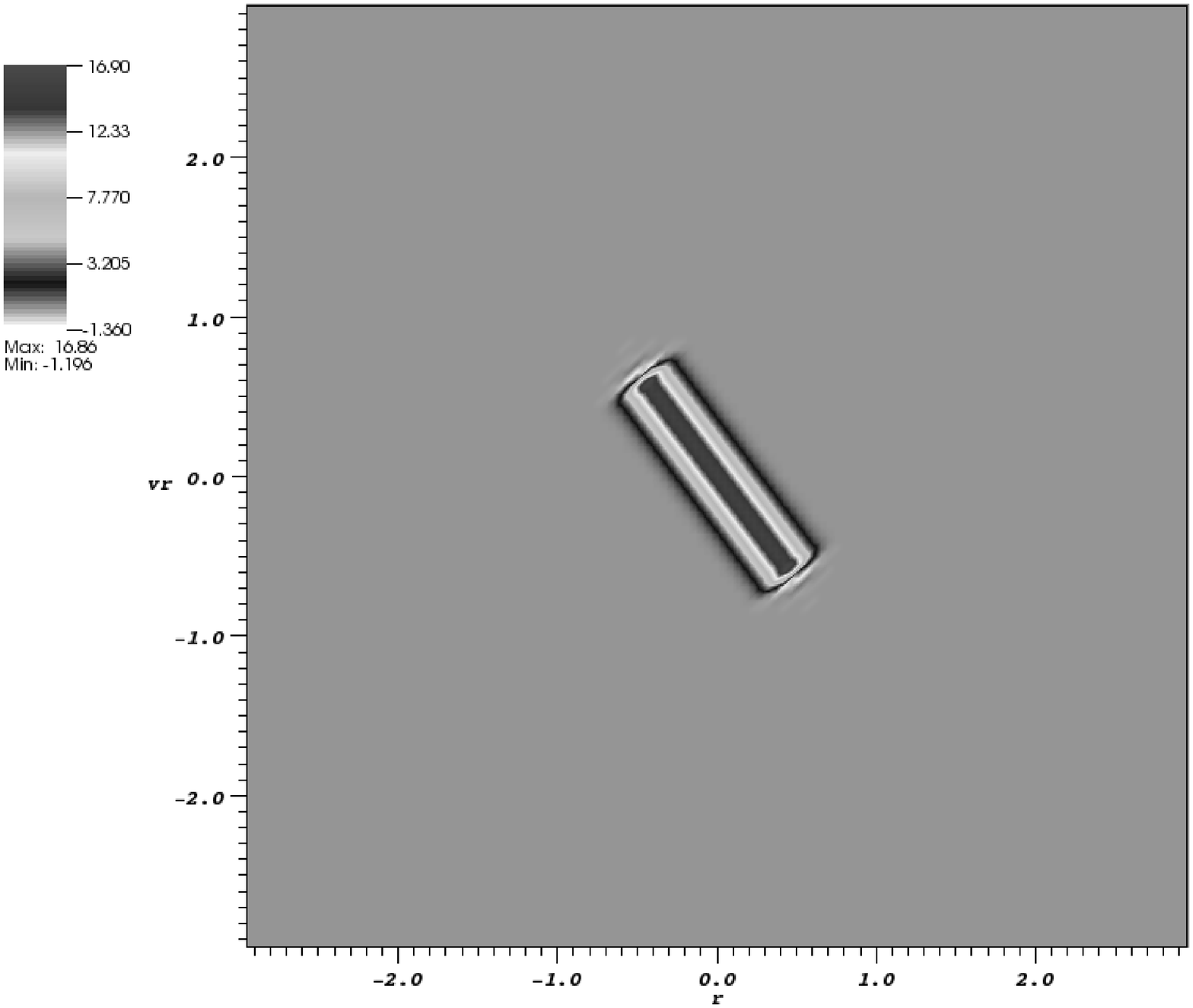} & \includegraphics[scale=0.15]{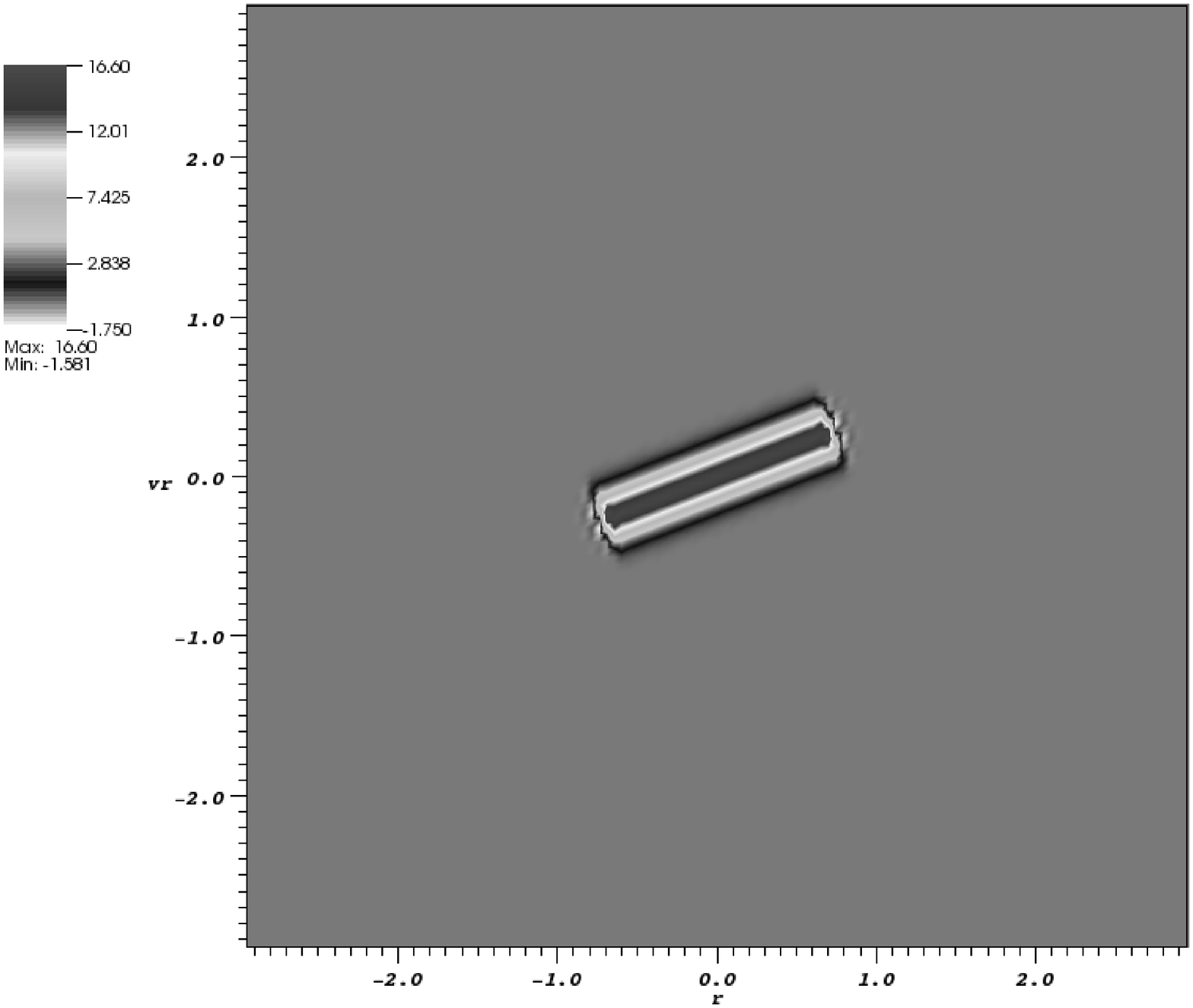} \\
$t = 1.1088$ & $t = 6.468$
\end{tabular}
\caption{Simulations of type (I) (first row), (II) (second row), (III) (third row), and (IV) (fourth row) for a semi-gaussian beam wihout self-consistent electric field and $\omega_{1} = 4\sqrt{2}$, $H_{1}(\tau) = \cos(\tau)$.}
\end{center}

\indent In a first example, we suppose that $\omega_{1} \notin \Q$, so we obtain that $G$ is stationary in $t$, and is equal to $\frac{1}{2\pi} \, f_{0}$. Then, the two-scale simulation reduces itself to the computation of
\begin{equation} \label{analytic_nonresonant}
(r,v_{r},t) \longmapsto f_{0}\Big( \cos\big(\frac{t}{\epsilon}\big) \, r - \sin\big(\frac{t}{\epsilon}\big) \, v_{r} , \sin\big(\frac{t}{\epsilon}\big) \, r + \cos\big(\frac{t}{\epsilon}\big) \, v_{r} \Big) \, ,
\end{equation}
for any function $H_{1}$. It is much simpler than simulating the model (\ref{NH-polar}) with the semi-lagrangian method described in the paragraph 3.2. In Figure 5, we observe such a case, with $\omega_{1} = 4 \sqrt{2}$, $H_{1}(\tau) = \cos(\tau)$, and $f_{0}$ given by
\begin{equation} \label{Semi-gaussian}
f_{0}(r,v_{r},t) = \frac{n_{0}}{\sqrt{2\pi}\,v_{th}} \, \exp\Big( -\frac{v_{r}^{2}}{2\,v_{th}^{2}} \Big) \, \chi_{[-r_{m},r_{m}]}(r) \, ,
\end{equation}
with $\chi_{[-r_{m},r_{m}]}(r) = 1$ if $|r| \leq r_{m}$ and $0$ otherwise. This corresponds to a semi-gaussian beam in particle accelerator physics. In these Figures, we suppose that $r_{m} = 0.75$, $v_{th} = 0.1$, $n_{0} = 4$ and $\epsilon = 10^{-2}$. Furthermore, the simulations (I), (II), (III) and (IV) correspond to
\begin{itemize}
\item simulation (I): we solve the system (\ref{NH-polar}) with a classical semi-lagrangian method, with $P_{r} = P_{v_{r}} = 64$ and $R = v_{R} = 3$,
\item simulation (II): we solve the system (\ref{NH-polar}) with a classical semi-lagrangian method, with $P_{r} = P_{v_{r}} = 128$ and $R = v_{R} = 3$,
\item simulation (III): we solve the system (\ref{H-polar}) with a two-scale semi-lagrangian method on a two-scale mesh, with $P_{r} = P_{v_{r}} = 64$, $P_{\tau} = 16$ and $R = v_{R} = 3$,
\item simulation (IV): we solve the system (\ref{H-polar}) with a two-scale semi-lagrangian method on a uniform mesh, with $P_{q} = P_{u_{r}} = 128$, $P_{\tau} = 16$, and $R = v_{R} = 3$, and we compute the approximation (\ref{approx_tsu}) on a uniform $129 \times 129$ grid in $(r,v_{r})$ on $[-R,R] \times [-v_{R},v_{R}]$.
\end{itemize}

\indent Since we can compute an analytic solution $G$ of (\ref{H-polar}), we can compare the approximations of the solution $f^{\epsilon}$ of (\ref{NH-polar}) to the function defined by (\ref{analytic_nonresonant}). These comparisons are summarized in Figure 6: in this figure, we present the $L^{1}$ norm of the difference between the function (\ref{analytic_nonresonant}) and the approximation $f_{h}$ of $f^{\epsilon}$ obtained with each of the simulations (I), (II), (III), (IV).
\vspace{-0.7cm}
\begin{center}
\includegraphics[scale=0.7]{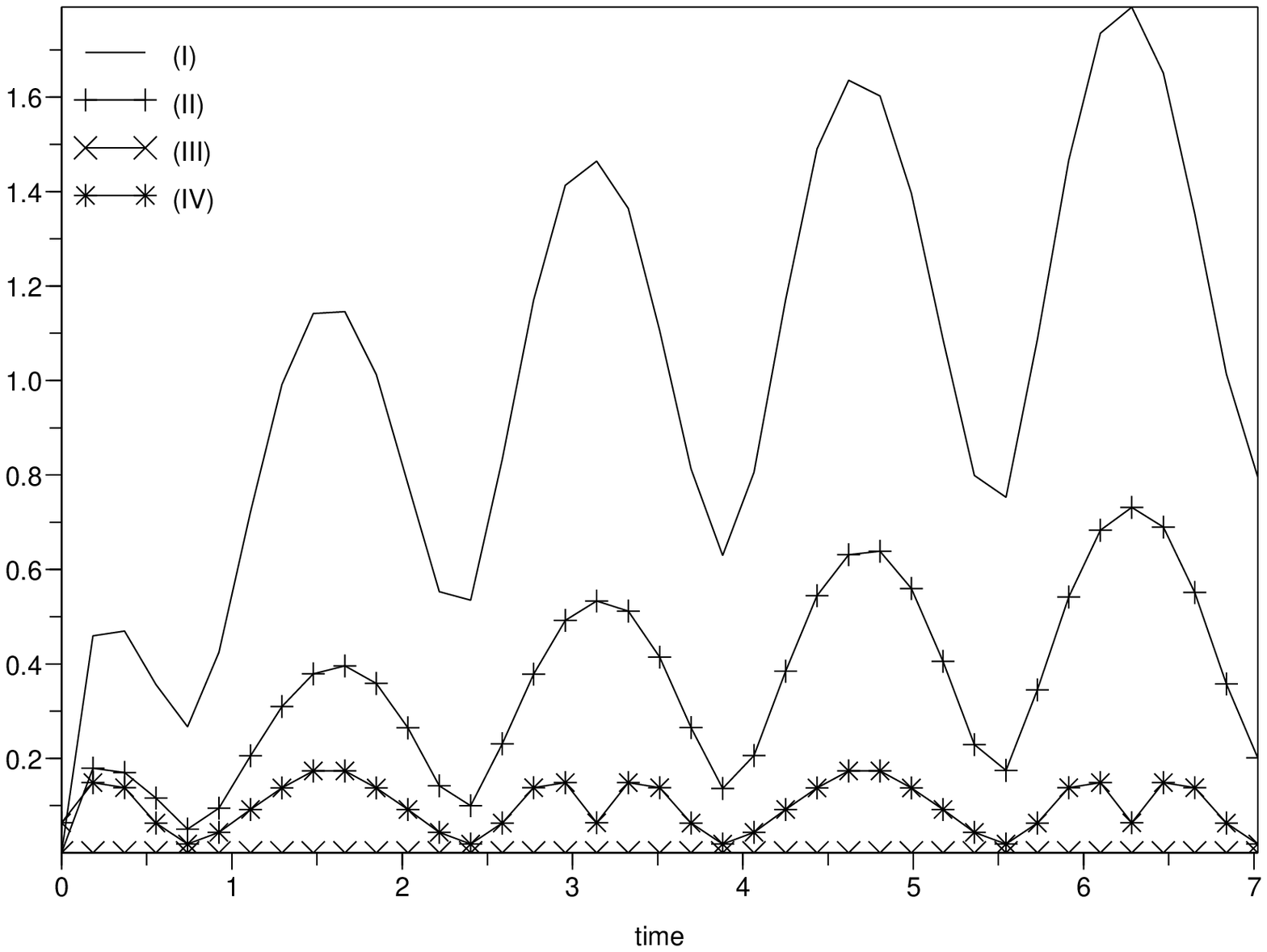}
\vspace{-0.3cm}
\caption{Evolution of the $L^{1}$ norm of the difference between the function (\ref{analytic_nonresonant}) and the approximation $f_{h}$ computed with the simulations (I), (II), (III), and (IV).}
\end{center}
\indent In a second example, we suppose that $\omega \in \N_{\geq \, 2}$. Then, if we assume that $H_{1}(\tau) = \cos^{2}(\tau)$, the model (\ref{H-polar}) reduces itself to
\begin{equation}
\D_{t} G - \frac{u_{r}}{4} \, \D_{q}G + \frac{q}{4} \, \D_{u_{r}}G = 0.
\end{equation}
Knowing the initial data $f_{0}$, we can write
\begin{equation}
G(q,u_{r},t) = \frac{1}{2\pi} \, f_{0} \Big( \cos(\frac{t}{4}\big) \, q + \sin(\frac{t}{4}\big) \, u_{r}, -\sin(\frac{t}{4}\big) \, q + \cos(\frac{t}{4}\big) \, u_{r} \Big) \, ,
\end{equation}
so the two-scale simulation reduces itself to the computation of
\begin{equation}
(r,v_{r},t) \longmapsto f_{0} \Big( \cos(\frac{t}{4}-\frac{t}{\epsilon}\big) \, r + \sin(\frac{t}{4}-\frac{t}{\epsilon}\big) \, v_{r}, -\sin(\frac{t}{4}-\frac{t}{\epsilon}\big) \, r + \cos(\frac{t}{4}-\frac{t}{\epsilon}\big) \, v_{r} \Big) \, .
\end{equation}
In Figure 7, we observe such a case with $\omega_{1} = 2$, $f_{0}$ given by (\ref{Semi-gaussian}) with $r_{m} = 0.75$, $v_{th} = 0.1$ and $\epsilon = 10^{-2}$.

\begin{center}
\begin{tabular}{cc}
\includegraphics[scale=0.15]{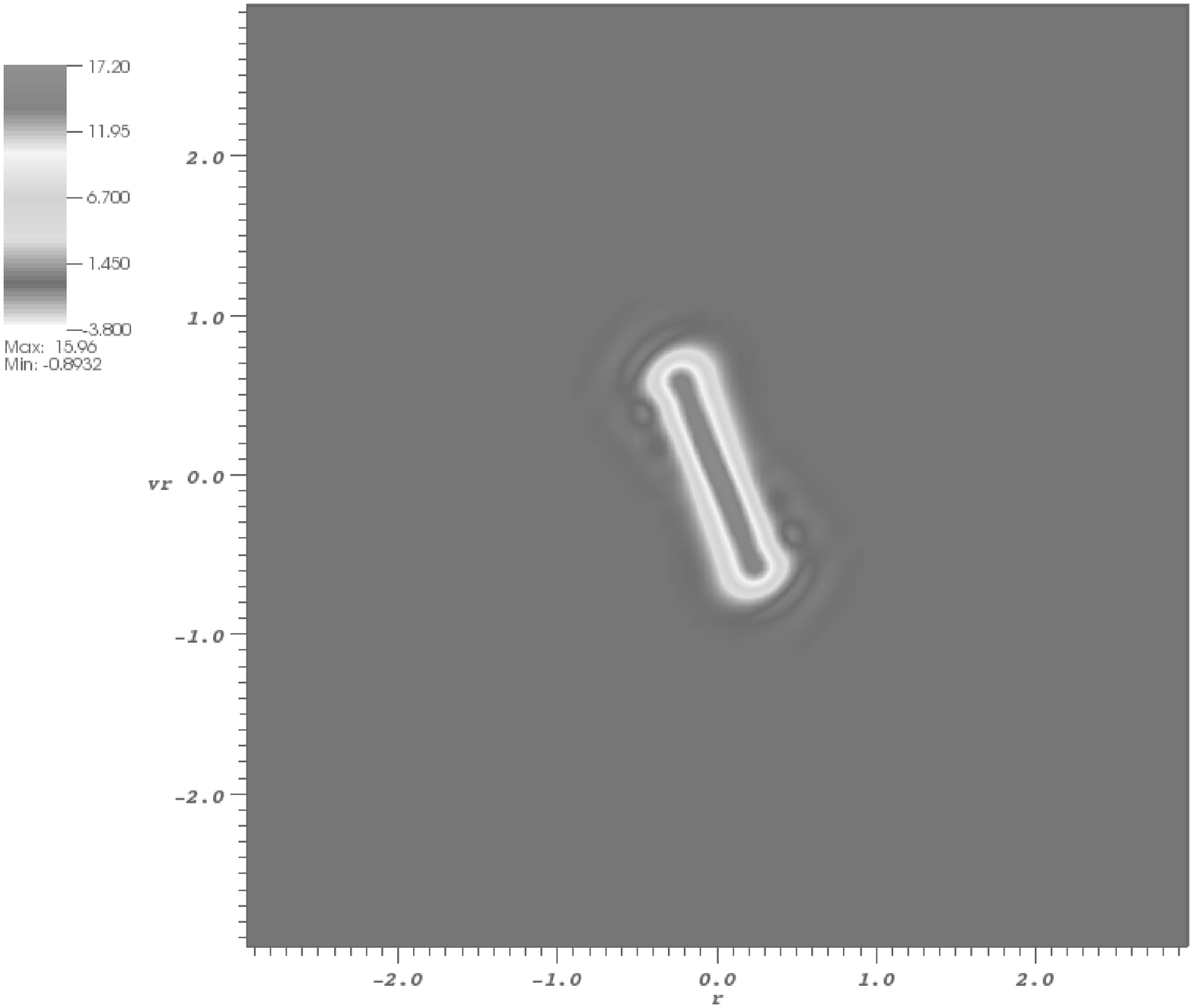} & \includegraphics[scale=0.15]{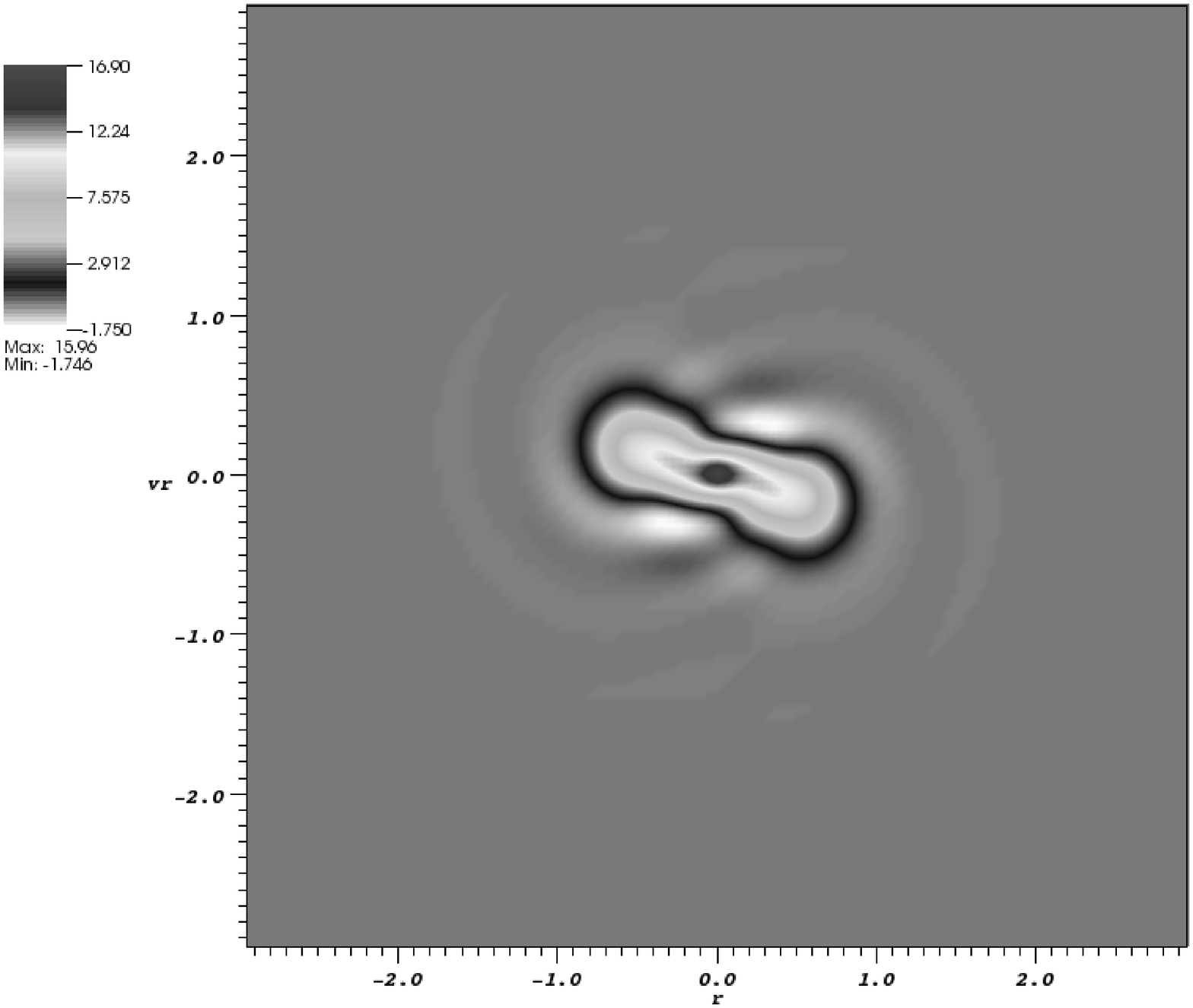} \\
\includegraphics[scale=0.15]{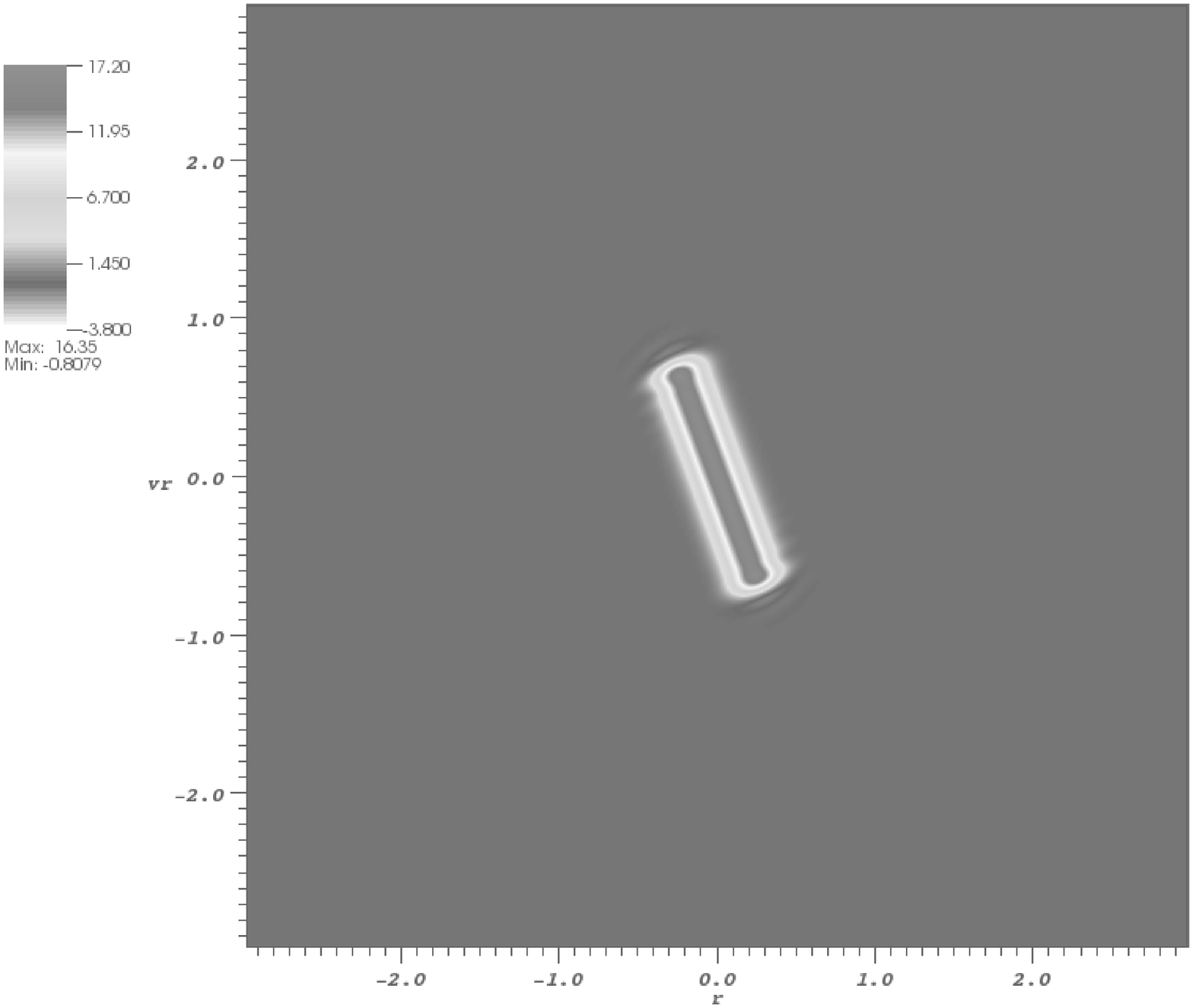} & \includegraphics[scale=0.15]{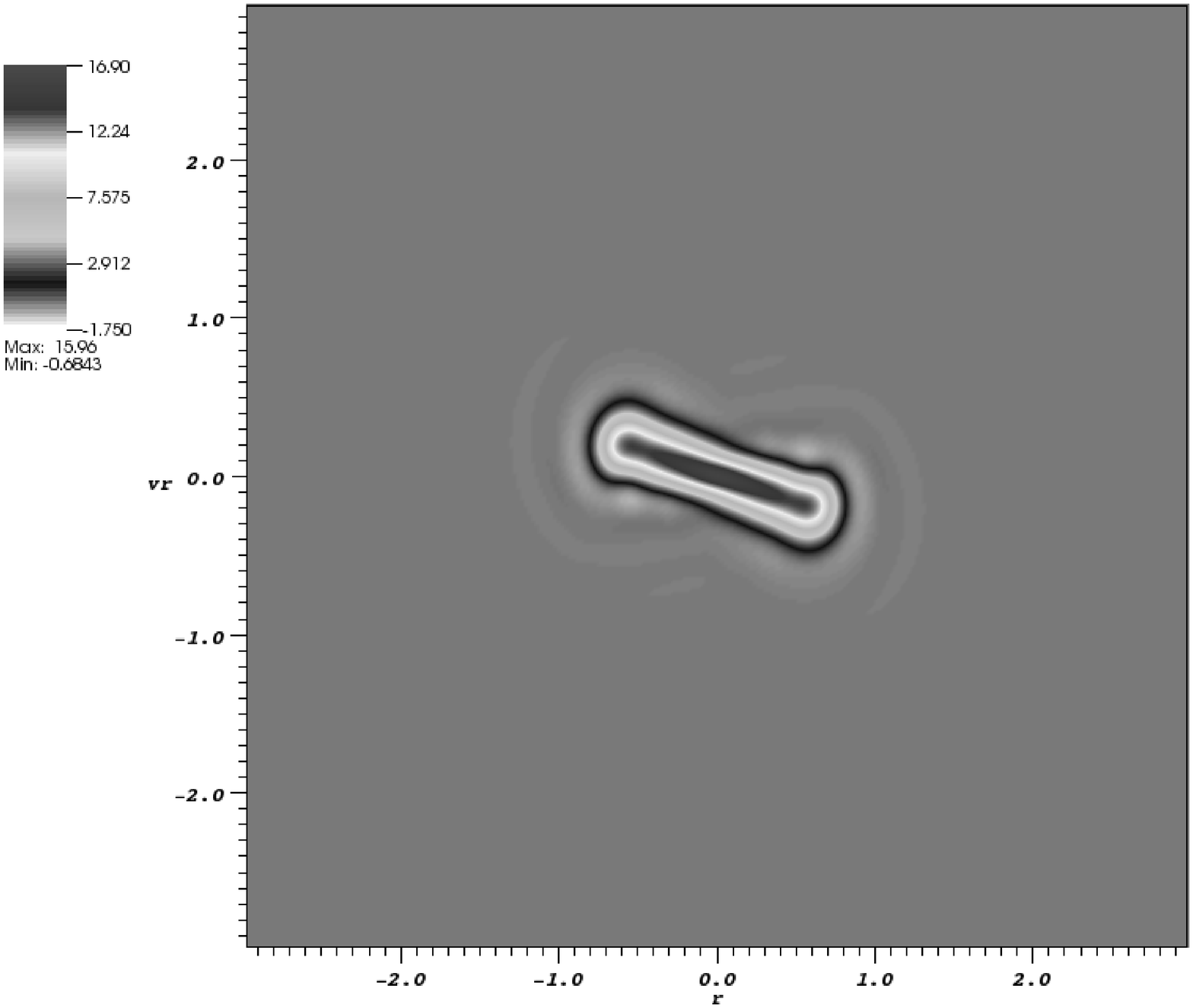} \\
\includegraphics[scale=0.15]{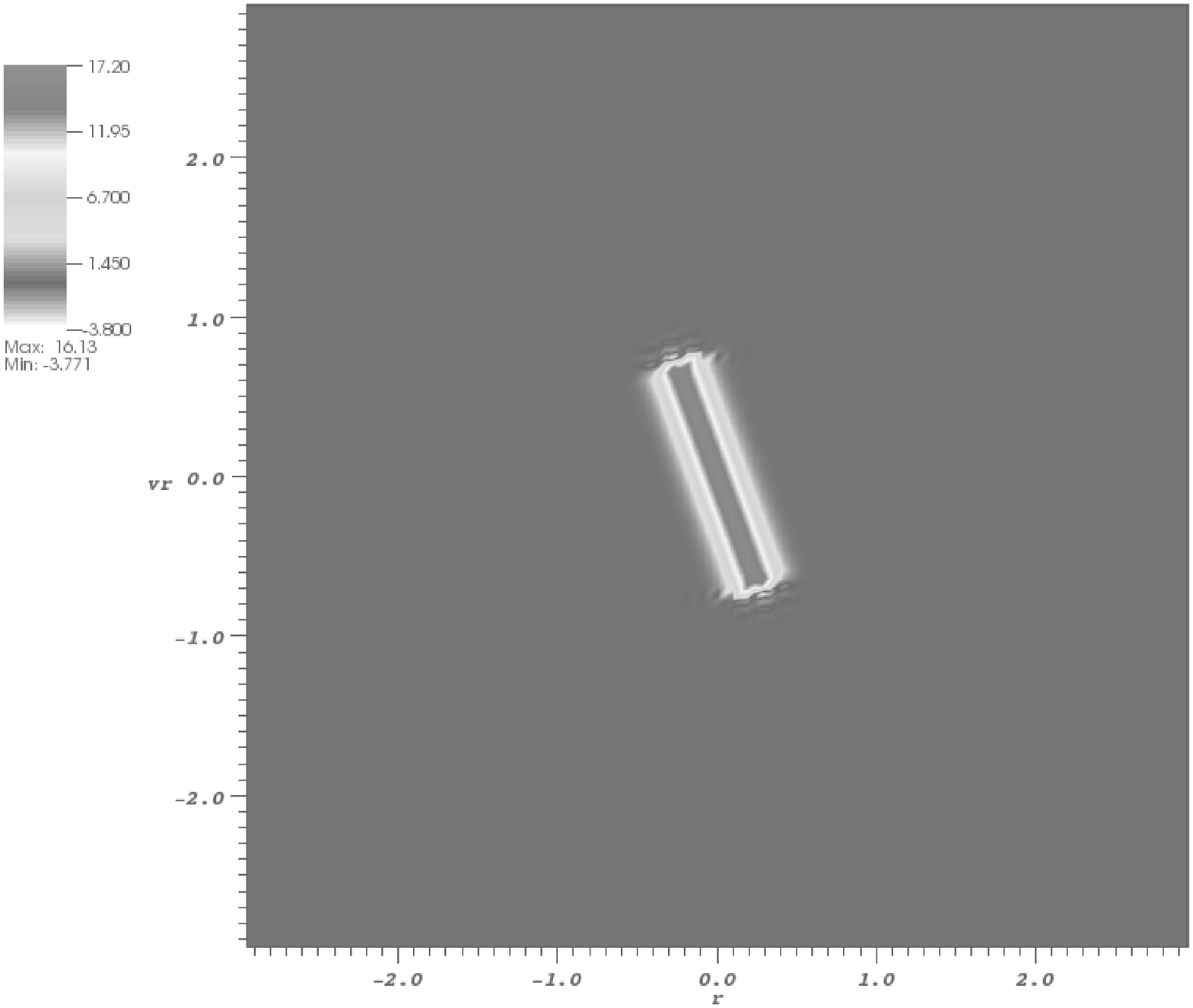} & \includegraphics[scale=0.15]{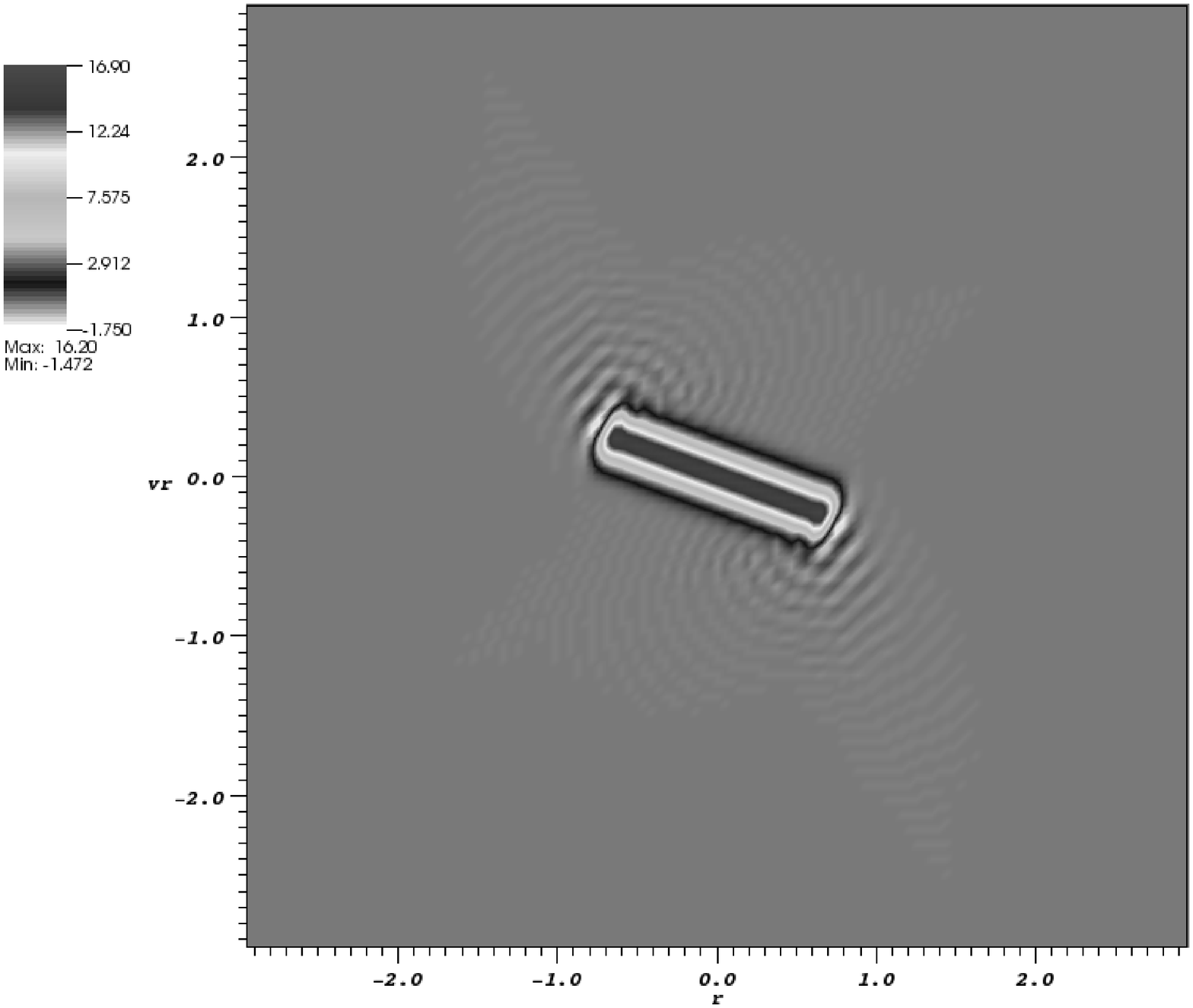} \\
\includegraphics[scale=0.15]{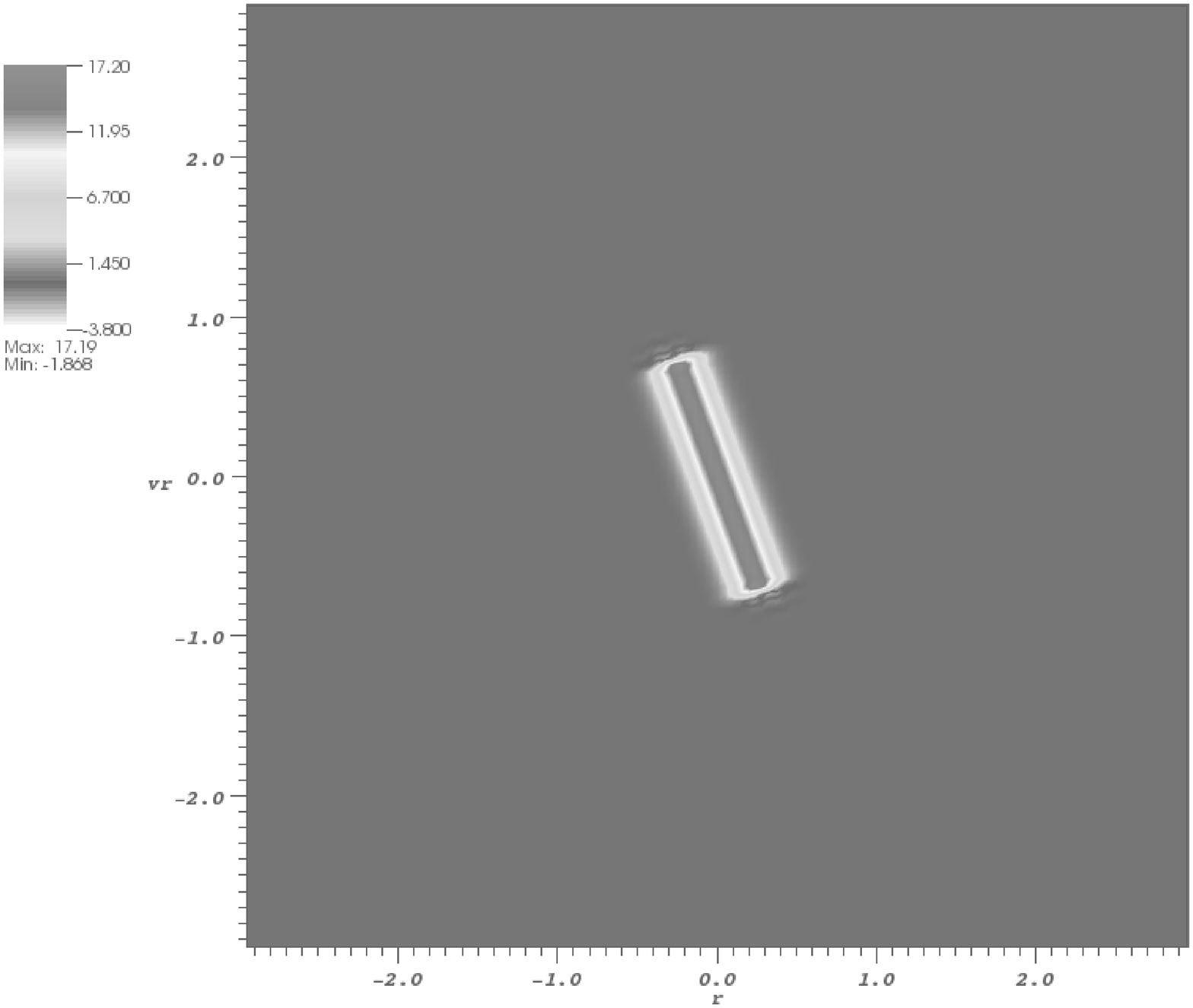} & \includegraphics[scale=0.15]{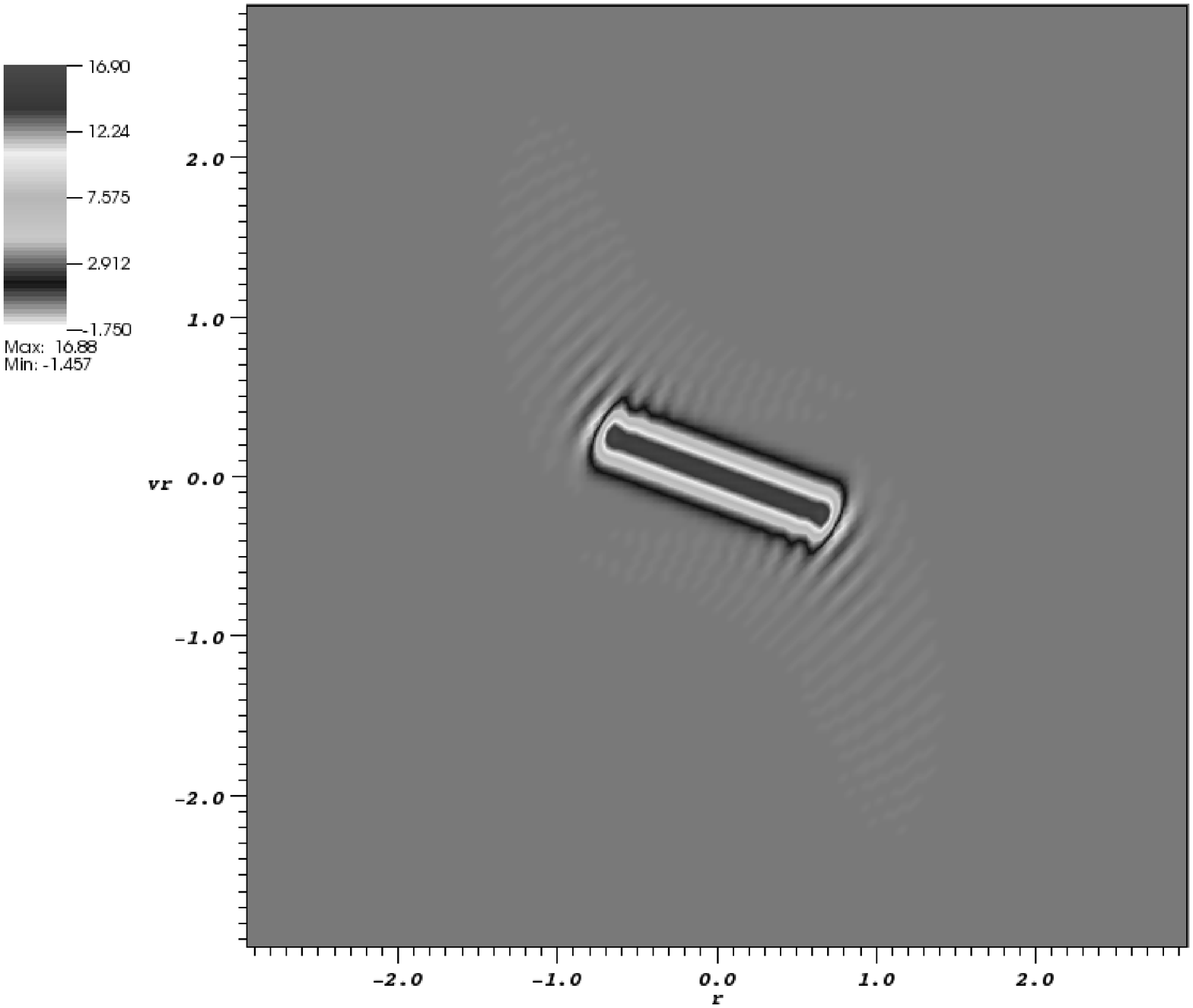} \\
$t = 0.2957$ & $t = 5.9875$
\end{tabular}
\caption{Simulations of type (I) (first row), (II) (second row), (III) (third row), and (IV) (fourth row) for a semi-gaussian beam wihout self-consistent electric field and $\omega_{1} = 2$, $H_{1}(\tau) = \cos^{2}(\tau)$.}
\end{center}

\indent As we can see in Figures 5 and 7, the classical semi-lagrangian method needs a very refined mesh in $r$ and $v_{r}$ directions in order to produce good results. Otherwise, it produces results which do not correspond to physics. We can explain it by the small time step induced by the condition (\ref{CFL-NH}) in order to guarantee the robustness of the method, and then a very high number of interpolations which introduce numerical diffusion. As an example, we have taken a time step $\Delta t_{NH} \approx 1.5 \times 10^{-4}$ for the simulation of type (I) of the second case. \\
\indent On the other hand, both two-scale numerical methods we have described in paragraphs 3.3 and 3.4 produce good results, even if we consider a non-refined mesh in $r$ and $v_{r}$ directions: this phenomenon can be explained by the independence of the time step in $\epsilon$. Consequently, we are allowed to take a time step $\Delta t_{H}$ much larger than we can in the classical semi-lagrangian context, and then, we can significantly reduce the number of interpolations and the numerical diffusion which is induced. For example, for simulations of type (III) and (IV) on the second test case, we have defined $\Delta t_{H}$ with the formula (\ref{def_Delta_t}) where $K = 2$, giving us $\Delta t_{H} \approx 7.4 \times 10^{-3}$. \\
\indent Finally, by observing Figure 6, we remark that the use of the two-scale mesh in the first case reduces significantly the error between the function defined by (\ref{analytic_nonresonant}) and the approximation of $f^{\epsilon}$ given by the discretization of the two-scale model (\ref{H-polar}): indeed, the $L^{1}$ norm of this error is nearby $10^{-6}$ when the approximation of $f^{\epsilon}$ is given by the simulation (III) whereas it oscillates between 0 and 0.75 when the approximation of $f^{\epsilon}$ is given by the simulation (IV). Furthermore, the oscillations of this error for simulations (I), (II), and (IV) are due to the fact that we discretize a non-smooth semi-gaussian distribution the support of which rotates in the phase space under the action of the external electric field and that we do it on a uniform phase space mesh.\\

\indent As a first conclusion, we can say that, on linear cases, the two-scale semi-lagrangian methods we have proposed give much better results than the classical semi-lagrangian method on the same mesh in $r$ and $v_{r}$, which is promising for non-linear cases. Furthermore, as it was announced in the paragraph 3.4, the use of the two-scale mesh reduces significantly the numerical diffusion linked with the global number of interpolations within the simulation.

\subsection{Non-linear cases}

\begin{center}
\begin{tabular}{ccc}
\includegraphics[scale=0.15]{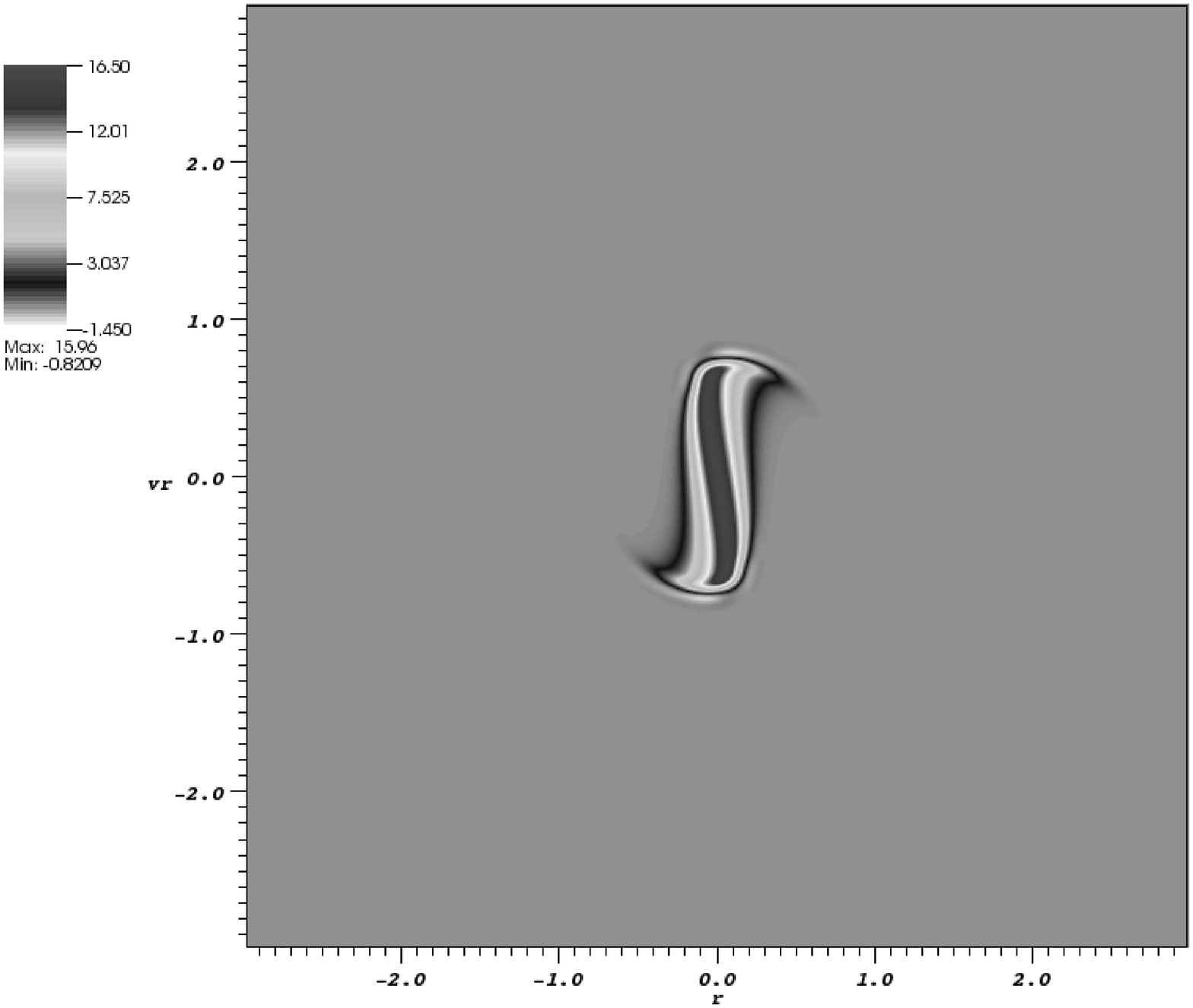} & \includegraphics[scale=0.15]{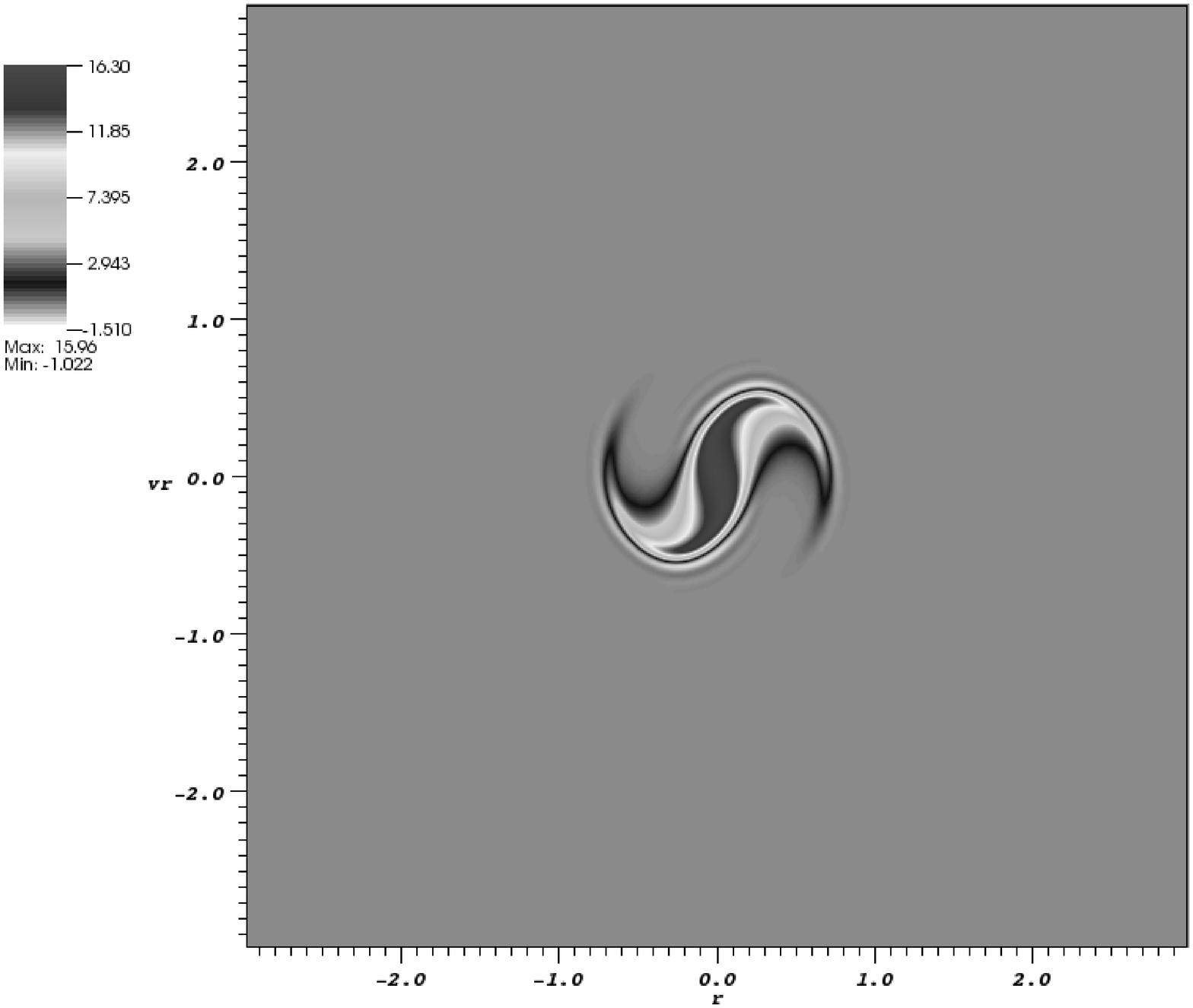} & \includegraphics[scale=0.15]{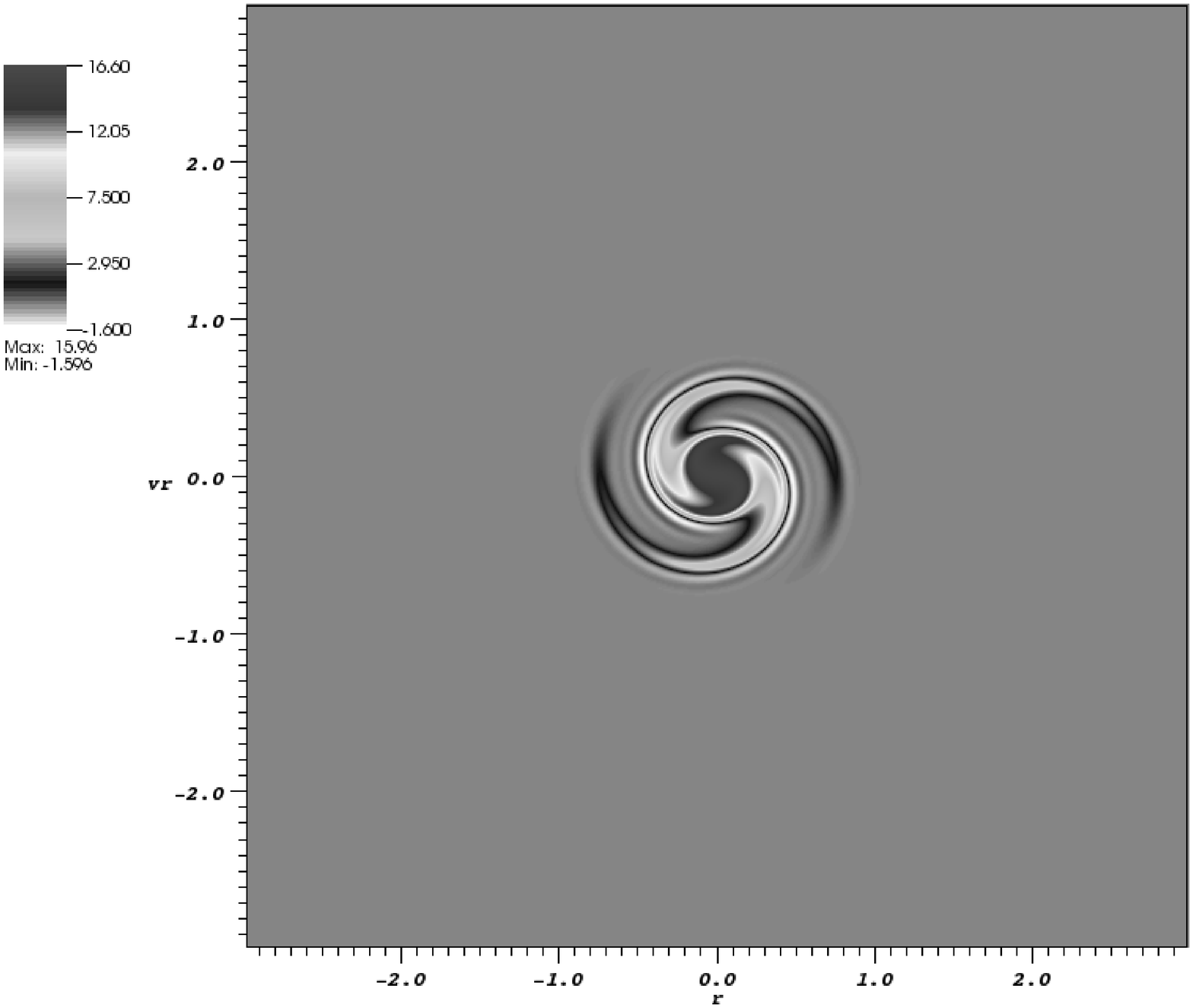} \\
\includegraphics[scale=0.15]{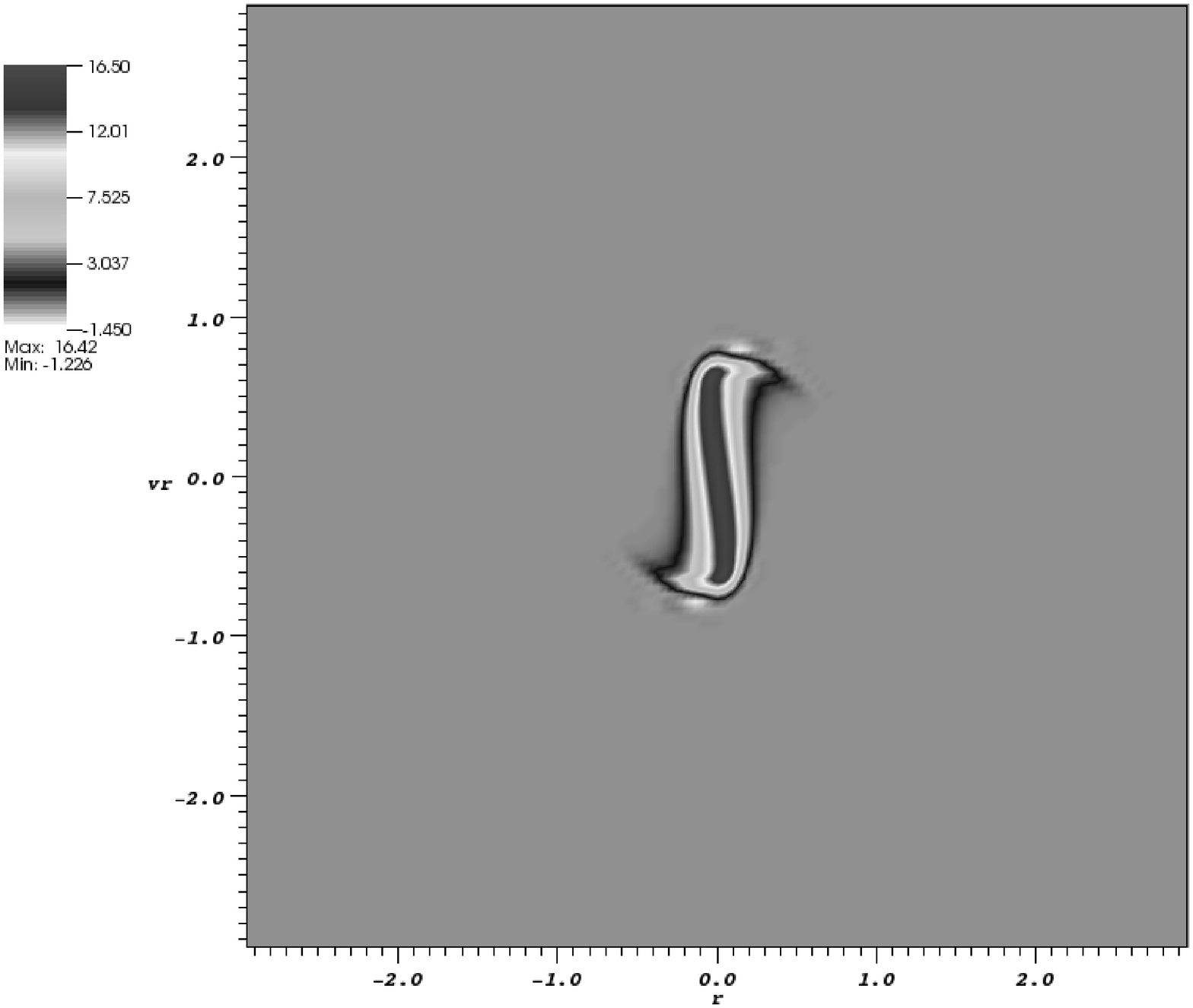} & \includegraphics[scale=0.15]{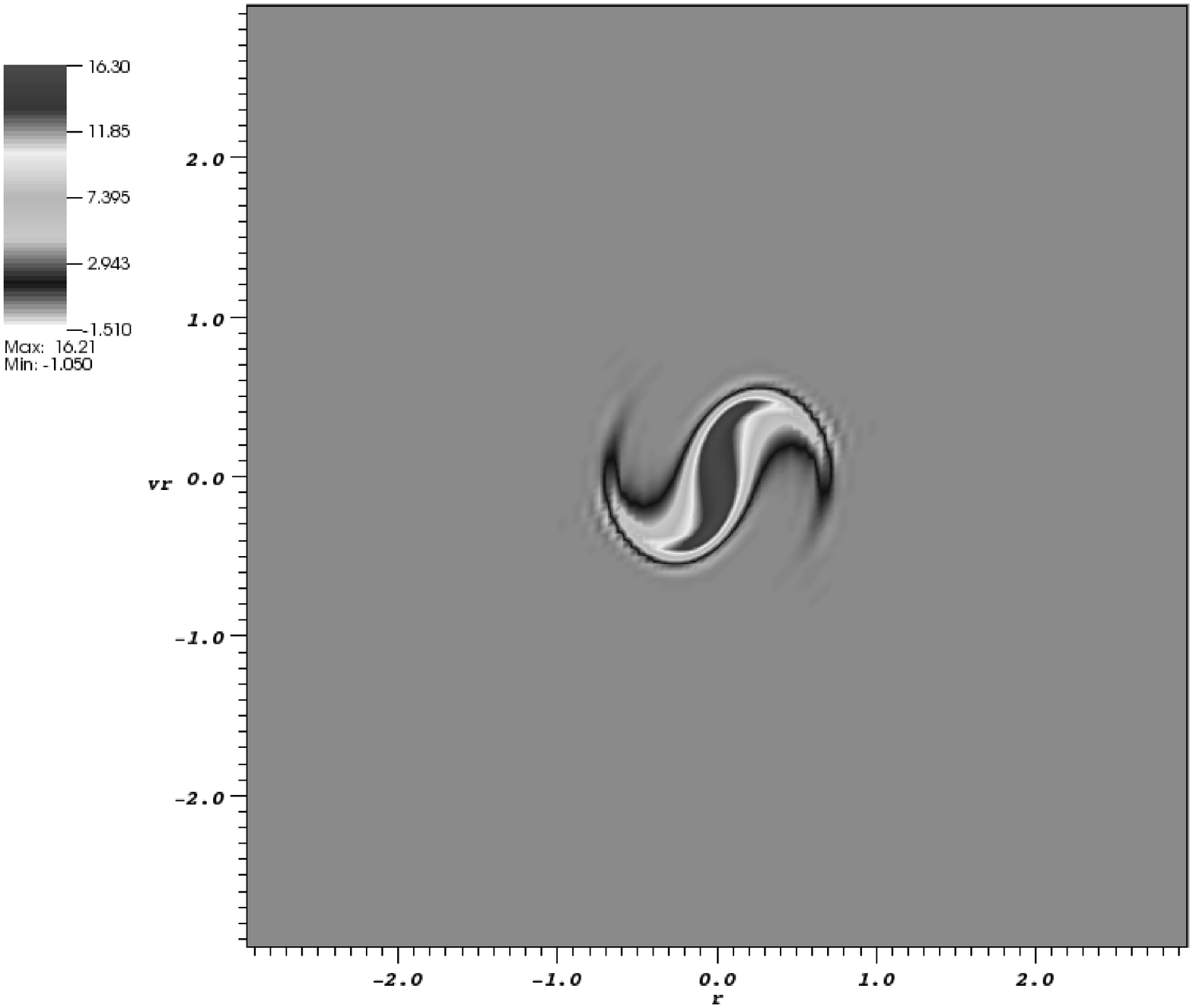} & \includegraphics[scale=0.15]{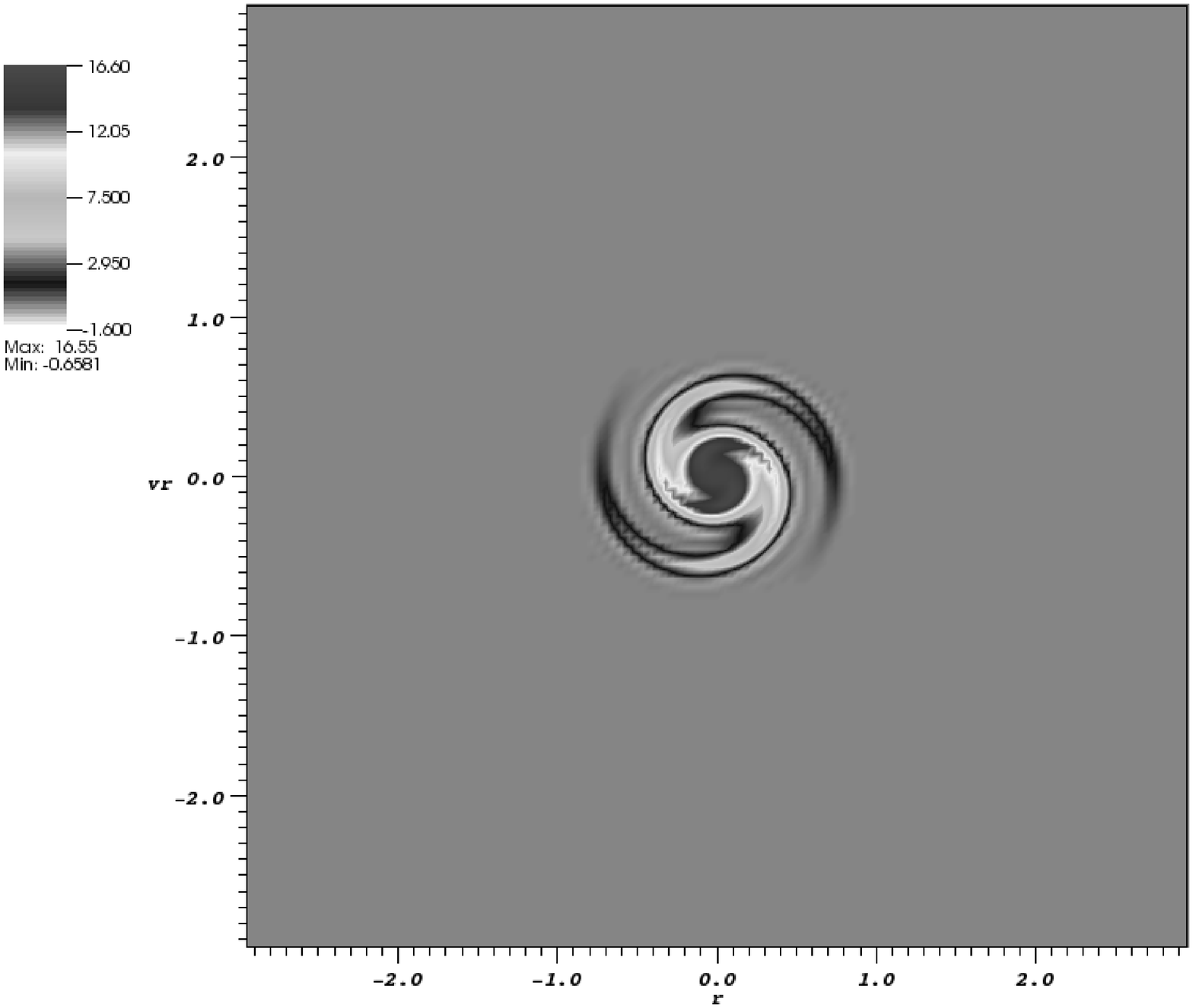} \\
\includegraphics[scale=0.15]{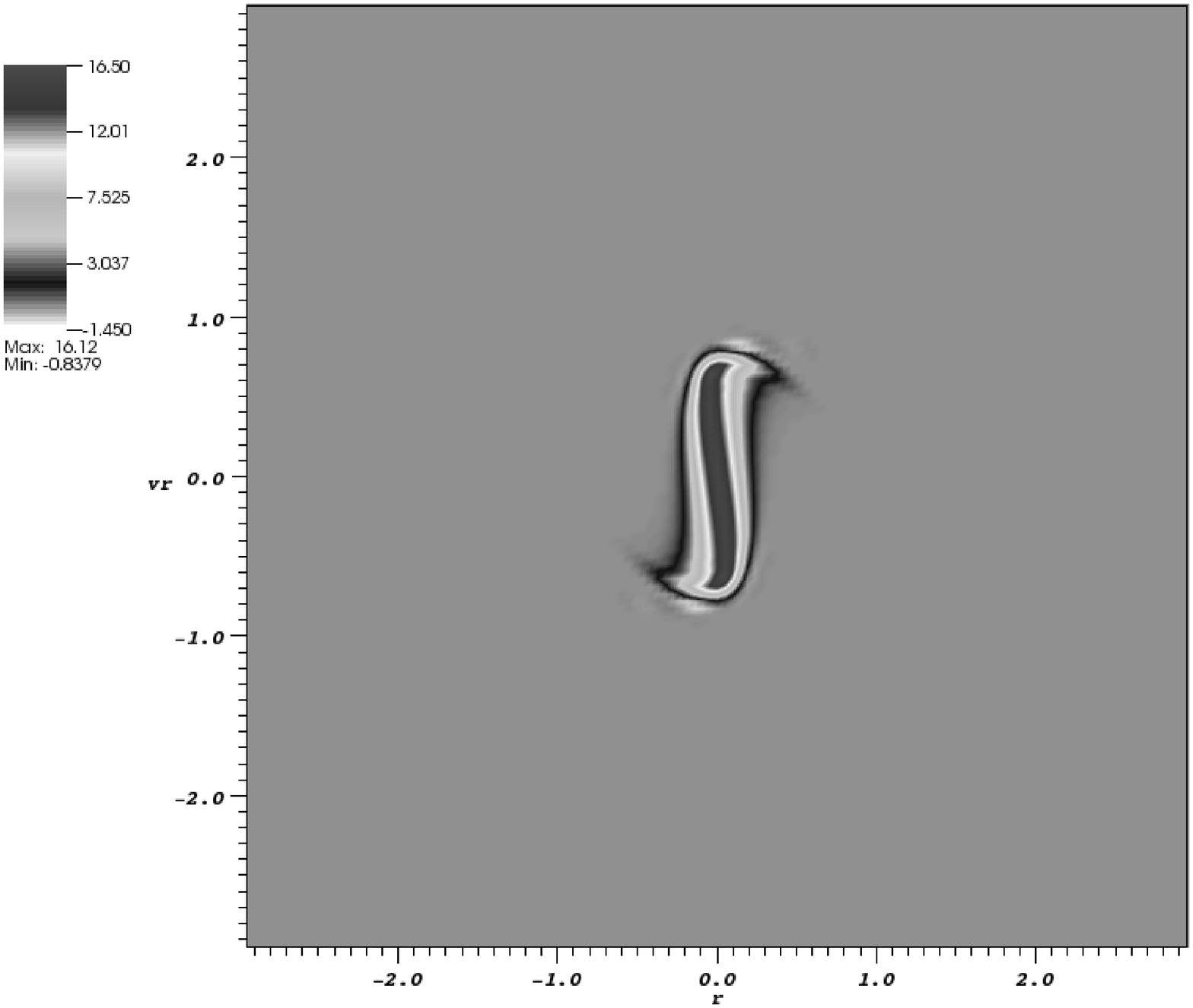} & \includegraphics[scale=0.15]{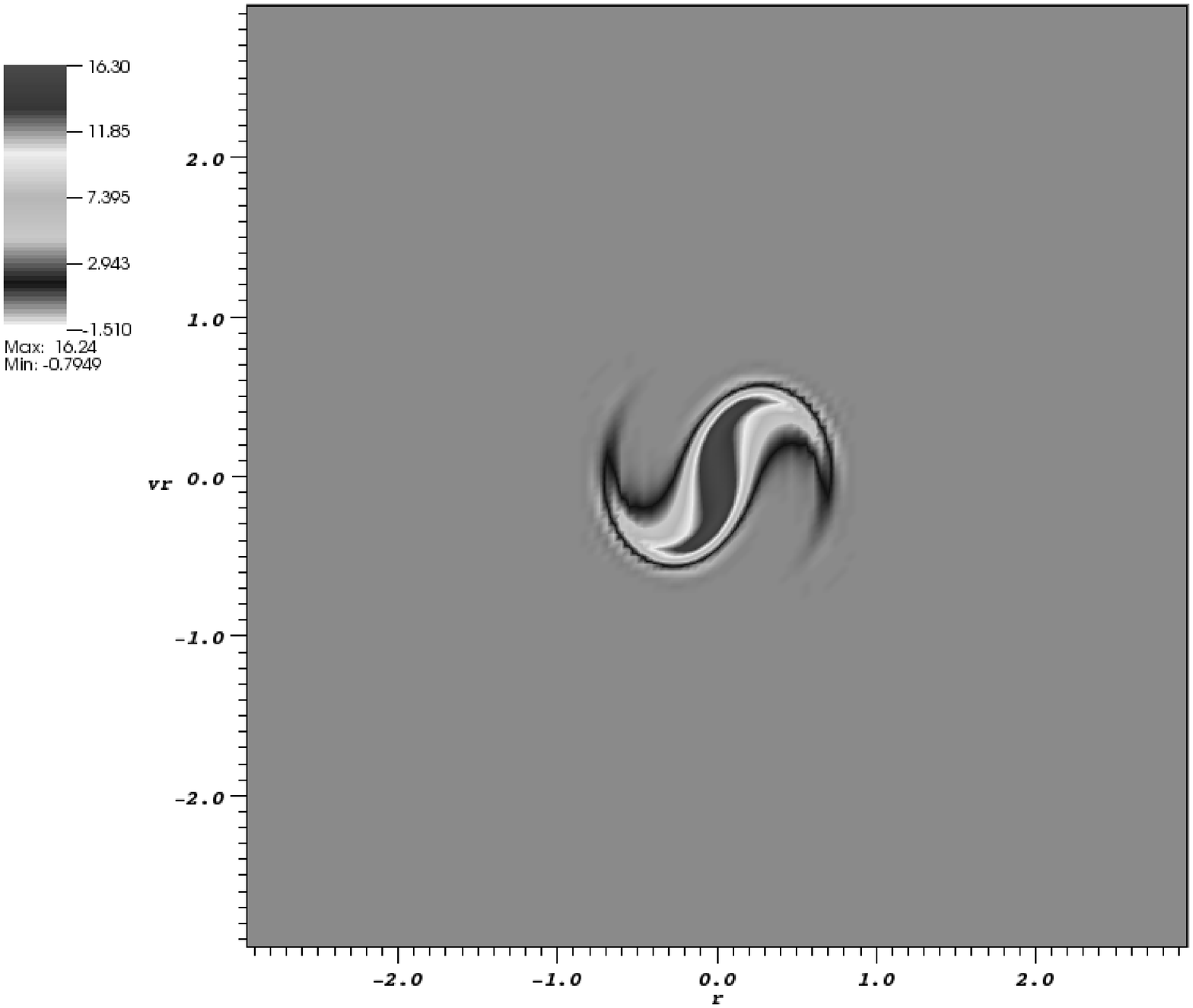} & \includegraphics[scale=0.15]{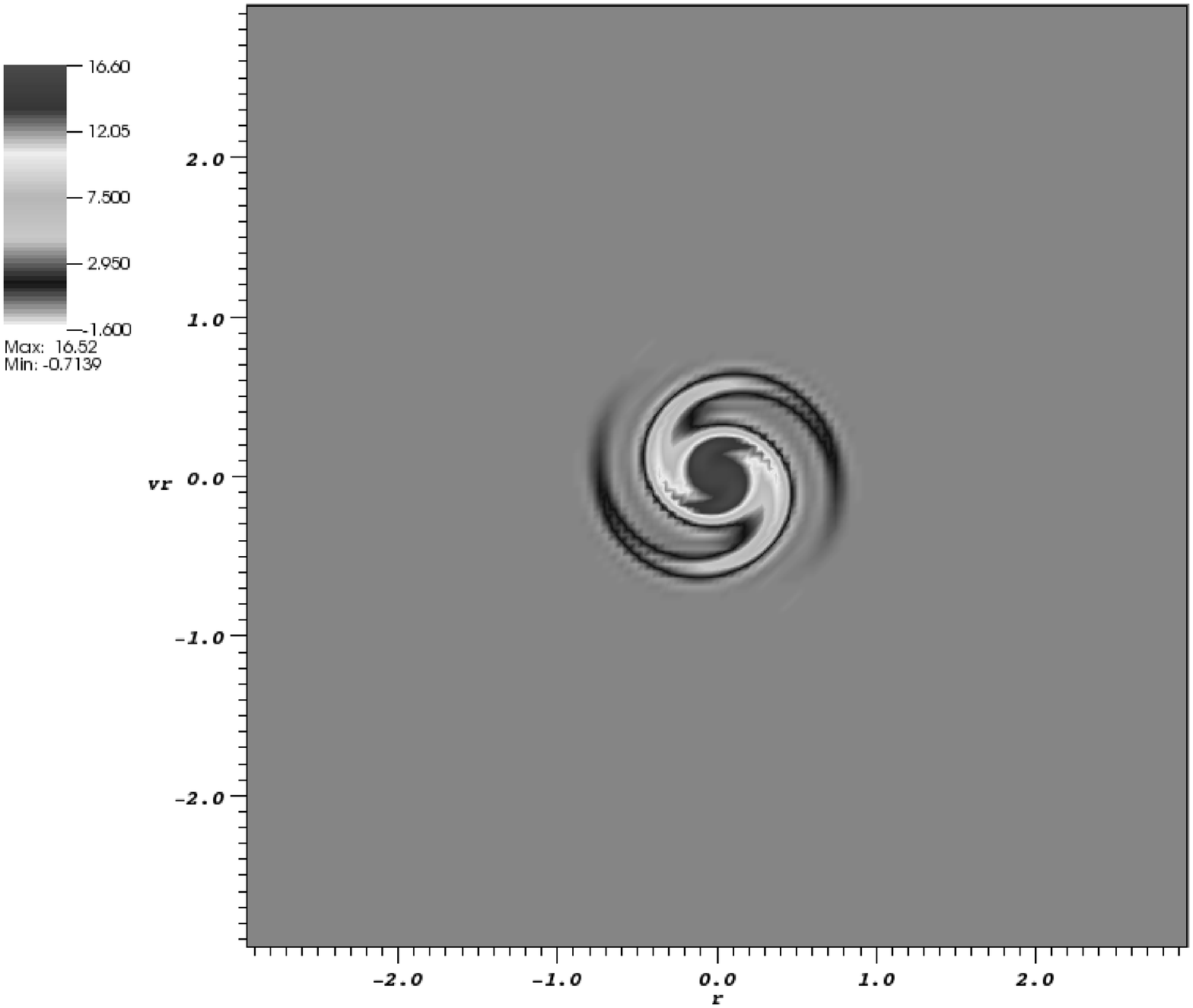} \\
$t = 1.4784$ & $t = 3.234$ & $t = 5.544$
\end{tabular}
\caption{Simulations of type (II') (first row), (III) (second row), and (IV) (third row) for a semi-gaussian beam with $H_{1}(\tau) = \cos(\tau)$ and $\omega_{1} = 4\sqrt{2}$.}
\end{center}

\indent In this paragraph, we do not assume that the self-consistent electric field vanishes. Then, in most of cases, we are not able to find an analytic expression of the solution $G$ of (\ref{H-polar}). So, in order to validate our two-scale methods, we have to compare their results on such a case to the results produced by a classical semi-lagrangian method on (\ref{NH-polar}) with the same initial data, and pay attention to the development of thin structures within the beam. For that, we consider the simulations (I), (II'), (III) and (IV), where the simulation type (II') corresponds to a classical semi-lagrangian method on the system (\ref{NH-polar}) with $P_{r} = P_{v_{r}} = 256$ and $r = v_{r} = 3$. Since the simulation (I) already gives bad results on linear cases, it is not really useful to present its results in terms of quality on non-linear cases. However, we have to pay attention to its CPU time cost in order to compare it to the other simulation types ones\footnote{All the simulations have been conducted on a Sun Fire X4600 server with an AMD Opteron 8220 processor (Dual Core 3000 Mhz) under SunOS 5.10 system.}. \\

\indent In a first case, we suppose that $H_{1}(\tau) = \cos(\tau)$, $\omega_{1} = 4\sqrt{2}$, $\epsilon = 10^{-2}$, and $f_{0}$ is given by (\ref{Semi-gaussian}) with $r_{m} = 0.75$ and $v_{th} = 0.1$. In order to guarantee the robustness of the two-scale schemes, we suppose that the time step $\Delta t_{H}$ is computed by using the the formula (\ref{def_Delta_t}) with $K = 5$, giving us $\Delta t_{H} \approx 0.0185$, whereas the time step $\Delta t_{NH}$ for the simulations (I) and (II') is of the form $\frac{\Delta t_{H}}{N}$ with $N$ large enough in order to verify the condition (\ref{CFL-NH}). In Figure 8, we observe the results obtained with simulations of type (II'), (III), and (IV), and we can find some results in terms of CPU time in the Table 1. \\
\indent In a second case, we suppose that $H_{1}(\tau) = \cos^{2}(\tau)$, $\omega_{1} = 2$, $\epsilon = 10^{-2}$, and $f_{0}$ is given by (\ref{Semi-gaussian}) with $r_{m} = 0.75$ and $v_{th} = 0.1$. In order to guarantee the robustness of the two-scale schemes, we suppose that $\Delta t_{H}$ is defined by (\ref{def_Delta_t}) with $K = 2$, which gives us $\Delta t_{H} \approx 7.392 \times 10^{-3}$, and $\Delta t_{NH} = \frac{\Delta t_{H}}{N}$ with $N$ large enough. In Figure 9, we observe the results obtained with simulations of type (II'), (III), and (IV) and CPU times costs for each simulation can be found in the Table 1. \\

\begin{center}
\begin{tabular}{ccc}
\includegraphics[scale=0.15]{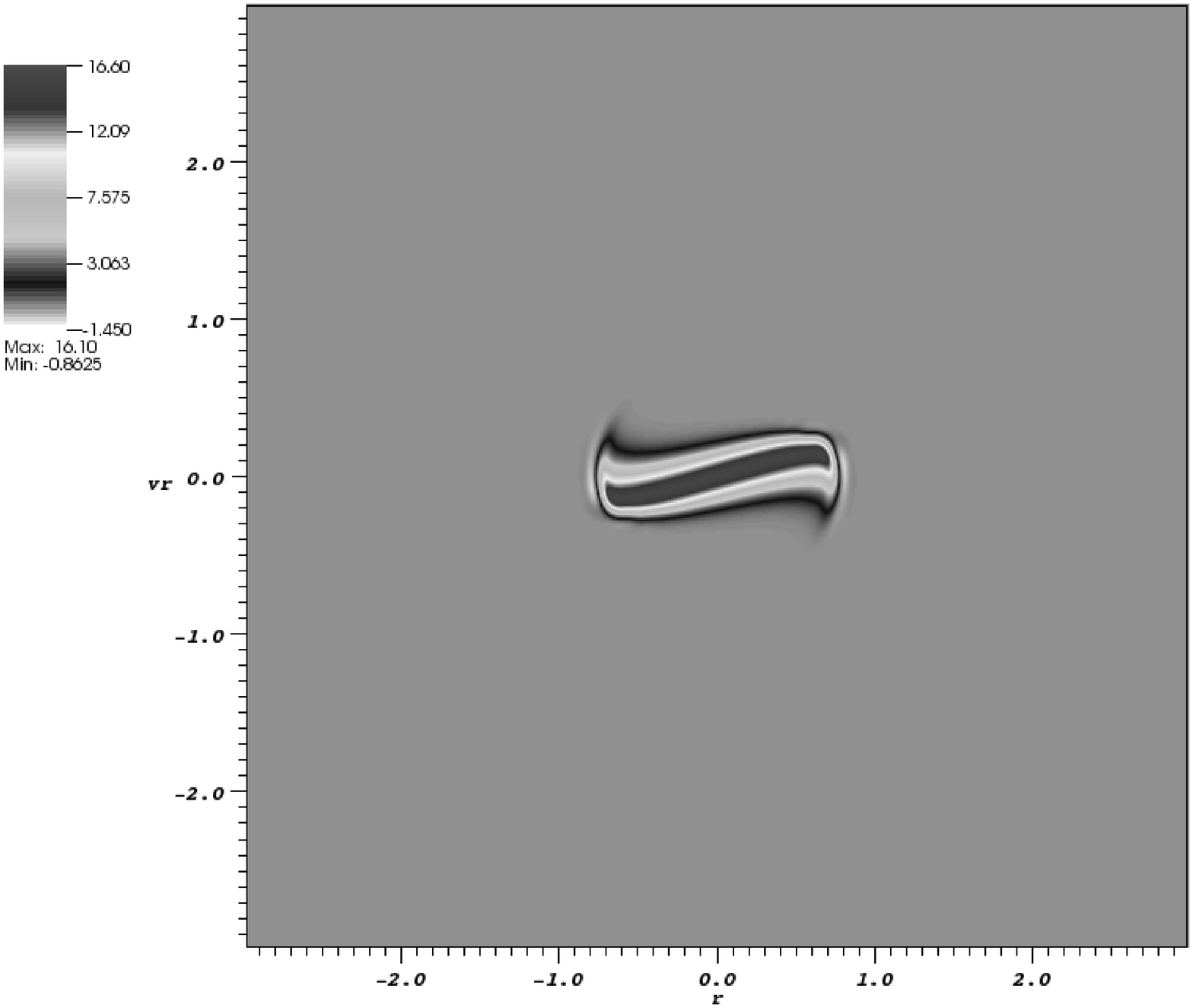} & \includegraphics[scale=0.15]{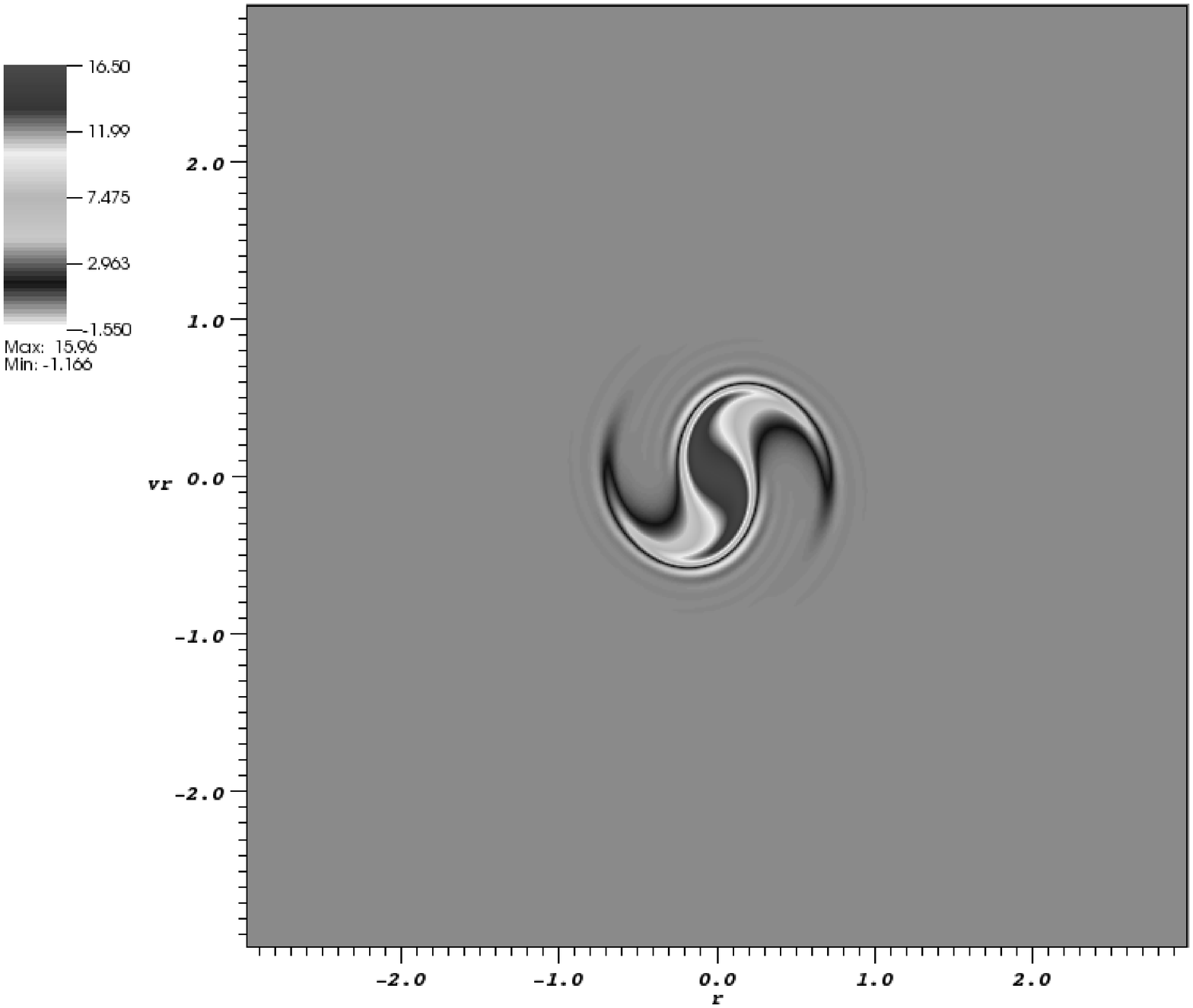} & \includegraphics[scale=0.15]{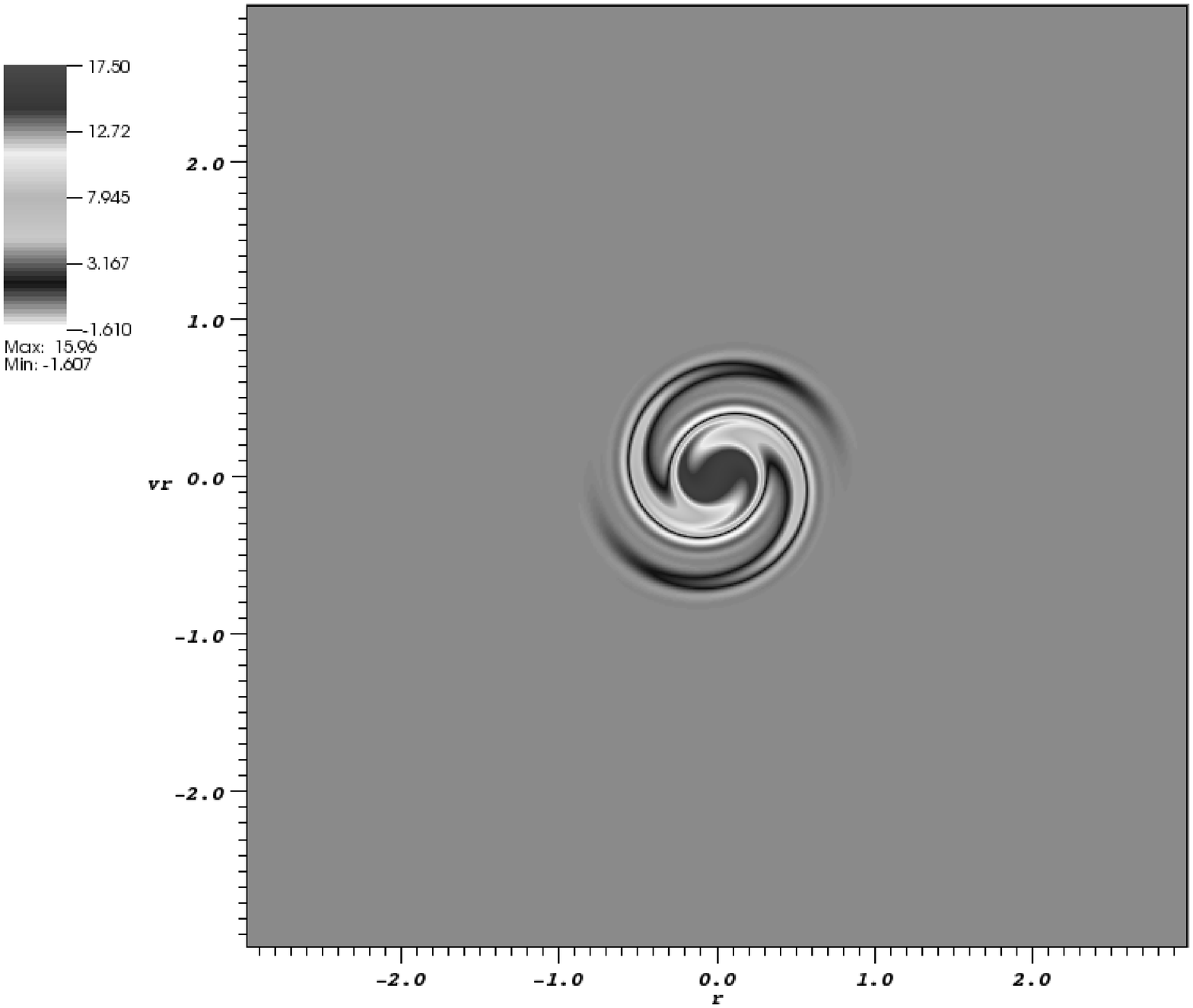} \\
\includegraphics[scale=0.15]{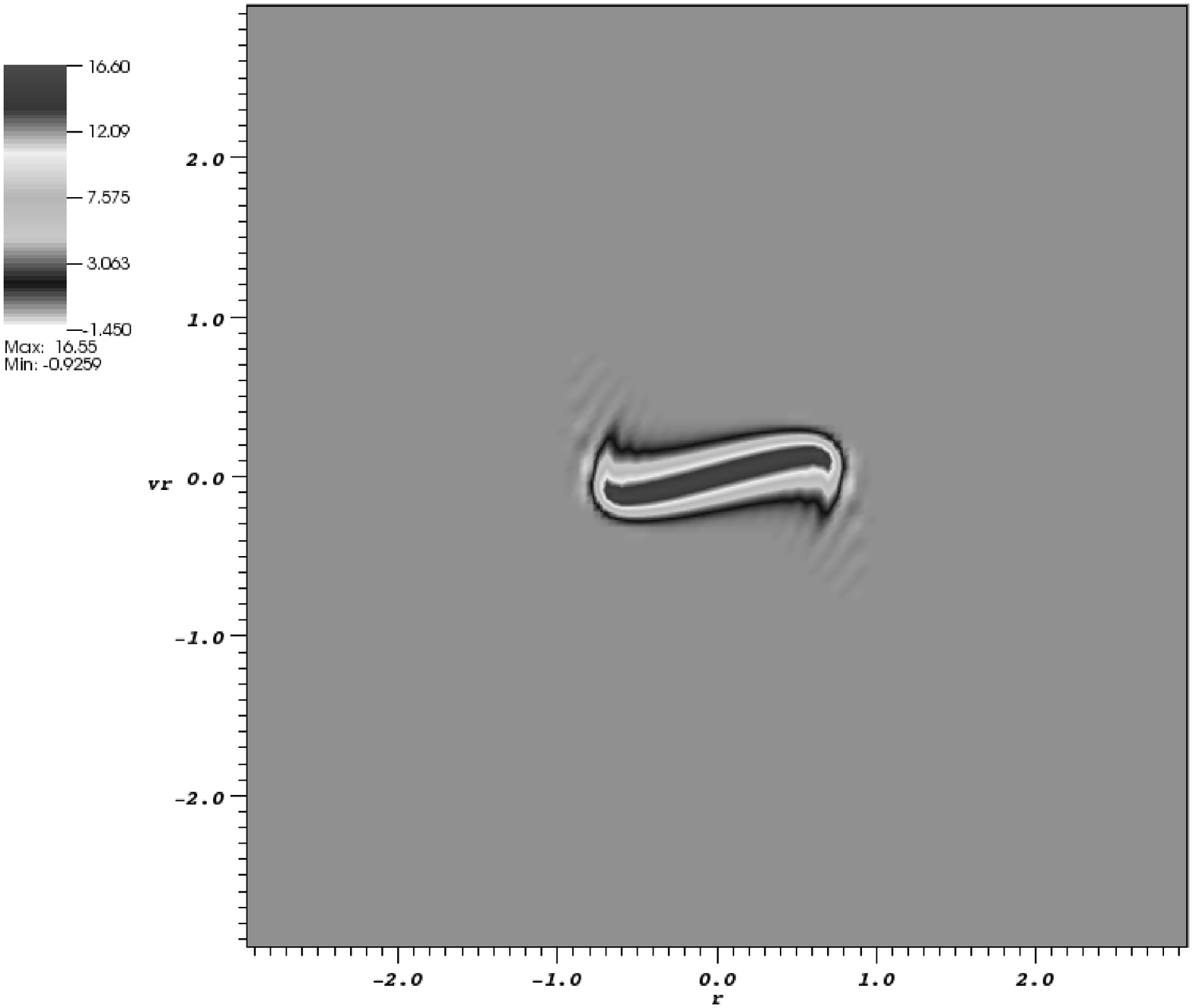} & \includegraphics[scale=0.15]{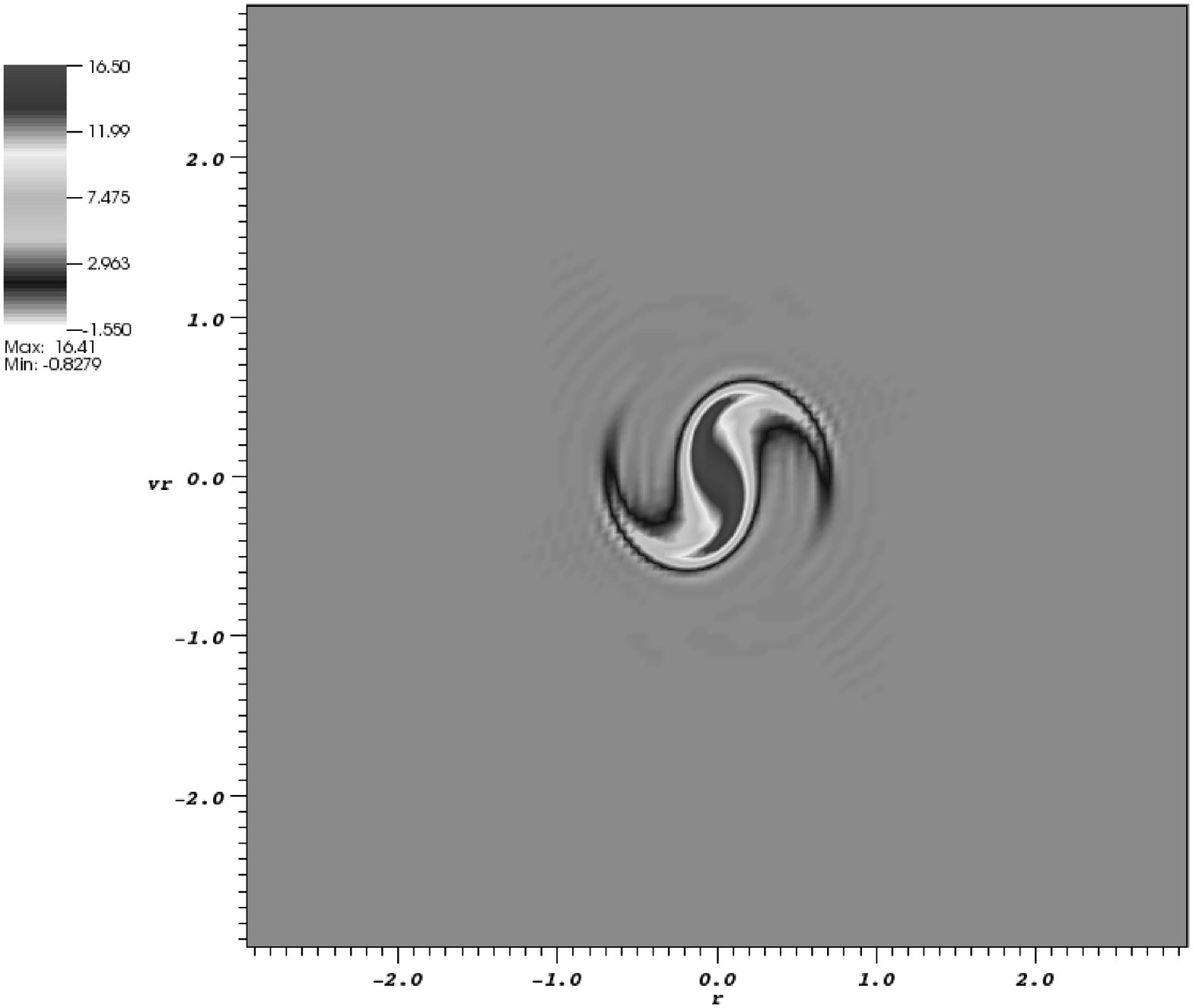} & \includegraphics[scale=0.15]{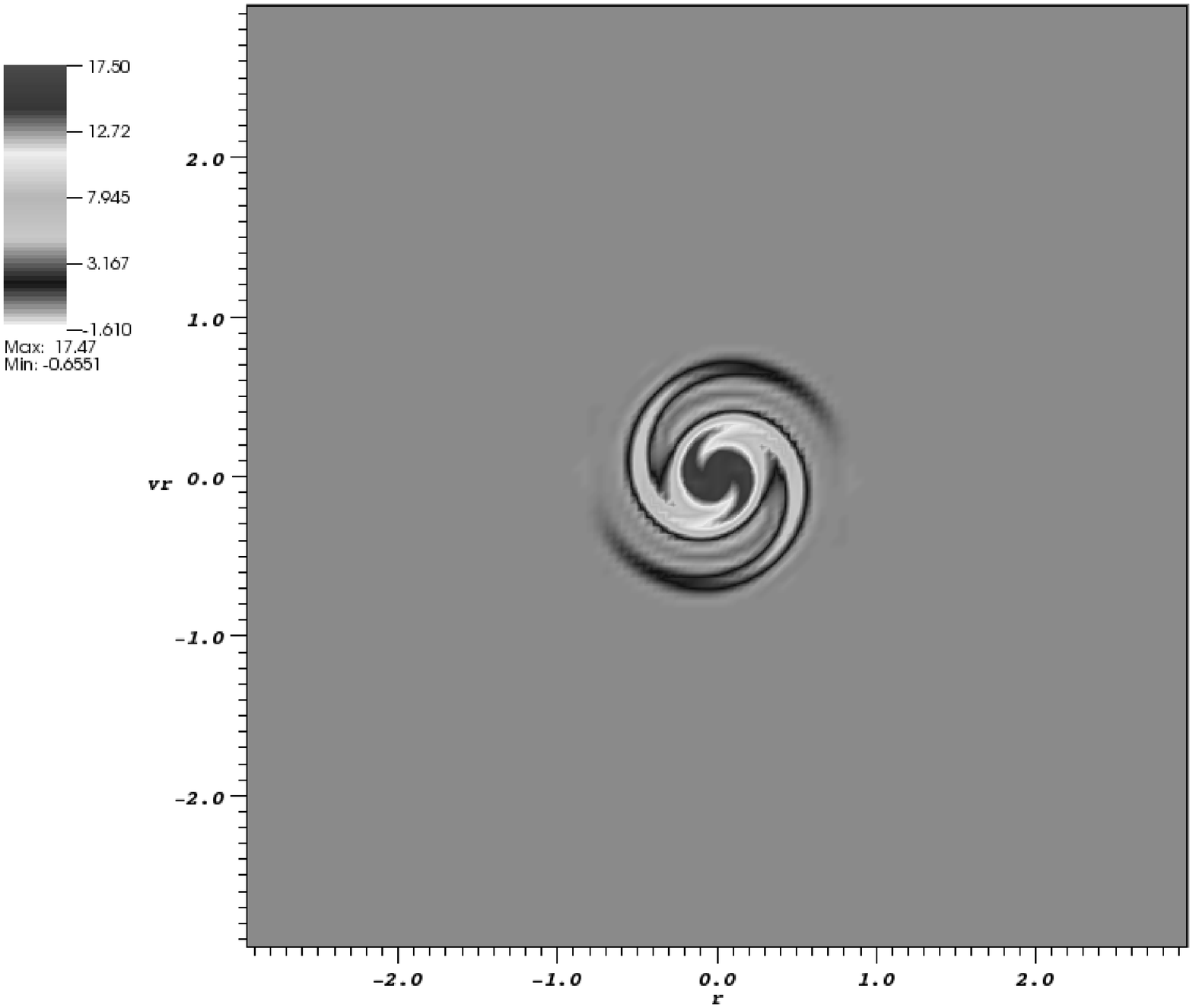} \\
\includegraphics[scale=0.15]{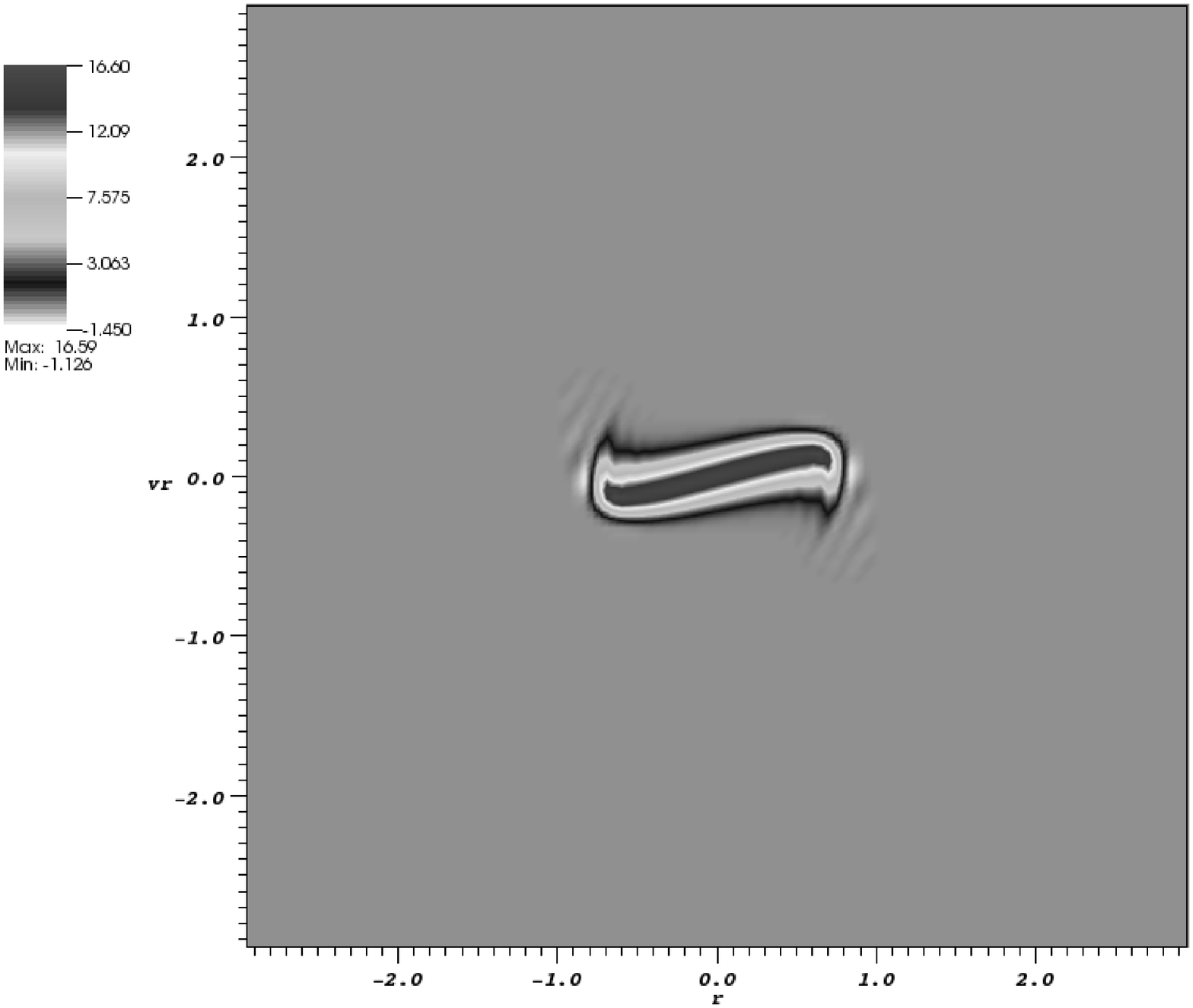} & \includegraphics[scale=0.15]{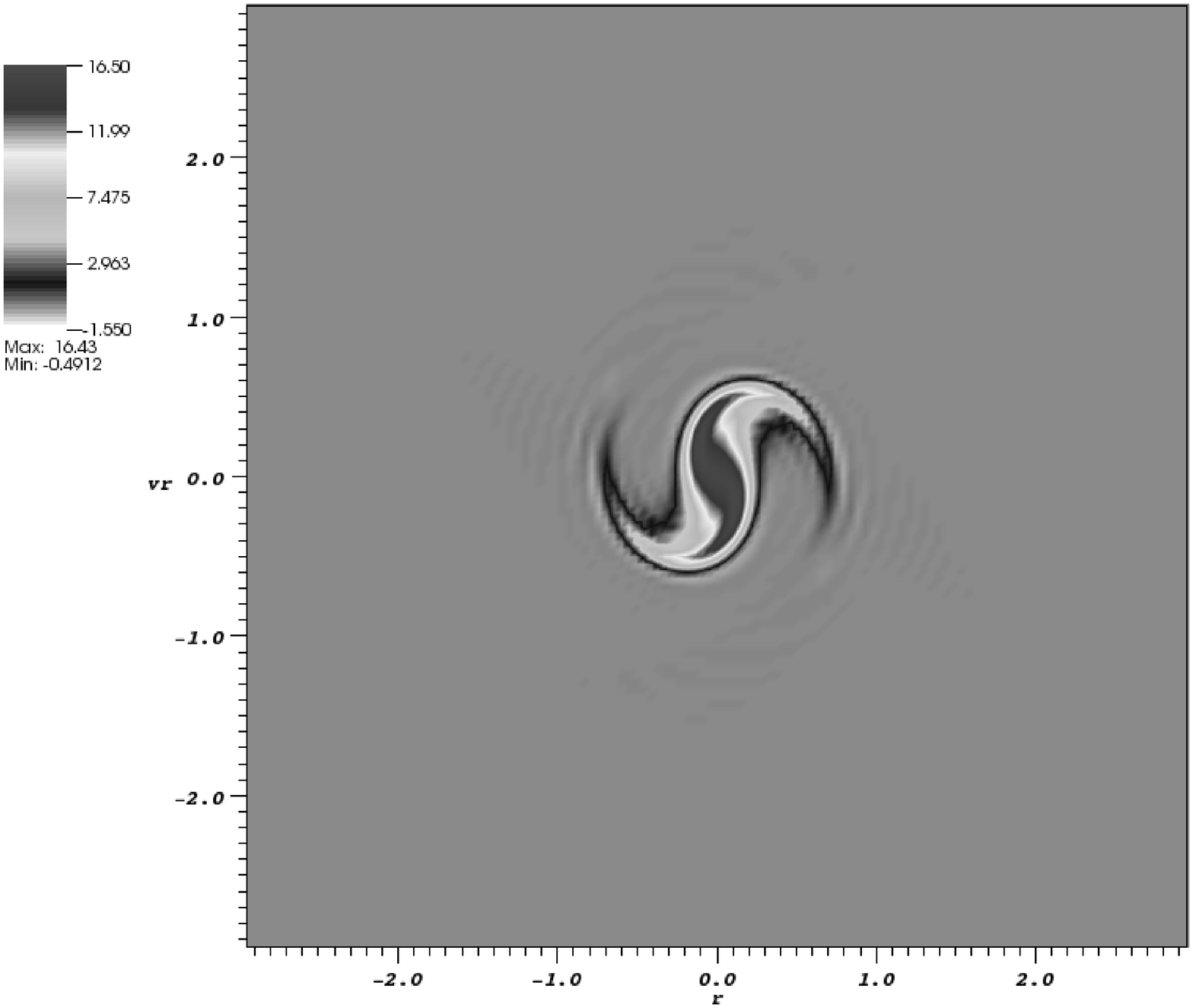} & \includegraphics[scale=0.15]{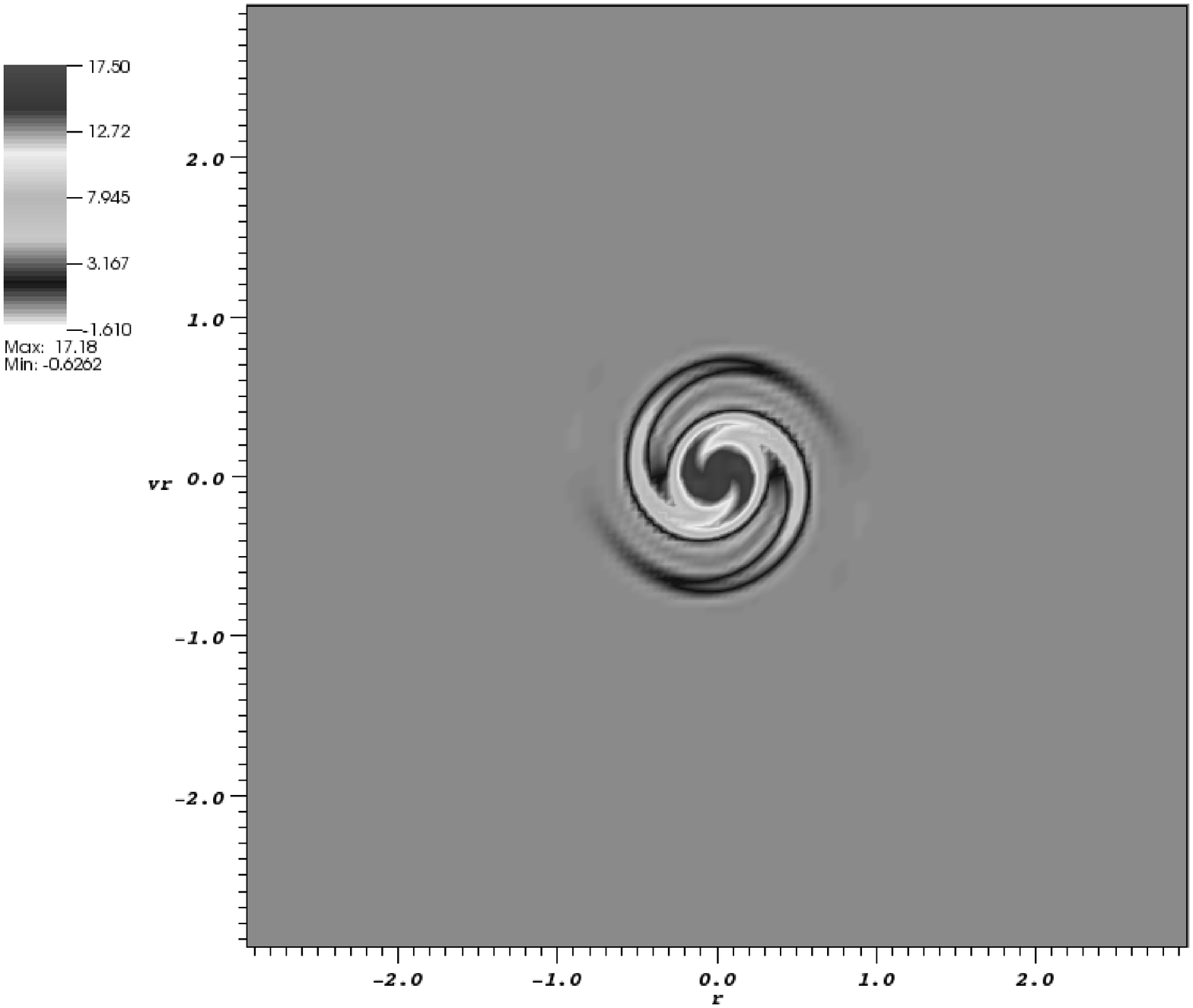} \\
$t = 1.1458$ & $t = 3.6221$ & $t = 5.8027$
\end{tabular}
\caption{Simulations of type (II') (first row), (III) (second row), and (IV) (third row) for a semi-gaussian beam with $\omega_{1} = 2$, $H_{1}(\tau) = \cos^{2}(\tau)$.}
\end{center}

\indent As we can see in Figures 8 and 9, the three simulations (II'), (III) and (IV) produce results of the same quality, even if the mesh in $r$ and $v_{r}$ directions used in the classical semi-lagrangian method is much more refined than the one used for the two-scale simulations, which enlarge to non-linear cases the conclusion we have established at the end of the previous paragraph. \\

\begin{center}
\begin{tabular}{|c|c|c|c|c|c|c|}
\hline
Case & \multicolumn{2}{|c|}{Simulation (I)} & \multicolumn{2}{|c|}{Simulation (II')} & Simulation (III) & Simulation (IV) \\
\cline{2-7} & CPU time & $N$ & CPU time & $N$ & CPU time & CPU time \\
\hline
$\omega_{1} = 4\sqrt{2}$, $H_{1} = \cos$ & 35m & 122 & 35h 6m 50s & 480 & 1h 43m 39s & 55m 3s \\
\hline
$\omega_{1} = 2$, $H_{1} = \cos^{2}$ & 37m 32s & 49 & 38h 7m 6s & 192 & 5h 45m 25s & 2h 37m 25s \\
\hline
\end{tabular}
\begin{table}[ht]
\caption{CPU time costs: the final time is $T = 6.93$ for the case where $\omega_{1} = 4\sqrt{2}$ and $H_{1} = \cos$, and $T = 6.9854$ for the case where $\omega_{1} = 2$ and $H_{1} = \cos^{2}$.}
\end{table}
\end{center}

\indent Furthermore if we observe the CPU times costs of each simulations (see Table 1), we remark that both two-scale numerical methods are much slower than the classical semi-lagrangian method when they are ran on a same mesh in $r$ and $v_{r}$ directions. It is not surprising because of the high number of fixed point problems we have to solve during the two-scale simulations. But on the other hand, if we compare the CPU time costs of the simulations (II'), (III) and (IV) which give the same quality of results, we remark that both two-scale numerical methods we have described are much faster than a precise classical semi-lagrangian method. One more time, this phenomenon can be explained by the condition (\ref{CFL-NH}) imposed to the time step $\Delta t_{NH}$ within the classical method: if we refine the mesh of the phase space, we diminish the time step, and then increase the number of iterations in time we need to reach the the final time of the simulation. As a conclusion, we can say that if we want high quality results, it is preferable in terms of CPU time cost to run one of the two-scale methods we have developed instead of the classical method. \\

\indent In a last case, we suppose that $H_{1}(\tau) = \cos^{2}(\tau)$, $\omega_{1} = 1$, $\epsilon = 10^{-2}$, and $f_{0}$ is given by (\ref{Semi-gaussian}) where $n_{0} = 4$, $r_{m} = 1.85$ and $v_{th} = 0.1$. As in the previous tests, we suppose that the time step $\Delta t_{H}$ for the two-scale methods is given by (\ref{def_Delta_t}) where $K = 1$, and that the time step $\Delta t_{NH}$ for the classical semi-lagrangian method is of the form $\frac{\Delta t_{H}}{N}$ with $N$ large enough. The goal of this numerical experiment is to observe the same structures as Fr\'enod, Salvarani and Sonnendr\"ucker have observed in \cite{PIC-two-scale}. Since the structures we want to observe are quite thin, we consider the simulation types (II''), (III') and (IV') corresponding to
\begin{itemize}
\item simulation (II''): we solve the system (\ref{NH-polar}) with a classical semi-lagrangian method, with $P_{r} = P_{v_{r}} = 256$ and $R = v_{R} = 3$,
\item simulation (III'): we solve the system (\ref{H-polar}) with a two-scale semi-lagrangian method on a two-scale mesh, with $P_{r} = P_{v_{r}} = 128$, $P_{\tau} = 20$ and $R = v_{R} = 3$,
\item simulation (IV'): we solve the system (\ref{H-polar}) with a two-scale semi-lagrangian method on a uniform mesh, with $P_{q} = P_{u_{r}} = 256$, $P_{\tau} = 20$, and $R = v_{R} = 3$, and we compute the approximation (\ref{approx_tsu}) on a uniform $257 \times 257$ grid in $(r,v_{r})$ on $[-R,R] \times [-v_{R},v_{R}]$.
\end{itemize}

\indent In the Figure 10, we can observe some results obtained with the simulations (II''), (III') and (IV'). One more time, we can remark that the two-scale methods do not need a mesh as refined as the classical method's one for producing high quality results. Furthermore, we can observe in the last column of the Figure 10 that the beam simulated with the classical semi-lagrangian method becomes unfocused in long time, even if we consider a highly refined mesh. This remark is confirmed by the Figure 11 where we notice the thin structures are better determined with a two-scale method. Moreover, the result given by the classical method is slightly out of phase and more diffusive. One more time, the main reason of this problem is the condition (\ref{CFL-NH}) imposed on the time step in the classical semi-lagragian simulation: if we consider a $513 \times 513$ grid in $(r,v_{r})$, this condition induces a so small value of $\Delta t_{NH}$ (nearby $3.84 \times 10^{-5}$) that the numerical diffusion introduced by the interpolations makes the beam unfocused in a long time simulation. On the other hand, the two-scale results match well with the expected long time behavior as described in \cite{PIC-two-scale}. 

\begin{center}
\begin{tabular}{ccc}
\includegraphics[scale=0.15]{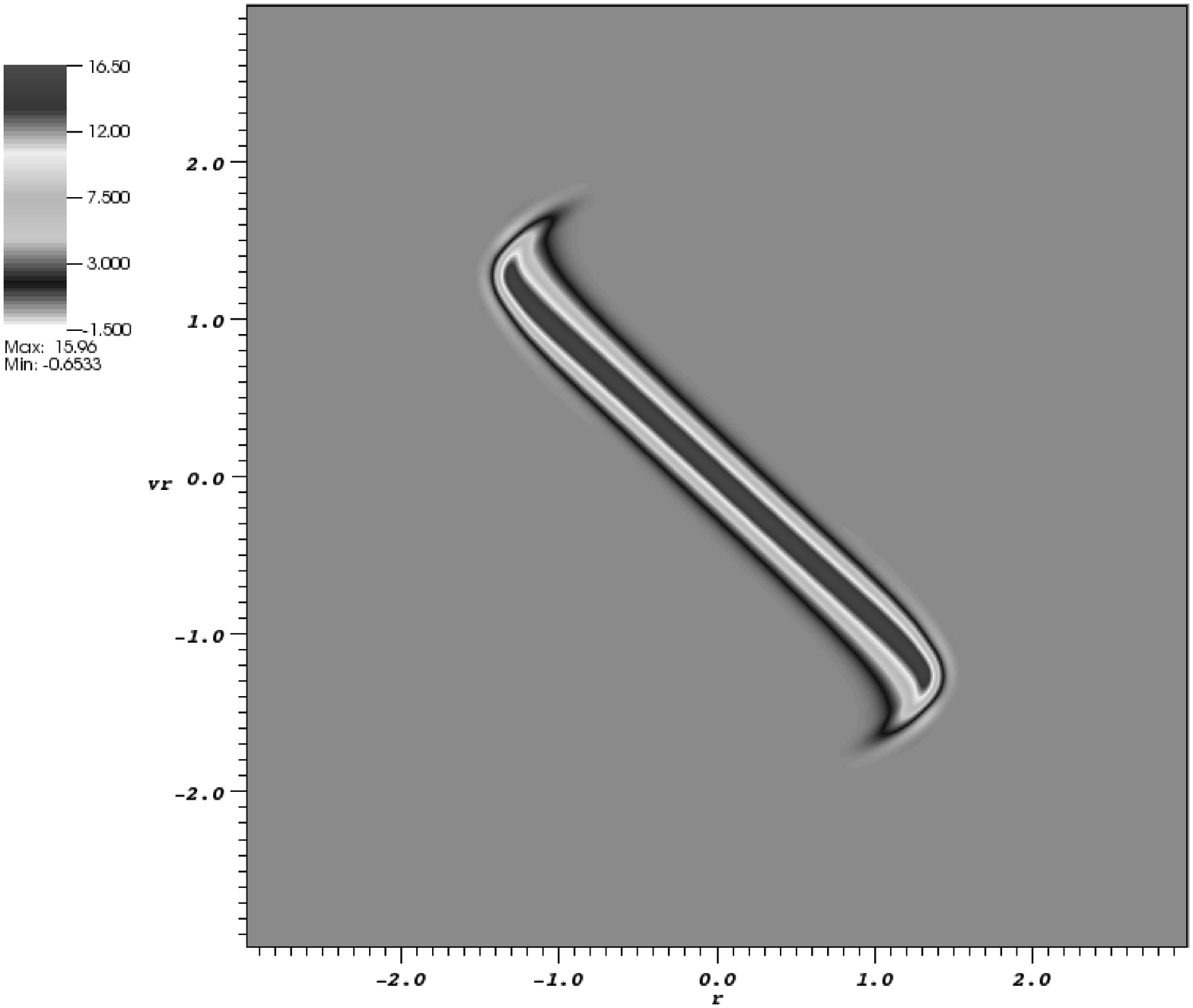} & \includegraphics[scale=0.15]{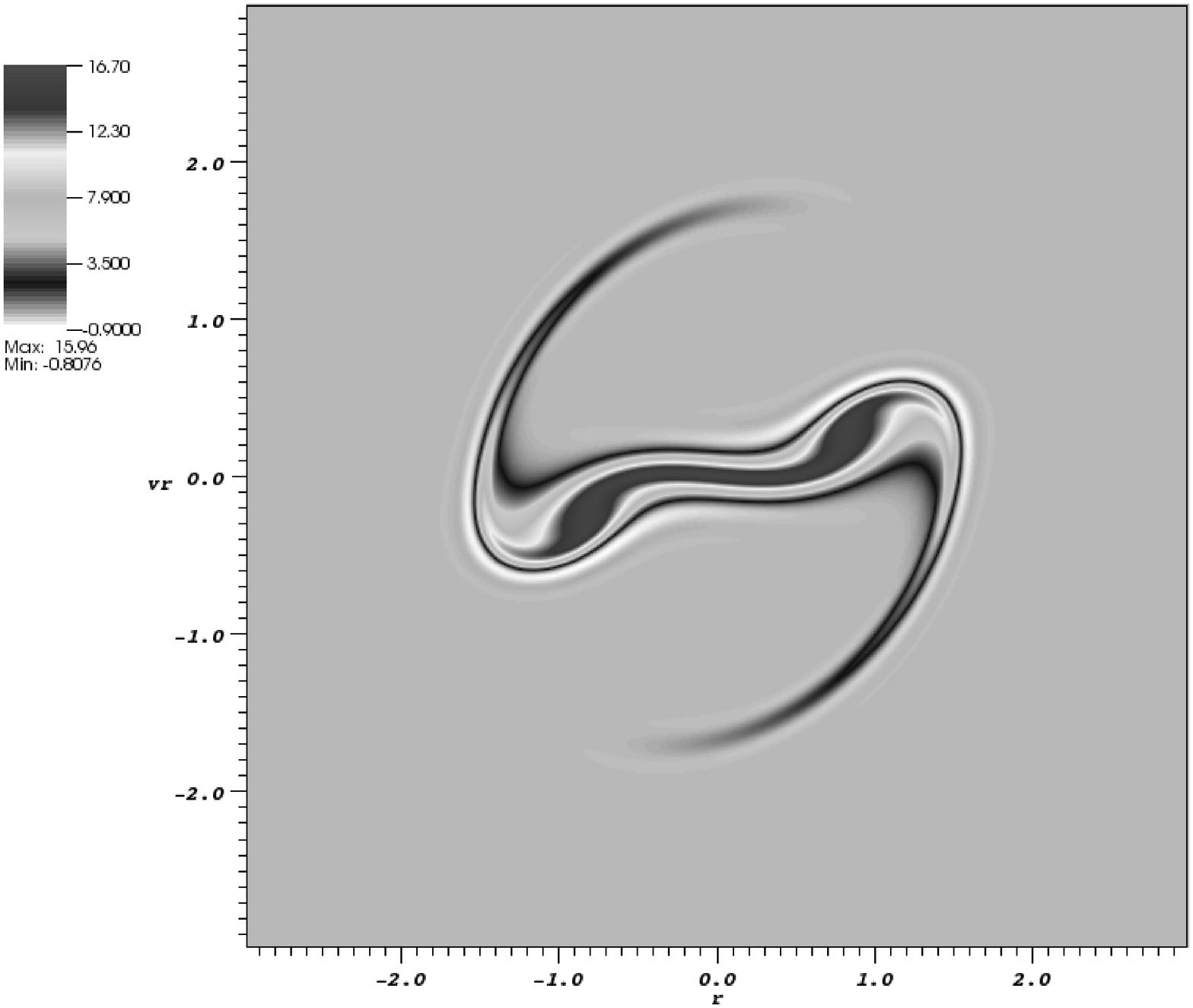} & \includegraphics[scale=0.15]{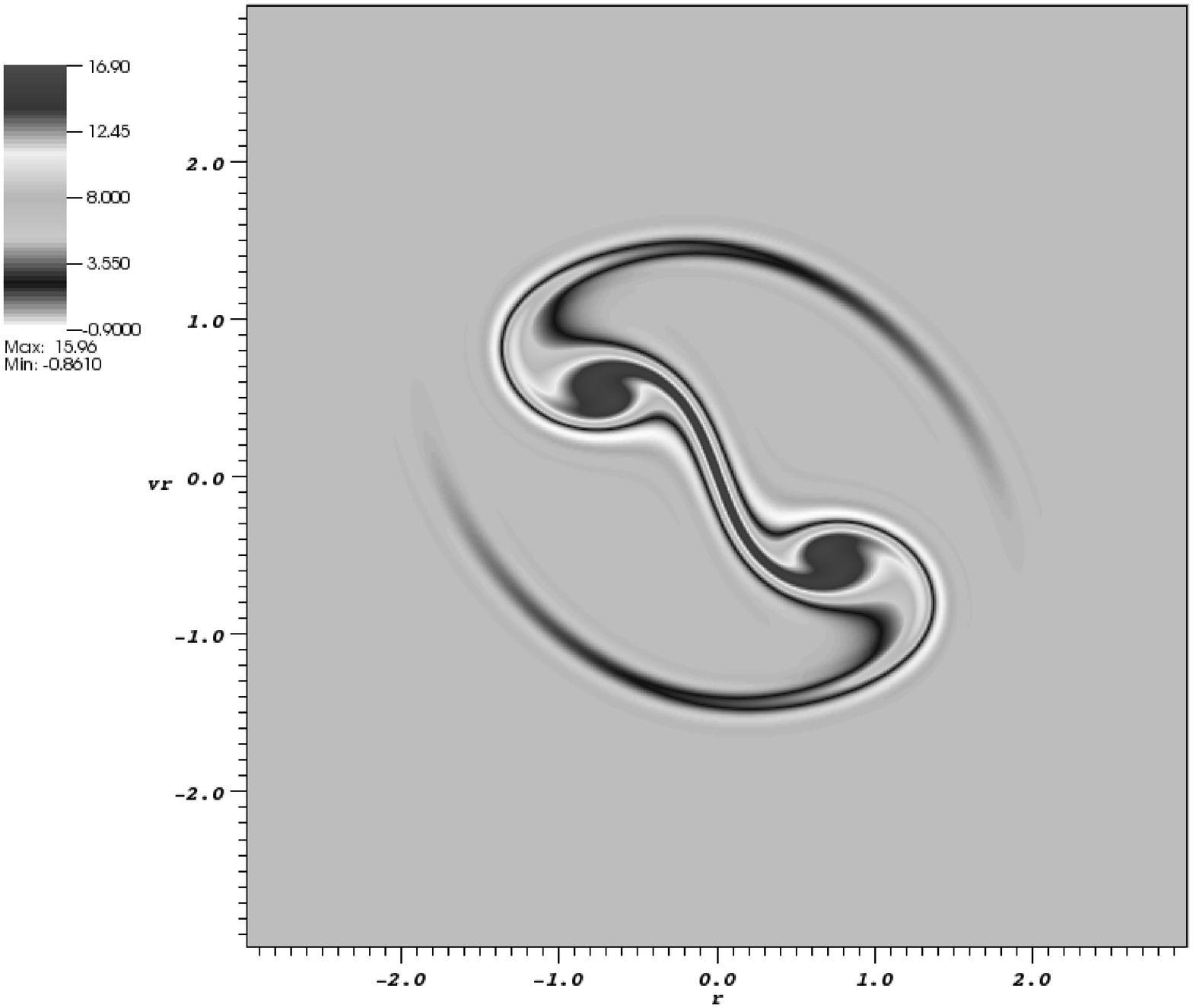} \\
\includegraphics[scale=0.15]{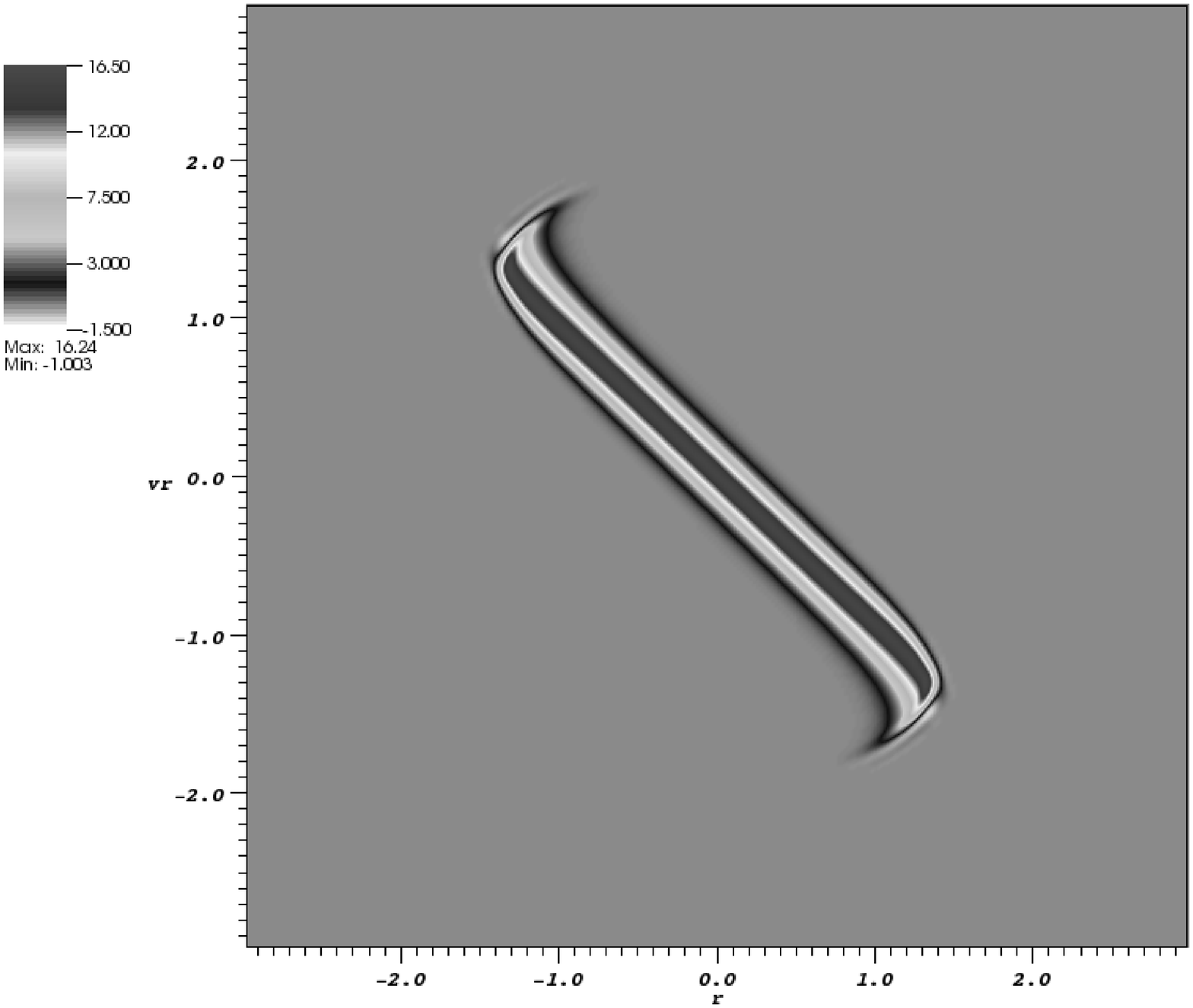} & \includegraphics[scale=0.15]{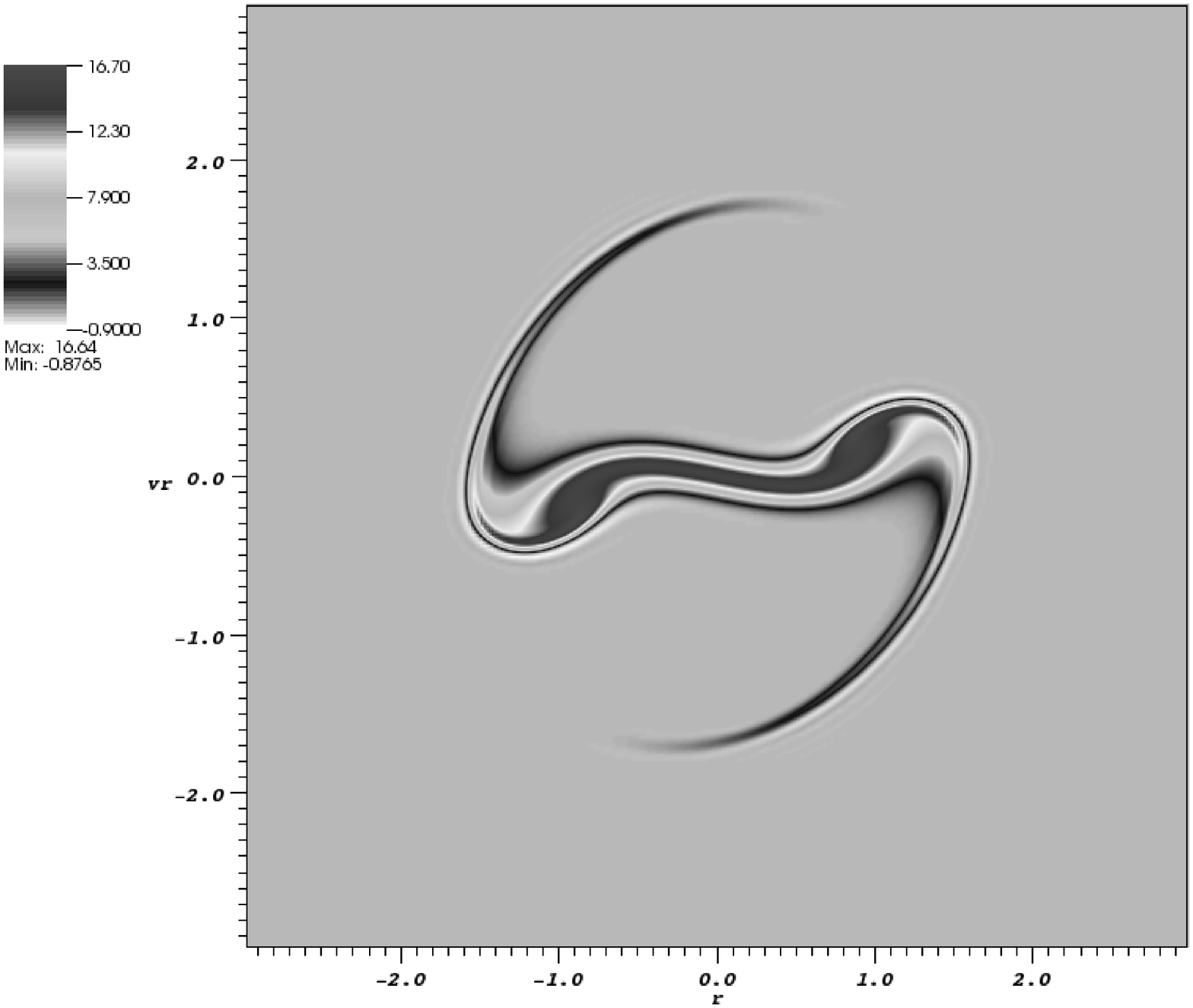} & \includegraphics[scale=0.15]{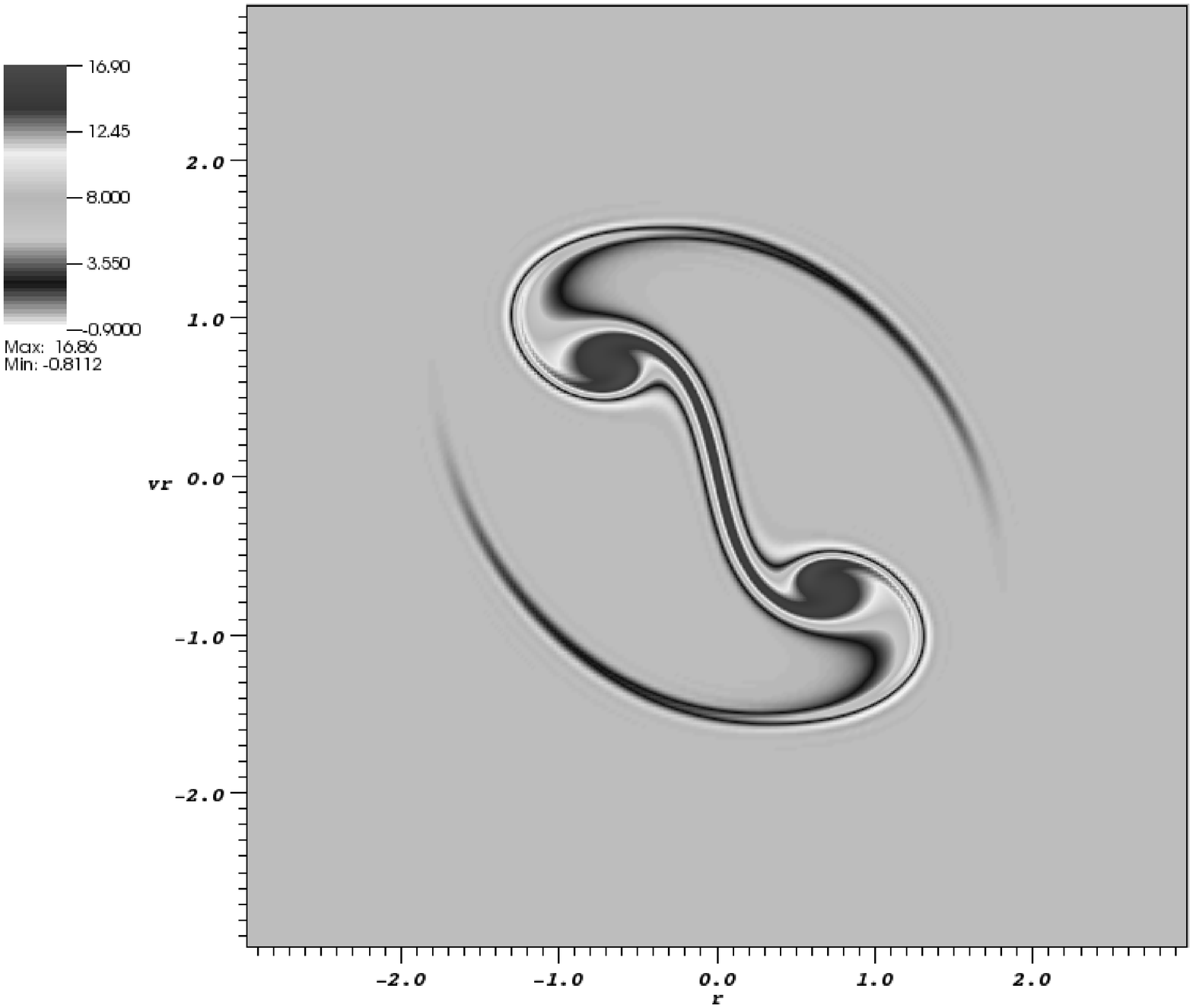} \\
\includegraphics[scale=0.15]{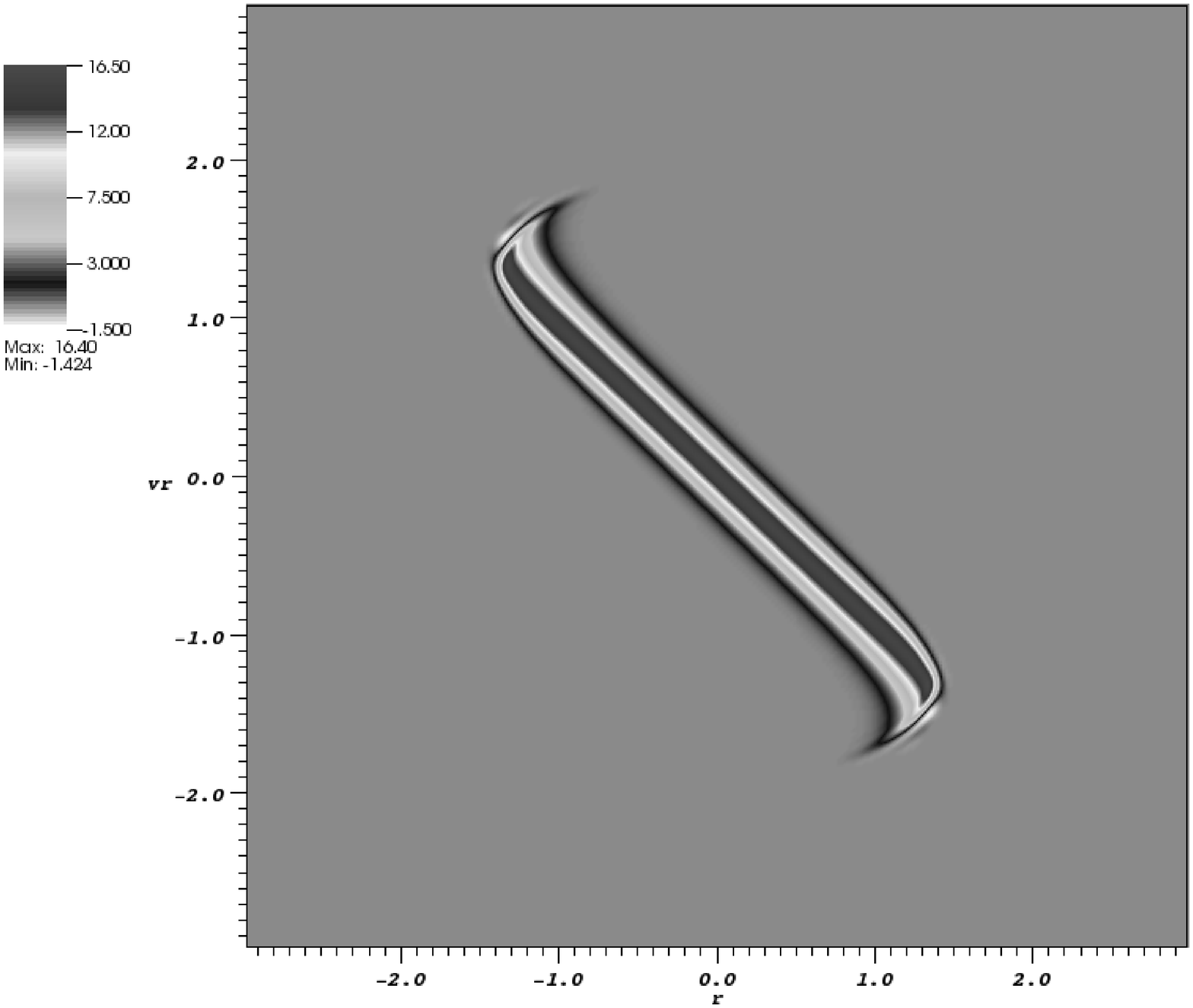} & \includegraphics[scale=0.15]{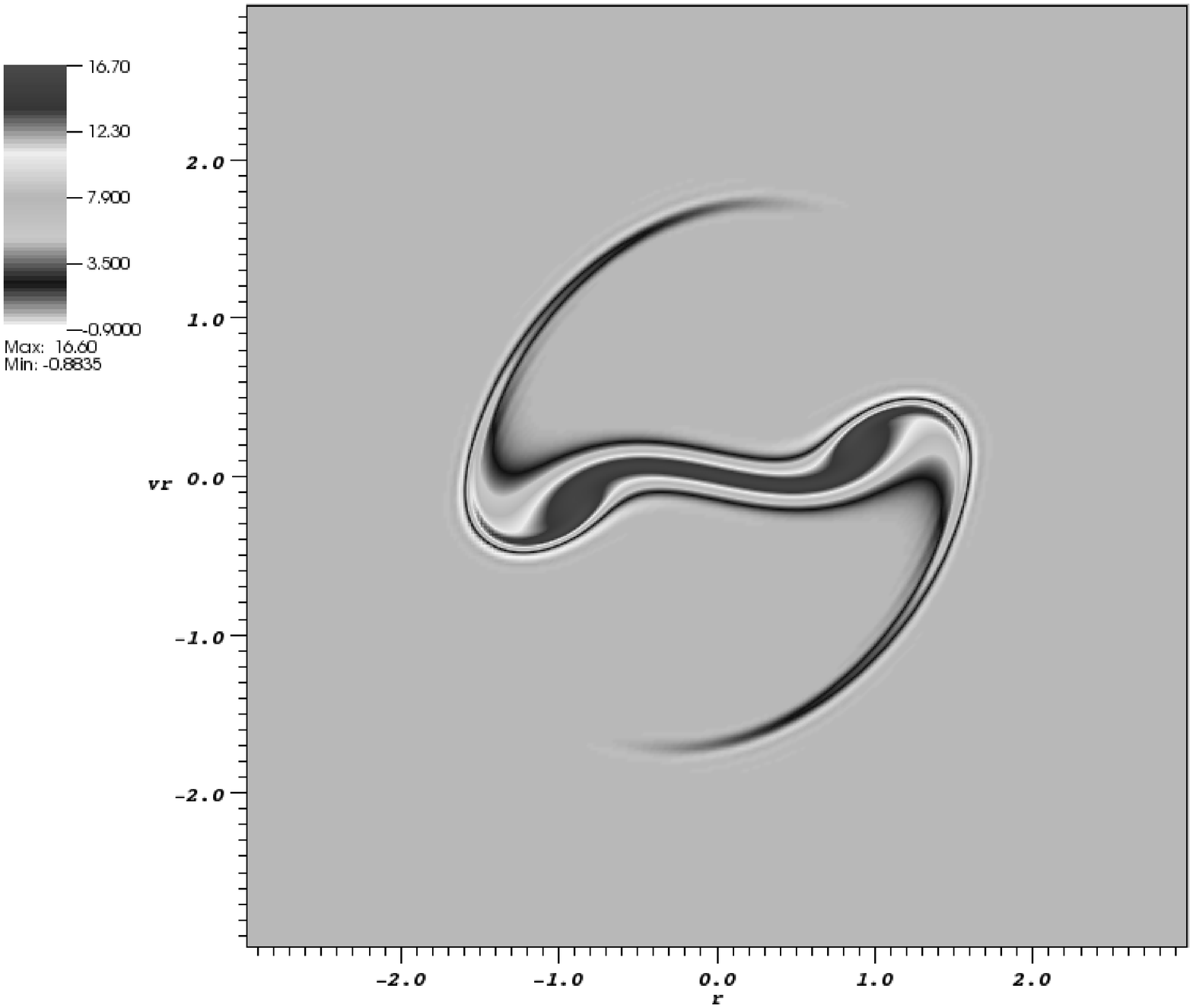} & \includegraphics[scale=0.15]{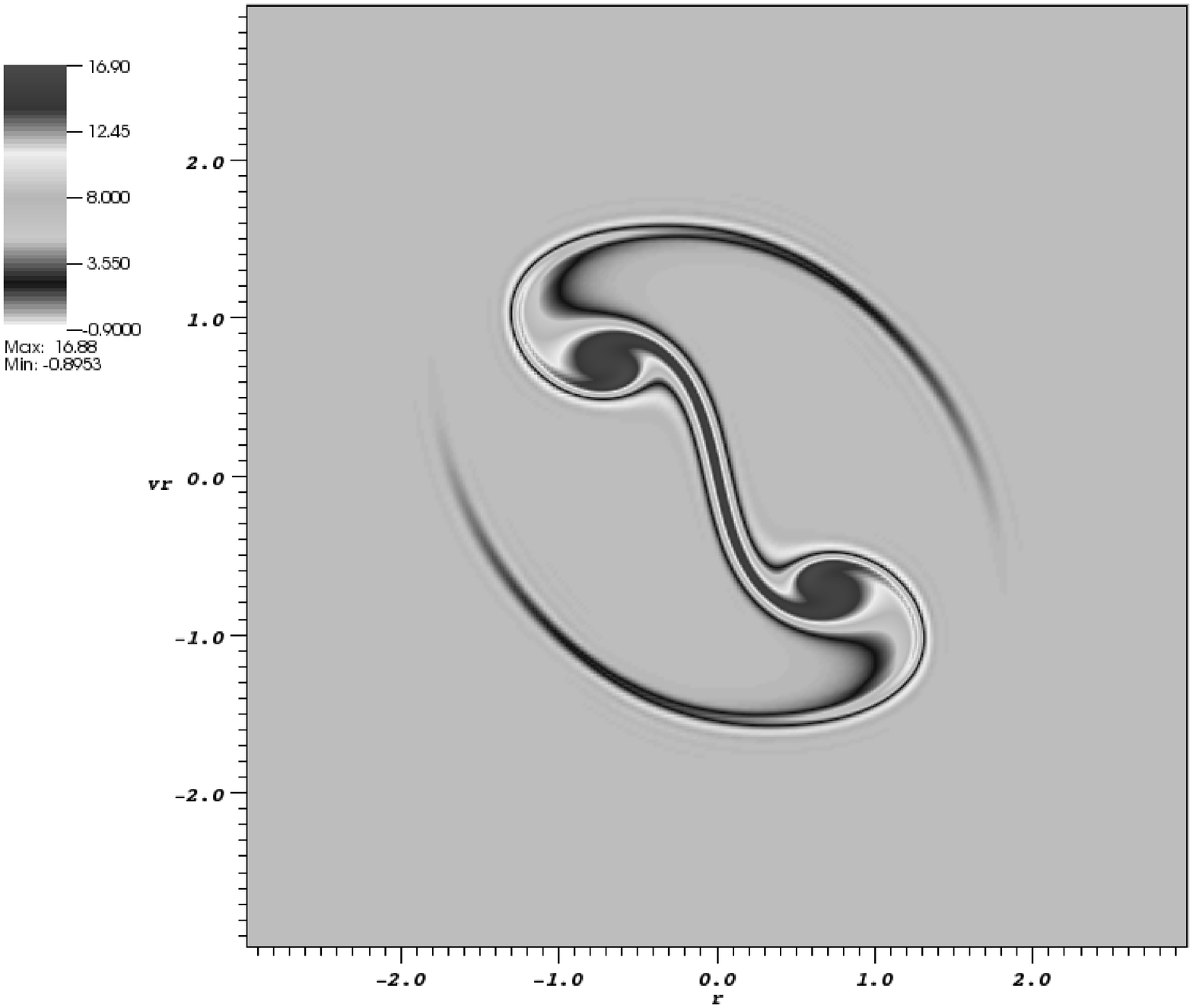} \\
$t = 1.3464$ & $t = 4.3388$ & $t = 5.1462$
\end{tabular}
\caption{Simulations of type (II'') (first row), (III') (second row) and (IV') (third row) for a semi-gaussian beam with $\omega_{1} = 1$, $H_{1}(\tau) = \cos^{2}(\tau)$ and $r_{m} = 1.85$.}
\end{center}

\begin{center}
\begin{tabular}{cc}
\includegraphics[scale=0.22]{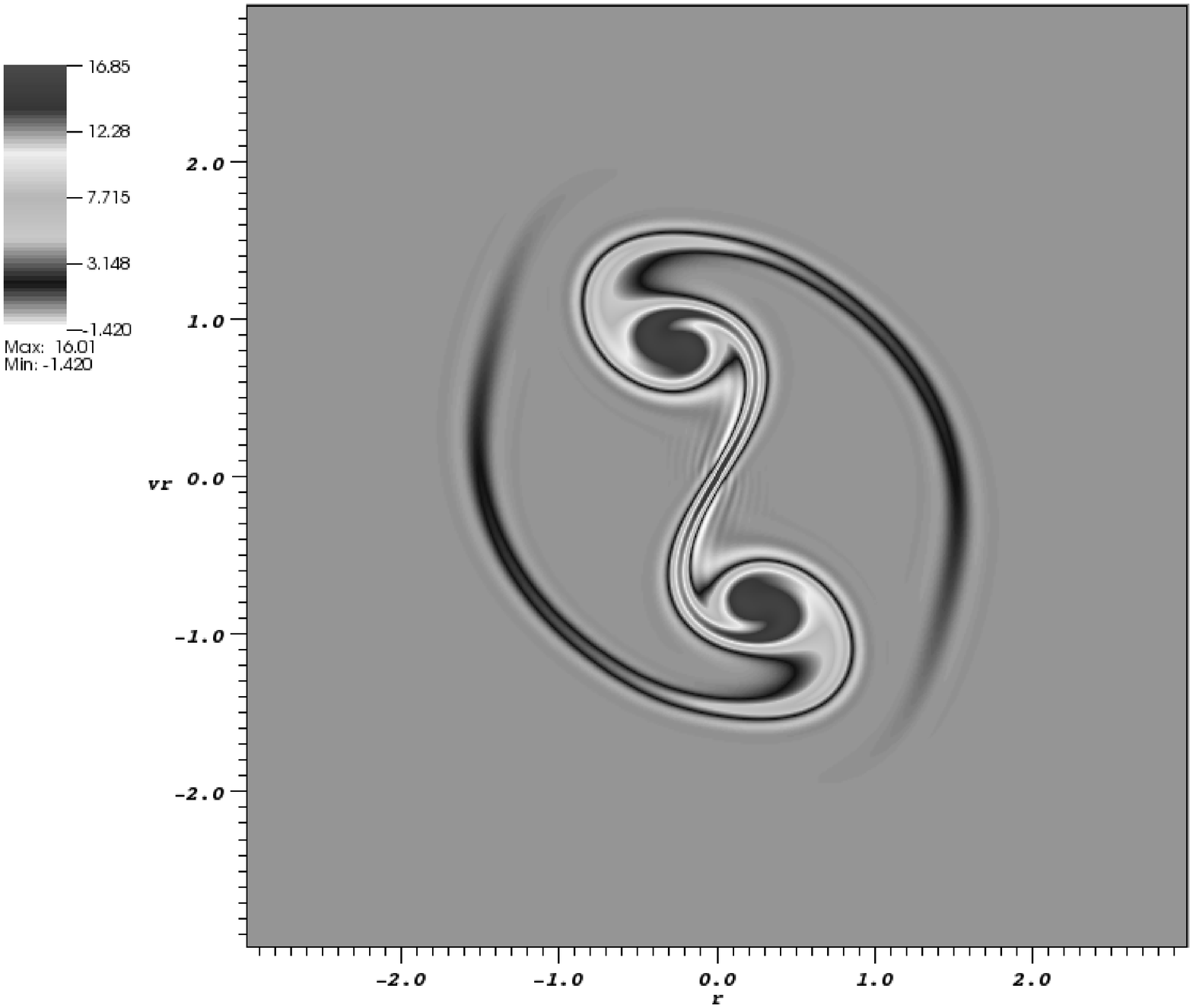} & \includegraphics[scale=0.22]{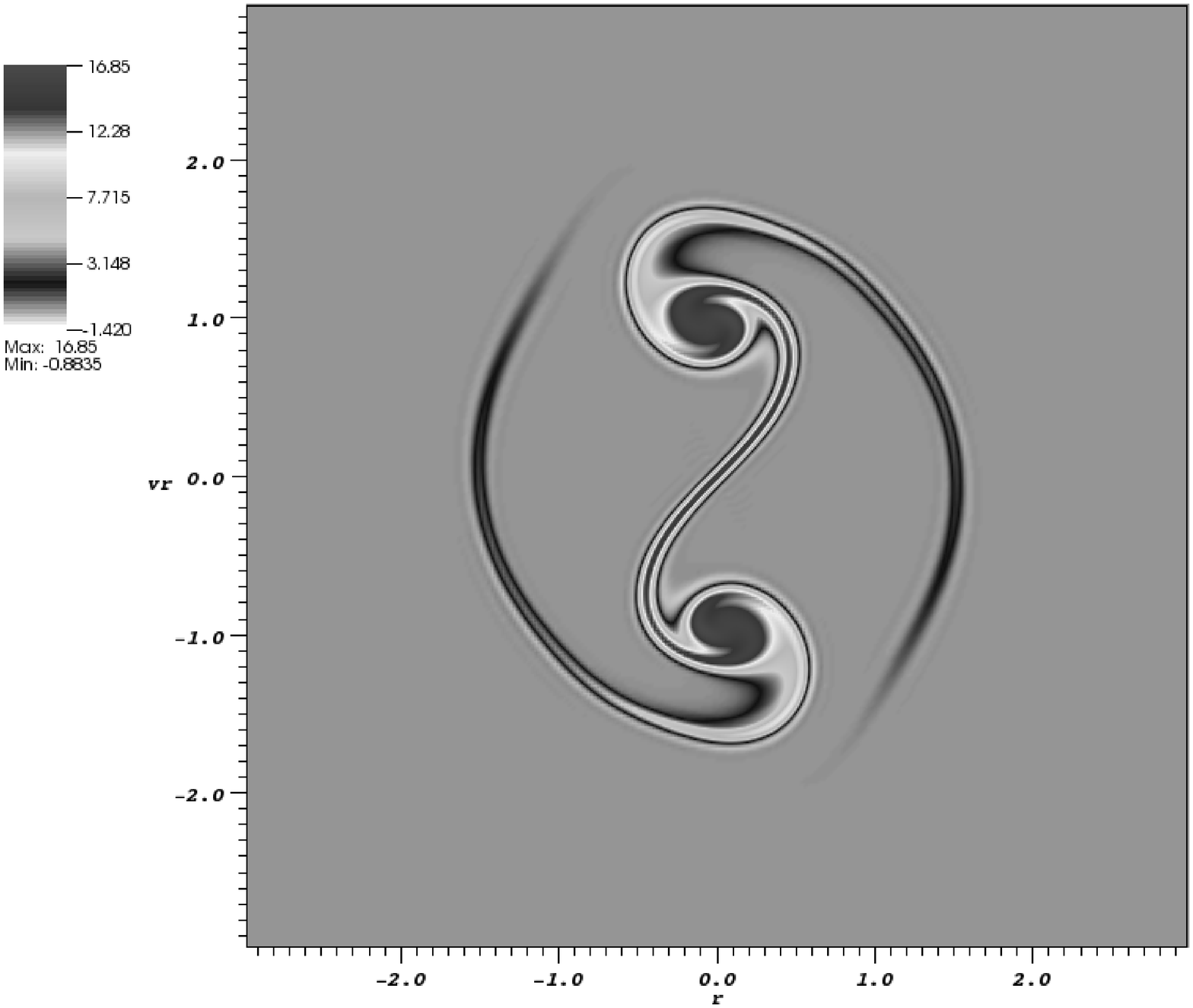} \\
(II'') & (III')
\end{tabular}
\caption{Simulations of type (II'') and (III') for a semi-gaussian beam at time $t = 5.984$ with $\omega_{1} = 1$, $H_{1}(\tau) = \cos^{2}(\tau)$ and $r_{m} = 1.85$.}
\end{center}

\section{Conclusion and perspectives}

\indent We have built a two-scale semi-lagrangian method and proposed a new mesh in order to simplify the computation of the electric field and, which leads to another two-scale semi-lagrangian method. These methods have been tested on non-smooth initial conditions, especially semi-gaussian beam initial conditions: on linear cases, we have concluded that these two-scale methods are very efficient, even if we consider a coarse mesh in $r$ and $v_{r}$ directions, contrary to a classical semi-lagrangian scheme on the model (\ref{NH-polar}) which needs much more points to produce good results. On non-linear cases, we have reached the same conclusion for short time simulation, not only in terms of quality of results, by also in terms of CPU time costs. These results are very promising for extensions to higher dimensional problems such as the two-scale limit models obtained by Frénod and Sonnendrücker in \cite{Homogenization_VP,Long_time,Larmor-radius}, the finite Larmor radius approximation obtained in \cite{Bostan_2007}, or other charged particle beam problems which cannot allow any time splitting. Furthermore, for long time simulation, both two-scale numerical methods we have developed give very good results, contrary to the classical semi-lagrangian method which is so penalized by its numerical diffusion for long time simulations that it produces results which do not not correspond to the expected behavior. These results are also promising since they consolidate the conclusions of Frénod, Salvarani and Sonnendr\"ucker in \cite{PIC-two-scale}, which are that a two-scale numerical method can be successfully used in a context of non-smooth initial data. \\
\indent Finally, we also have remarked that, even if they are faster than the classical semi-lagrangian method for obtaining the same quality of results, both two-scale numerical methods need a very high CPU time cost. Since this time is essentially spent for computing the fixed point problems within the methods, it can be interesting to find a way to improve this part of the methods.

\paragraph{Acknowledgements} The author wishes to thank E. Fr\'enod and E. Sonnendr\"ucker for the discussions on the topic.

\textit{ \\ }

E-mail adress: \verb|mouton@math.u-strasbg.fr|

\end{document}